%% file: root.tex
\documentstyle[12pt]{report} 

\newtheorem{lemma}{Lemma}
\def\blm{\begin{lemma}}
\def\elm{\end{lemma}}
\newtheorem{propos}[lemma]{Proposition}
\newtheorem{theorem}{Theorem}
\newtheorem{cor}{Corollary}
\newtheorem{defi}{Definition}
\def\beq {\begin{equation}}
\def\eeq {\end{equation}}
\def\btab{\begin{tabbing}}
\def\etab{\end{tabbing}}
\def\barr{\begin{array}}
\def\ear{\end{array}}

\def\bea{\begin{eqnarray}}
\def\eea{\end{eqnarray}}

\def\bpp{\begin{propos}}
\def\epp{\end{propos}}
\def\btm{\begin{theorem}}
\def\etm{\end{theorem}}
\def\ben{\begin{enumerate}}
\def\een{\end{enumerate}}
\def \H {{\cal H}}
\def \f {{\cal F}}
\def \fo {{\cal F}^0}
\def \fl {{\cal F}^1}
\def \fla{{\cal F}^1_{\alpha}}
\def \flb{{\cal F}^1_{\beta}}
\def \li {{\cal L}}

\def\s{\bigskip}
\def \ep{\hfill{$\Box$}}

\def\isto {\widetilde{\longrightarrow}}
\def\bdot {\hbox{\raise .4ex\hbox{\large\bf .}}}
\def\ub{\underline}

\def\lra {\longrightarrow}

\def\id{{1\mkern-5mu {\rm I}}}

\def\cc {c\mkern -1.2mu \ell}

\def\diff
{{\cal
D}i\!f\!f}
\def\h {{\tt g}\,}

\def \hh {{\sf h}}

\mathchardef\mouth="315E
\mathchardef\ffc= "334B
\mathchardef\nose="333E

\def\smile{\raisebox{2.1ex}{$\ffc$}\mkern
-20 mu\raisebox{ -.35 ex}{$\mouth$}\mkern -15 mu \raisebox{1 ex}{$.\,.$}
\mkern -13.5 mu \raisebox {.95ex}{$\nose$}\quad}

\mathchardef\bbigcap="335C
\mathchardef\lleg="3237

\def\ccap {\raisebox {2ex}{$\bbigcap$}\mkern-14.7mu\raisebox
{0.1ex}{$\bigcap$}\mkern -12.7mu \raisebox {-.85ex}{$\lleg\mkern 9.5mu \lleg$}
\,}

\def \frL{\Diamond\mkern -12.5mu
\raisebox {.45ex}{\mbox{\small-\tiny 1}}\,}
\def \FLO{\Diamond\mkern -10mu
\raisebox {.45ex}{\mbox{\tiny 1}}\,}
\def \FLN{\Diamond\mkern -10mu
\raisebox {.45ex}{\mbox{\tiny n}}\,}

\def \ct {Cob_3}
\def \wct {\widetilde {Cob_3}}
\def \dfll {\leaders \hbox to 1em {\hss.\hss}\hfill}

\def\twoheadrightarrow{\to\mkern-20mu\to}
\def\thrafill{$\mathsurround=0pt \mathord- \mkern-6mu \cleaders\hbox{$\mkern-2mu
\mathord- \mkern-2mu$}\hfill \mkern-6mu\mathord\twoheadrightarrow$}

\def\clr { \raisebox {0ex} {\mbox{$\Box \mkern -3.4mu
\raisebox {.3ex}{=} \mkern -3mu \Box$}}} 

\def \ub {\underline}

\input{epsf.sty}

\begin{document}

\include{BLcont}

\include{BL0}

\setcounter{chapter}{1}
\include{BL1}

\setcounter{chapter}{2}
\include{BL2}
\setcounter{chapter}{3}
\include{BL3}

\setcounter{chapter}{4}
\include{BL4}

\include{BLbibl}

\end{document}

%% file: BLcont.tex
\ 
\hfill{{\small To appear in } {\em Adv. Math.}}

\vspace*{.5cm}

\begin{center}

{\LARGE\bf 
Bridged Links and \\ Tangle Presentations of Cobordism Categories
}

\vspace*{1cm}
\bigskip

{\Large Thomas Kerler\\
}

\vspace*{1.5cm}

First Version July 1994 (This  Version May 1997)\\

\end{center}
\vspace*{3.5cm}

\noindent{\bf Abstract :}{\small 
\ \ We develop a calculus of surgery data, called {\em bridged links}, which involves
besides  links also pairs of balls that describe one-handle attachements.

As opposed to the usual link calculi of Kirby and others this
description uses only elementary, {\em local} moves(namely modifications and
isolated cancellations ), and it is valid also on non-simply connected and
disconnected manifolds. In particular, it allows us to give a 
presentation of a 3-manifold by doing surgery on any other 3-manifold with the
same boundary. 

Bridged link presentations on unions of handlebodies are used to give a
Cerf-theoretical derivation  
of presentations of 2+1-dimensional cobordisms categories in terms of 
planar ribbon tangles and their composition rules. As an application we
give a different, more natural proof of the Matveev-Polyak presentations of
the mapping class group, and, furthermore, find systematically surgery 
presentations of general mapping tori.

We discuss a natural  extension of the Reshetikhin Turaev invariant to
the calculus of bridged links. Invariance follows now - similar as for knot 
invariants -  from simple identifications of the elementary moves with 
elementary categorial relations for invariances or cointegrals, respectively. 
Hence, we avoid the lengthy computations and the unnatural Fenn-Rourke
reduction of the original proofs. Moreover, we are able to start from a much weaker
``modularity''-condition, which implies the one of Turaev.

Generalizations of the presentation to cobordisms of surfaces with
boundaries are outlined.}

\newpage

\hphantom{xxxx}

\pagenumbering{roman}\

\vspace*{-1.7cm}

\section*{Contents}

\paragraph{0) Introduction and Survey of Results} \hfill \pageref{pg-0}

\paragraph{1) Framed Cobordisms}\hfill \pageref{pg-1}

\subparagraph{1.1) Construction of Bounded  Cobordisms }\hfill \pageref{pg-11}

{\em 1.1.1) 2+1 Cobordisms}\dfll\pageref{pg-111}

{\em 1.1.2) Bounding 3+1 Cobordisms}\dfll\pageref{pg-112}

{\em 1.1.3) The Category $\wct$ }\dfll\pageref{pg-113}

{\em 1.1.4) Moves between Four Manifolds}\dfll\pageref{pg-114}

\subparagraph{1.2) A Central Extension by $\Omega_4$ }\hfill\pageref{pg-12}

{\em 1.2.1) The $\Omega_4$-Action and Anomalous TQFT's}\dfll\pageref{pg-121}

{\em 1.2.2) Another Two-Cocycle}\dfll\pageref{pg-122}

\subparagraph{1.3) Ordering and Connectivity}\hfill\pageref{pg-13}

\paragraph{2) Enhanced Surgery Presentations}\hfill\pageref{pg-2}

\subparagraph{2.1) Some Cerf Theory}\hfill\pageref{pg-21}

{ \em 2.1.1) Excellent and Codimension One Functions}\dfll\pageref{pg-211}

{\em 2.1.2) Paths and Deformations of Paths of Functions}\dfll\pageref{pg-212}

\subparagraph{2.2) Bridged Link Presentations}\hfill\pageref{pg-22}

{\em 2.2.1) Surgery Presentation from Bridged Links}\dfll\pageref{pg-221}

{\em 2.2.2) Moves for Bridged Links}\dfll\pageref{pg-222}

\subparagraph{2.3) Equivalences of Bridged Links}\hfill\pageref{pg-23}

{\em 2.3.1) Completeness of Moves for Bridged Links}\dfll\pageref{pg-231}

{\em 2.3.2) Reduction of Moves}\dfll\pageref{pg-232}

\subparagraph{2.4) Bridged Ribbon Graphs}\hfill\pageref{pg-24}


\paragraph{3) Standard and Tangle Presentations of Cobordisms}\hfill\pageref{pg-3}

\subparagraph{3.1) Standard Presentations on $S^3$}\hfill\pageref{pg-31}

{\em 3.1.1) Standard Presentation of $\breve M\,$ and the class $\cal U$}\dfll\pageref{pg-311}

{\em 3.1.2) Standard Links in $H^{\pm}_g\,$}\dfll\pageref{pg-312}

{\em 3.1.3) The $\sigma$-Move and the Lemma of Connecting Annuli}\dfll\pageref{pg-313}

{\em 3.1.4) Existence of Standard Presentations and a Projection on $\,\cal S$}\dfll\pageref{pg-314}

{\em 3.1.5) Moves in a Standard Presentation}\dfll\pageref{pg-315}

\subparagraph{3.2) Tangle Presentation of Cobordisms}\hfill\pageref{pg-32}

{\em 3.2.1) From Standard Presentations to Admissible Tangles}\dfll\pageref{pg-321}

{\em 3.2.2) Moves for Admissible Tangles}\dfll\pageref{pg-322}

{\em 3.2.3) Compositions of Admissible Tangles}\dfll\pageref{pg-323}

{\em 3.2.4) Na\"\i ve Compositions, Connecting Annuli, and Closed Tangles}\dfll\pageref{pg-324}

\paragraph{4) Applications and Implications}\hfill\pageref{pg-4}

\subparagraph{4.1) Invertible Cobordisms and Presentations of Mapping Tori}\hfill\pageref{pg-41}

{\em 4.1.1) Tangle Presentation of the Mapping Class Groups}\dfll\pageref{pg-411}

{\em 4.1.2) Presentations of Manifolds from ${\cal
T}=\pi_0\bigl(\diff(\Sigma)\bigr)\,$}\dfll\pageref{pg-412} 

\subparagraph{4.2) On the Reshetikhin Turaev Invariant}\hfill\pageref{pg-42}

{\em 4.2.1) Invariants  of Three Manifolds from Bridged Ribbons}\dfll\pageref{pg-421}

{\em 4.2.2) The Connecting Annuli and Selfconjugate Objects}\dfll\pageref{pg-422}

\subparagraph{4.3) Punctured Cobordisms and Glue-$\otimes$}\hfill\pageref{pg-43}

%% file: BL0.tex
\newpage

\setcounter{page}{1}
\pagenumbering{arabic}

\section*{0) Introduction and
Survey of Results}\label{pg-0}
\
\s

The discovery of new algebraic structures in the form of quantum groups and 
braided tensor categories (BTC's) has led to new insights in several areas of 
mathematical physics, and, in particular, has stimulated renewed interest in
low dimensional topology. This development was set off with the observation
that the braid group representations obtained from quantum groups or BTC's 
can be used to construct knot and link invariants. The general strategy
followed here 
is to associate to the ``regular'' singularities of a generic flat projection
of a link, namely crossings, maxima, and minima, the defining  braid- and
rigidity-morphisms of a BTC. To the higher order singularities of
projections of codimension one, (triple crossings, crossings at equal
heights, saddles, tangential strands, etc.)
which connected  all generic projections to each other, one associates the respective relations 
that occur in the axioms of a 
BTC, i.e., the Artin relations (or hexagonal- and pentagonal equations), the 
rigidity constraints, inverse braid isomorphisms, etc. 
\medskip

Soon thereafter it was realized that the same structures can be used to define invariants of closed, compact three manifolds. The works of Turaev-Viro and
of Reshetikhin-Turaev brought forward two types of invariants constructed from
tensor categories, which are related to each other. In this article we shall 
focus on the second type described in [RT], since it is closer related to link invariants. 
\medskip 

The starting point of the construction there are surgery presentations of a
manifold by a framed link $\li\hookrightarrow S^3\,$. In the description
of [Wc] the components of the link are the curves along which two handles  
$\hh^2$ are attached to a four ball. The resulting three fold is then given as
 the boundary  $M=\partial\bigl(D^4\cup\hh^2\cup\ldots\cup\hh^2\bigr)\,$.
From $\li$ the three fold  $M=M_{\li}\,$ can also be
described  directly by a surgery
along the components of $\li\,$, as in [Li].
\medskip

The manifold invariant of [RT] is then defined as a weighted sum over the
link invariants, described previously, evaluated on $\li\,$. The fact that
still needs to be checked is  that links presenting the same manifold
yield the same invariant. In particular we need to know how two links
$\li_1$ and $\li_2$ with  $M_{\li_1}=M_{\li_2}$ 
can be related.
\medskip

An answer to this question is given by Kirby's theorem, [Ki]. 
It states that $\li_2$ 
can be obtained from  $\li_1$ by a series of two types of moves in the
class of links in  
$S^3\,$. They are the two-handle slides ( ${\cal O}_2$-moves ) and the signature move (or ${\cal O}_1\,$-move ), which corresponds to the replacement 
$W\to W\#{\bf CP}^2\,$ for the bounding four fold. 
 This result is improved by the result in [FR], where it is shown 
that it suffices to consider only a special type of  ${\cal O}_2$-moves,
namely the $\kappa$-moves. 
\medskip

The main technical concern in both works is 
that the moves can be chosen such that one never leaves the class of 
presentations that use only two handles. For this purpose one also has 
to restrict the 
class of three folds on which we surger to connected and simply connected ones
(like $S^3\,$). If simply-connectedness is given up, Proposition
\ref{pp-red1} tells us that the Kirby formulation
(but not that of [FR] )  can 
be salvaged by including an additional move, which is introduced as  the $\eta$-move in
Section 2.2.2). Clearly, for disconnected manifolds not even the results of
[Wc] and [Li] hold, since surgeries along links do not change connectedness.  
\medskip

Consequently, the point of view chosen in the mentioned 
references involves a number
difficulties when dealing with the constructions of [RT]:
First the  algebraic verification of  the $\kappa$-move is somewhat involved,
since several pieces of the  $SL(2,{\bf Z})\,$ representation, that is associated  to any modular BTC, have to be constructed by hand.  
Moreover, the verification is not like in the case of links a more or less 
obvious  identification of elementary topological and algebraic constituents. 
Instead, it relies on the non trivial result in [FR], that presents all two handle slides as compositions of $\kappa$-moves.  
In the study of non semisimple invariants interpretations of two handle slides 
as the defining equation of an integral of a Hopf algebra have been found in [HKRL]. This, however, does not simplify the verification that the formulae 
given in [RT] are actually integrals. 
\medskip

In [RT] not only an invariant but more generally  the construction of a TQFT, i.e., fiber functors on cobordism categories, is proposed. 
A rigorous construction of presentations of cobordisms, which generalizes
the theory for closed manifolds, has so far been missing.  
For the special case of invertible cobordisms, which are given by elements
of the mapping class group, tangle presentations  have been constructed in
[MP]using  
explicit presentations in terms of generators and relations, as in [Wj].
 In order to generalize the surgery presentation as in [Wc], [Li], [Ki], and [FR], and for nice descriptions of the composition of cobordisms it is most
convenient to start surgery on unions of handlebodies. These  are neither 
connected nor simply connected and therefore do not fit in the 
conventional Kirby calculus.\footnote{After completion of the 
first version of this papaer the author received the preprint [Sw].
The approach presented there is similar to the description of [RT], where
the spines of embedded handlebodies are treated on the same footing as surgery
ribbons. Presentations are thus reduced to the ordinary Kirby calculus. Unlike 
the description we derive here, the embeddings of the boundaries are variable, which makes the formulation of composition rules (let alone presentations of cobordism categories)  more difficult.}
\medskip

In order to shed more light on both of the outlined problems, it is most 
instructive, if we refrain from insisting on 
presentations, where only two handles are attached, but also include 
one handles into the surgery description.
The manifolds we shall consider are thus of the form   $M\,=\,\partial\bigl(N^{(4)}\cup\hh^1\cup\ldots\cup\hh^1\cup\hh^2\cup\ldots\cup
\hh^2\,\bigr)\,$, and we make no connectedness assumptions on the boundary of 
$N^{(4)}$. The attaching data is now encoded into, what we call
 a ``bridged link'',
which is embedded in the  manifold $\li\hookrightarrow M_0=\partial N^{(4)}\,$,
to be surgered on. The additional data that enters a bridged link are 
pairs of small spheres in $M_0$, that indicate, where the one handles are to
 be attached. In analogy with presentations of groups this is like including 
more  generators into the description and thereby providing more freedom in 
choosing a practical set of relations.
\medskip

The one surgered manifold, 
$\breve M\,=\,\partial\bigl(N^{(4)}\cup\hh^1\cup\ldots\cup\hh^1\bigr)\,$, 
is obtained by gluing the spheres of each pair together along an
orientation reversing 
homeomorphism, after the interior three-balls have been removed.
 In particular, this type of surgery allows for a change of connectedness. 
To get from $\breve M\,$ to $M\,$ we attach the remaining two handles, i.e., 
we surger along an ordinary, framed link $\breve \li$ in $\breve M\,$. The preimage of the components of $\breve \li$ in $M_0\,$ are ribbons that may
end in one surgery sphere and emerge at the respective spot on the partner 
sphere. A typical example of a bridged link is given in figure~\ref{fig-typp}.
\medskip

Manipulations with similar types of attaching data have appeared in the papers
of, e.g., [FR]. A frequently used notation for one handle attachement due
to Kirby is that of a ``dotted circle''. However, the isoptopy classes of
such attaching data are generally smaller, because associated surgery
spheres cannot be moved independently. 
\medskip

With our definition, the analogue of the presentation 
theorems in [Wc] and [Li] is the following:
\begin{theorem}\label{ONE}
Suppose $M$ and $M_0$ are compact, oriented three folds, and $M$ is
connected. Also assume there is a homeomorphism of boundaries:
$$
\psi\,:\,\partial M_0 \isto \partial M\;.
$$
Then there is a bridged link $\li$ in the interior of $M_0$ such that $\psi$
can be extended to a homeomorphism
$$
\tilde \psi\,:\,\bigl(M_0)_{\li}\isto M\;.
$$
\end{theorem}
\medskip

The ``relations'' between presentations, which we propose to choose,  are
given by the following five moves. See also, Proposition \ref{pp-quint-mv}
in Section 2.3).

\begin{theorem}\label{TWO}
If for two bridged links $\li_1$ and $\li_2$ in the same compact manifold
$M_0\,$ we have that $\bigr(M_0\bigl)_{\li_1}\,=\,\bigr(M_0\bigl)_{\li_2}\,$,
then they are related by a sequence of the following moves:

\begin{enumerate}

\item \ub{Isotopies}: Regular isotopies of the bridged link $\li\hookrightarrow M_0\,$, where the ends of the ribbons stay attached to the spheres. 
\item \ub{Signature or ${\cal O}_1$ -move}:  (as in [Ki])

\item \ub{ Isolated Cancellation} : If a component of the bridged link consists 
of  a pair of spheres and a single ribbon which penetrates the spheres in
only
one pair of points, then this component can be discarded from the diagram.
(This corresponds to figure~\ref{fig-canc},  
where the outer ribbons, $a$,$b$, and $c$, are omitted.)

\item \ub{Modification (or Handle Trading)} : This corresponds to replacing
a one handle of the four fold by a two handle.
The operation on the bridged link is shown in
figure~\ref{fig-modi} and explained in Section 2.2.2.5).

\item \ub {One Slides and Isotopies over Components} : If a pair of surgery spheres (say $Lo$ and $Lo'$) lies in two different
  components of $M_0$, then we can push another surgery sphere ( $Hi$) through as shown in figure~\ref{fig-1slid}. An  isotopy of the link in a fixed $\breve M\,$ gives rise isotopies of ribbons through spheres as indicated in figure
\ref{fig-isot}. 

(In the description of planar ribbon diagrams those are broken into opposite braid group actions on the 
the ribbons attached to the two surgery spheres, and pushing through loops.)

\end{enumerate}

\end{theorem}

Since the original motivation of this work was to find cobordism
presentations, the language and organization used in the derivation of
these results sometimes suggests a
specialization to the case where $M_0\,$ is a union of handle bodies
$\H\,$. As no references to the special properties of this choice are
ever made, the proofs may be read literally in full generality.   
\medskip

All of the proposed moves in Theorem~\ref{TWO} are local. In particular, we have no two handle 
slides or $\eta$-moves 
in the list. Note also, that the fifth move can be omitted if $M_0\,$ is
connected. 
\medskip

 Diagrammatically,
the modification-move is similar to the $\kappa$-move, but its topological
derivation and its interpretation are far more elementary. 
Also, the algebraic computations involved
in the invariance proof of [RT] become much easier, see Section 4.2.
The proof of invariance we give here is thus closer in spirit to the 
proof of invariance for links as described in the beginning of the 
introduction. The relations between surgery and algebraic data are 
summarized in the table at the end of Section 4.2.1). In [KL] this
correpondence is sytematically put to use, in order to construct 
extended TQFT's and the elementary cobordisms that belong to the
algebraic generators are more explicitly identified.  
\medskip

Interestingly, we find - analogous to the correspondence in [HKRL] between 
2-handles and integrals of Hopf algebras - an interpretation of 1-handles
as cointegrals. Details and implications for non-semisimple TQFT's will
be discussed in separate papers.
\medskip

As a byproduct of purely topological considerations, we will show that the modularity condition of 
[Tu] on the abelian 
BTC we start with can be equivalently replaced by the weaker condition:
$$
1\in im(S) \qquad\qquad \bigl(\;{\rm or,\quad  equivalently},\; S^2(1)=1\;, etc.\bigr)\;,
$$
where the $S$-matrix is as usual defined by  the traces over the monodromies.
If we think of $S$ as a generalized Fourier transform, this condition can 
be seen as an analogue of the point separation condition  of the 
Stone-Weierstrass theorem. For the non-semisimple version of this
condition see again [KL].
\medskip

The second application of the bridged link calculus 
is the derivation of tangle presentations of 
connected cobordisms - or a central extension, $1\to\Omega_4\to
\wct\to\ct\to 1\,$, thereof described in Chapter 1. 
We introduce so called ``standard presentations''  on unions of standard
handle bodies with  a fixed one-handle structure of the bounding
four fold. The moves under these 
restrictions are derived in Chapter 3. An additional move that has to
be added to those for closed manifolds is  the so called the $\sigma$-move,
which also appeared in the
combinatorial description of the mapping class group in [MP].
In fact, our presentation will entail a 3+1 dimensional, Cerf-theoretical
 proof of the results in [MP],
which does not use the explicit presentations from [Wj]. The category
\mbox{${\cal T}\!\!g\,$} is described in the following theorem:
\medskip

\begin{theorem}\label{THREE}
The category $\,\wct\,$ is naturally isomorphic to the category \mbox{$\,\widetilde{{\cal
T}\!\!g\,}\,$} of ordered sets of grouped, admissible, ribbon tangles in ${\bf
R}\times [0,1]\,$, modulo relations.

\noindent
The composition rules in \mbox{$\widetilde{{\cal T}\!\!g\,}$} are as descibed in 
Sections 1.3) and 3.2.3). 

\noindent 
The five equivalence relations are
\begin{enumerate}
\item[] \quad Isotopies
\item[] \quad The $\tau$-move \quad\qquad (see figure \ref{fig-tau-move})
\item[] \quad The $\sigma$-move \quad\qquad (see figure \ref{fig-sigma})
\item[] \quad The $\kappa$-move \quad\qquad (as a special ${\cal O}_2\,$, see [FR])
\item[] \quad The $\,^0\bigcirc\mkern-15mu\bigcirc^0\,$-move \quad(a special $\eta$-move,
see Section 2.2.2.8))
\end{enumerate}
The functor $\wct\to\ct$ is presented by a corresponding quotient
\mbox{$\widetilde{{\cal T}\!\!g\,}\to {\cal T}\!\!g\,$} of tangle
categories. 
In \mbox{${\cal T}\!\!g\,$} we have in addition the ${\cal O}_1\,$-move, which
allows us to omit the $\,^0\bigcirc\mkern-15mu\bigcirc^0\,$-move.

\end{theorem}

In Section 4.1) tangle presentations are used to derive presentations of
general mapping tori.  
We also discuss a diagram that is of interest in connection with the 
composition rules of the tangle categories. This tangle element gives rise to
a canonical idempotent in a given quasitriangular Hopf algebra $H$, 
whose image on a general representation of $H$ is the maximal, self conjugate subrepresentation.
\medskip

The proofs of Theorem \ref{TWO} and related statements are based on the 
theory of stratifications of function spaces and their topology as
developed by Cerf in [Ce]. The relevant facts are reviewed in Section 2.1). 
Another indication that the inclusion of one handles yields a more coherent
picture is given by the observation that the space of codimension one Morse
functions on a four fold with singularities of index one and two is
connected, where as 
the corresponding space with only index two singularities is disconnected. 
(see Lemma~\ref{lm-12-cerf}). 
\medskip

In Theorem \ref{TWO} the listed moves have specific meanings as to
which part of the presentation they affect: The moves 1), 3) and 5)
correspond to the elementary deformations of functions on a fixed
bounding four fold. The ${\cal O}_1$-move is as usual the elementary move
that changes the cobordism class of the four fold in $\Omega_4$, i.e., the 
class in the central extension $\wct\,$, by connected summing with  a ${\bf
C}P^2\,$. 
Only the modification changes the four fold. The elementary operation it
involves is the attachement of a five dimensional 2-handle to
$W\times[0,1]\,$ along a curve ${\cal C}\subset W\times\{1\}\,$. If
$\cal C$ intersects the attaching data of a ( four dimensional ) 1-handle
$\hh^1\subset W$ in exactly one point, then the surgered manifold $(W)_{\cal
C}\,$ is given by trading $\hh^1\,$ for a two handle which is attached
along an annulus surrounding $\cal C\,$.

\medskip   

In Section 4.3) we also describe the presentations for categories
$\wct(N)\,$ of compact
surfaces with $N$ boundary components, which we work out completely
in [KL]. The associated
tangles contain $N$ additional strands.
A new non-trivial operation on the set $\{\wct(N)\}_N\,$  is the glue tensor product
$\otimes_{glue}\,$, where we do not just take the
disjoint unions of surfaces
but also sew them  along some of their boundary components. The main
difficulty we encounter here is that the glue tensor of two standard
surfaces is not canonically identified as a standard surfaces
anymore. 
We briefly discuss the role of double categories, introduced in [KL]
to describe the two gluing operations. 
\s

\paragraph{Acknowledgments:}\  

I thank Daniel Altschuler, Dror Bar-Natan,
Raoul Bott, Dale Husemoller, Kiyoshi
Igusa, Vaughan Jones, Louis Kauffman, David Kazh\-dan, Robion Kirby,
Volodimir Lyu\-ba\-shen\-ko, David Radford, Nicolai Re\-she\-tik\-hin,
Justin Roberts,
Stephen Sawin, Volker Schomerus, and Vladimir Turaev  for interest, comments, and discussion.  

I wish to thank Jim Harmon for helping with the conversion of eps-files.
Pictures were created with {\em FreeHand} and  {\em Xfig}. 

The author is partially supported by NSF grant DMS-9305715.

%% file: BL1.tex
\section*{1) Framed Cobordisms}\label{pg-1}


\subparagraph{1.1) Construction of Bounded  Cobordisms }\label{pg-11}
\
\s

For the surgery presentation of closed three manifolds it is essential to
find a bounding four manifold. The purpose of this section is find an analogous,
practical notion of a bounded  $2+1$ cobordism category.
We will define it as a special subcategory of $3+1$ dimensional
cobordisms. In order to organize the presentations according to their
signature we introduce $\wct$, which is
the subcategory of $3+1$ cobordisms modulo five cobordisms. 

The resutling category can equivalently be viewed as the one of 
surfaces and {\em two-framed} 2+1-dim cobordisms, by a result of 
Atiyah. 
\s

{\em 1.1.1) 2+1 Cobordisms}\label{pg-111}
\s

We shall use the following conventions to describe  a 2+1 dimensional
cobordism category. 
An  object is  given
by a sequence $\bar g\,=\,(g_1,\ldots,g_K)\,$ of non-negative integers.
(We admit $K=0$, i.e., $\bar g=\emptyset\,$). To any such sequence we
associate a two dimensional surface
$\Sigma(\bar g):=\amalg_j^K\Sigma_{g_j}\,$ with $K$ components. Here 
the $\Sigma_g$ are fixed, oriented coordinate surfaces, one for each genus. 
 
A morphism from $\bar g^-$ to $\bar g^+$ is an  (oriented) three
dimensional  cobordism. It consists of an oriented three manifold $M$ with
boundary
$\partial M=B_+\amalg B_-$ together with an orientation reversing 
homeomorphisms $\psi=\psi_+\amalg\psi_-\,$ with $\psi_{\pm}:
B_{\pm}\isto\pm\Sigma(\bar g^{\pm})\,$,
which  we shall call a {\em chart}.

The composition 
of (oriented) cobordisms is given by identifications along the
boundaries, using the charts of the start and end boundary component
respectively. We consider two cobordisms to be equivalent $(M,\psi) \sim (M',\psi')\,$,
if there is a homeomorphism  $\chi: M\isto M'\,$ 
such that $\psi'\circ \chi|_{\partial M}\,=\,\psi\,$.

The category also has an obvious tensor product, $\otimes\,$, which is given by the
disjoint union of surfaces and cobordisms. 
  
\begin{defi} We denote by  $Cob_3$ the category of 
 2+1 dimensional cobordisms. The objects are  sequences of non-negative
 integers $(\bar g)$ and the morphisms $[M,\psi]$ are equivalence classes of cobording
three manifolds. The composition structure is induced by the identification
along boundaries.
\end{defi}

\s


{\em 1.1.2) Bounding 3+1 Cobordisms}\label{pg-112}
\s

The classical link presentation for closed manifolds is 
given on the standard manifold $S^3\,$. Similarly, a presentation 
of manifolds with boundary $\Sigma=\partial M$ should be given on
a standard manifold with the same boundary $\Sigma\,$. A nearby choice are
unions of standard  
 handlebodies $H_g$ (one for every genus $g$). Since we also wish to
describe compositions of presentations of cobordisms, we shall also
consider the complementary handlebodies.
\s

For a precise definition  we fix for any non-negative integer $g$ 
an unknotted embedding of a  standard handlebody $H^+_g=H_g$ of genus $g$
into $S^3$ with a fixed orientation. We denote the opposite handlebody
$H^-_g:=S^3-H_g\,$ and make the identification $\Sigma_g:=\partial H_g\,$,
with induced orientation. Also, we write $\H(\bar g^+,l,\bar g^-)$ for the 
(ordered) union of the handlebodies and $l$ copies of $S^3\,$, giving rise
to a standard cobordism from $\bar g^-$ to $\bar g^+\,$. 
\s

Having this standard manifold for a given boundary we construct for every
2+1 cobordism $\bigl(M,\psi\bigr)$
a  closed, oriented three manifold $M^{\rm cl}\,$ by
\beq\label{clM}
 M^{\rm cl}\,=\,\Bigr(\H(\bar g^+,l,\bar g^-)\Bigl) {\coprod} _{\sim}
\Bigr( \{\Sigma(\bar g^+)\amalg-\Sigma(\bar g^-)\}\times
[0,1]\Bigl)\coprod_{\psi} M\;.
\eeq   

Here $\sim$ stands for the standard identification of the lower boundary of
the cartesian product with the boundaries of the handlebodies.  
By a classical theorem of [Ro] we know that there always exists a compact
four fold $W\,$ such that
\beq\label{bd-pres}
M^{\rm cl}\,\cong\,\partial W\;.
\eeq
\s

This allows us to consider a subcategory of 3+1 dimensional cobordisms,
which admits a full functor onto $\ct\,$. Its objects are of the form
$\Sigma\times [0,1]\,$ for a closed, compact, and oriented two fold $\Sigma$
and can thus be identified with the ones from $\ct\,$. Without referring to
the bounded three fold as in (\ref{clM}) we define the morphisms as four
folds $W$ with special functions at the boundary. Specifically, we
assume the existence of a height function $\h\,$ and a chart $\psi$;

$\hphantom{which satisf}\h \,: \, U(\partial W)\to [0,1]\,$ in a vicinity of the boundary

$$\psi\,:\,\Bigl( \bigl\{\Sigma(\bar g^+)\amalg-\Sigma(\bar g^-)\bigr\}\times
[0,1]\Bigr){\coprod} _{\sim}\Bigl(\H(\bar g^+,l,\,\bar
g^-)\Bigl)\hookrightarrow \partial W\,,$$

which satisfy the conditions:

\begin{enumerate}
\item $\h$ is smooth and has no  critical points .
\item  $\h^{-1}(0)=\,\psi\Bigl(\H(\bar g^+,l,\bar g^-)\Bigl)$
\item $\h^{-1}(1)\subset \, \partial W$
\item $im(\psi)\,=\,\overline {\h^{-1}\bigl([0,1[\bigr)}\cap \partial W$
and $\h\circ \psi $ is the projection onto $[0,1]$ if restricted 
to $\Bigr( \Sigma(\bar g^+)\amalg-\Sigma(\bar g^-)\Bigr)\times
[0,1]$ 

\end{enumerate}

Also, we shall always assume that the bounding four cobordism
has components in one to one correspondence with the  
components of the three cobordism, i.e., we require
 
\beq\label{concon}
\pi_0\bigl(\h^{-1}(1)\bigr)\,\longrightarrow
\,\pi_0(W)
\eeq
to be an isomorphisms.

\s

\begin{figure}[ht]
\begin{center}\ 
\epsfbox{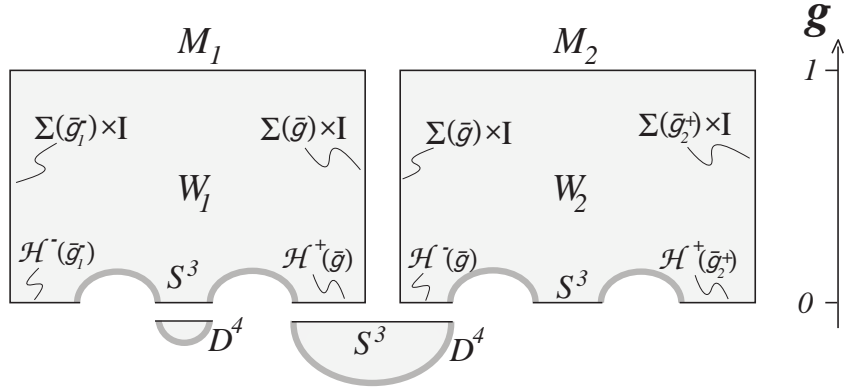}
\end{center}
\caption{3+1 Cobordism}\label{fig-cob-cat}
\end{figure}

The composition of two such manifolds $W_i\,,\,i=1,2\,,$ is given by
identifications along  the
common boundary pieces $\Sigma(\bar g)\times[0,1]\,$ using the
charts $\psi\,$. The new functions, $\h$ and $\psi$, are the restriction of the old
ones. If $l_i$ is the number of $S^3$ components in the boundaries of
$W_i\cap \h^{-1}(0)$ the respective number for the comopsition is
$$
l=l_1+l_2+\sum_{j=1}^K g_j\,.
$$
Since we wish to give presentations in only one $S^3\,$, we often redefine the
composition by gluing in four balls, $D^4$,  along the excessive $S^3$'s, at the
expense of introducing additional index 0  singularities for the composite
of the Morse
functions on the $W_i\,$'s. Typically, if $M$ is connected,  these will be
cancelled with other index 1 
singularities. A schematic picture of the composition is given in
Figure~\ref{fig-cob-cat}.

\s


{\em 1.1.3) The Category $\wct$ }\label{pg-113}
\s

Suppose a four fold $W_1$ bounding $M$ can be obtained from another four fold
$W_2$ by surgeries in the interior. If we had in addition a 
Morse function $f$ on $W_1$ we can easily find ``modified'' Morse functions
on the surgered manifolds, which coincide with $f$ outside the surgered
pieces, and whose singularity structure differs from that of $f$ in a specific way.
As Morse functions on $W$ define presentations of a three fold
$\subset\partial W\,$, surgeries on the 
four folds  yield 
moves between presentations of the same $M$ as boundary of either $W_1$ or $W_2$.
The surgery, and thus  the moves, only exist if  $W_1$ and $W_2$ are
cobordant with proper identifications of the boundaries.
This leads us to consider a wider notion of equivalence of the special 3+1
cobordisms than just the homeomorphy type:
\s

In order for two connected cobordisms $(W_1,\psi_1)$ and $(W_2,\psi_2)$ to be considered
equivalent
there shall be an orientation preserving  homeomorphism 
$\hat \chi :\partial W_1\isto \partial W_2\,$ which is compatible with
 the charts $\psi_1$ and $\psi_2$ . As in three dimensions (see
(\ref{clM})), we construct a closed four
dimensional manifold $\hat W$ from the pieces $W_1\,$,$W_2\,$ and
$[0,1]\times\partial W_i\,$ using the identification $\hat \chi\,$. 
This manifold shall be the boundary of a five dimensional manifold $Q\,$.
Since $\hat \chi$ was assumed to be orientation preserving an orientation 
on $\hat W$ is opposite to the orientation of exactly
one constituent. Thus, by e.g. [J], the signature  $\sigma(\hat W)$ is the
difference of signatures of $W_1$ and $W_2$ and it follows from [Wa2]
that $Q\,$ exists if and only if $\sigma(W_1)\,=\,\sigma(W_2)\,$.
If $(W_1,\psi_1)$ and $(W_2,\psi_2)$ are disconnected than we say that they are
equivalent iff all of their connected components are equivalent. 

\begin{defi} We denote $\wct$ the categories of 2+1 cobordism
bounding special 3+1 cobordisms, modulo 4+1 cobordisms. 
The objects are sequences of non-negative integers and
the morphisms are the equivalence classes of fure manifolds $[W,\psi ]\,$, for which smooth 
boundary functions $\h$ exist. The composition is given by identification 
of respective boundary pieces and shrinking of excessive $S^3\,$'s.
\end{defi}           


{\em 1.1.4) Moves between Four Manifolds}\label{pg-114}
\s

In the presentation of a three dimensional manifold we can change either
the bounding four fold $W$ or the Morse function on it that describes a
handle decomposition of $W$: The move that result from changing the Morse
function and an elementary modification of $W$ are described in the next
chapter. In this subsection we wish to give the possible changes in the
handle decomposition of $W\,$, if we change representatives of
$[W,\psi]\,$. Other than in [Ki] we need to account for the fact that
$\pi_1(W)$ may be non trivial and surgeries on $W$ may thus be non local.
\footnote{I am indebted to Justin Roberts for remarking the lack of treatment 
of   this  point in an earlier version of this paper.})
Here we shall use also techniques and definitions that will be explained 
later in this paper. 

\medskip

If we give a surgery description of a {\em connected} cobordism $M$
 starting from $M_0$  (e.g. $\,=\,\H(\bar g^+,l,\bar g^-)\,$)
the four fold $W$ is given by $M_0\times I$ with four dimensional
$k$-handles attached to the upper boundary. It is clear, using cancellations
and connectivity arguments, that we may omit all $0$- and $4$-handles.
In fact we are interested in four folds, which admit a handle
decomposition of the form 
\beq\label{eq-PPP}
W\,\cong\,M_0\times I\, \cup\hh^1\cup\ldots\cup\hh^1\cup\hh^2\cup\ldots\cup
\hh^2\;,
\eeq
where the handles are attached to the boundary piece
$int(M_0)\times\{1\}\,$. 

\medskip

Let us also introduce $\breve W = M_0\times I\,
\cup\hh^1\cup\ldots\cup\hh^1\,\subset \,W$ and the three fold
$\breve M$ which is the upper boundary component such that 
$\partial\breve W\,=\,M_0\amalg_{\partial M_0}\breve M\,$.
One useful feature of $\breve M$ is that it contains all
homotopy, more precisely:

\blm\label{lm-homot}

For $j=0,1$ the following are an isomorphism and an epimorphism:
$$
\pi_j(\breve M)\;\isto\;\pi_j(\breve W)\,\hbox to 25pt{\thrafill}\,\pi_j(W)
$$

\elm 

{\em Proof :} In fact for $\pi_0$ both maps are isomorphisms. To show that the inclusion
of $\breve M$ into $\breve W$ is an isomorphism it suffices to check that
for any four fold $W\,$, $\,\pi_j\bigl(\partial(W\cup\hh^1)\bigr)\,\isto\,
\pi_j\bigl(W\cup\hh^1\bigr)\,$ is an isomorphism whenever
$\,\pi_j\bigl(\partial W\bigr)\,\isto\,
\pi_j\bigl(W\bigr)\,$ is one (for every component). Since $(\partial W)\cup \hh^1\,$ is obtained
from $\partial\bigl(W\cup\hh^1\bigr)\,$ by gluing in  a $D^3\times I\,$
along the 0-,1-connected piece $S^2\times I\,$, it follows immediately
that $\pi_j\bigl((\partial W)\cup
\hh^1\bigr)\,\cong\,\pi_j\bigl(\partial(W\cup\hh^1)\bigr)\;$ for
$\,j=0,1\,$.  In the case where $\hh^1$ is attached to two different connected components
$\partial W^{\alpha},\,\partial W^{\beta}\,$ of $\partial W$, with
corresponding connected components $W^{\alpha},\,W^{\beta}\,$ of $W\,$,
 the spaces $W\cup \hh^1\,$ and  $\partial(W)\cup \hh^1\,$ are given 
- up to homotopy type - by replacing the respective components by $\partial
W^{\alpha}\vee\partial W^{\beta}\,$ and $W^{\alpha}\vee W^{\beta}\,$. If
$\hh^1$ is attached  to the same component of $\partial W$ we end up with 
$\partial W\vee S^1$ and $W \vee S^1\,$. In both cases the assertion follows
easily, since $\pi_1$ is freely generated by known parts. A simple Seifert van
Kampen argument shows that $\pi(W\cup\hh^2)\cong \,\pi_1(W)/[C]\,$, where
$C$ is the attaching curve of the two handle, completing the
 verification  of Lemma~\ref{lm-homot}.\hfill$\Box$

\medskip

For a four fold $W$ with a handle decomposition as in~(\ref{eq-PPP}) we may
attach a pair of two handles without changing the boundary. More precisely,
assume we have attached an $\hh^2$ to $\breve W$ along any curve $C\subset \breve M$.
Suppose $D$ is a small disk which intersects the attaching data of all
other two handles exactly once in $C$. Then $\partial\bigl(\breve W\cup
\hh^2\cup\hh^2\bigr)\,=\partial \breve W\,$, where the second two handle is
attached along $\partial D$ with framing induced by $D$. In the language of
bridged links this corresponds to the ``$\eta$-move'', and will be discussed in
more detail in Section 2.2.2.8). We may now describe the precise relation
between four folds giving rise to the same classes in $\wct$.

\bpp\label{prop-conn}
        
\begin{enumerate}
\item For any connected cobordism class $[W,\psi]$ there exists a
representative $(W_0,\psi)\,$, such that $W_0$ is of the
form~(\ref{eq-PPP}).
\item Suppose two connected four cobordisms $(W_j,\psi_j\,)$ with handle
decompositions as in~(\ref{eq-PPP}) for $j=1,2$
give rise to the same class in $\wct\,$. Then there exists a representative
$(W_3,\psi_3)$ which is of the form (\ref{eq-PPP}) and can be obtained from
either $W_1$ or $W_2$ by a sequence of $\eta$-moves.
\end{enumerate}
\epp

{\em Proof :} As explained in the proof of Lemma~\ref{lem-one} any $W$ with
connected upper boundary $M$ is representable without 0- or 4-handles. To
show the first assertion it thus remains to replace the 3-handles. Suppose
the upper parts $M$ and $M'$ of the boundaries  $\partial W$ and  $\partial
\bigl(W\cup \hh^3\bigr)$ are connected. Reading the cobordism piece from
$M$ to $M'$ in different directions we have $W':=\bigl(M\times I\bigr)\cup
\hh^3\,\cong \,\hh^1\cup\bigl(M'\times I\bigr)\,$. Using a modification as 
described in Section 2.2.2) we have another cobordism between the same
manifolds $W''=\hh^2\cup\bigl(M'\times I\bigr)\,=\,\bigl(M\times
I\bigr)\cup \hh^2\,$. We may connect $W''$ to a bunch of ${\bf C}P^2$'s
or $\overline{{\bf C}P^2}$'s (which is a special type of two-handle
attachement) such that $\sigma(W''')=\sigma(W')\,$. The closed manifold 
$W'''\amalg_{\partial W'}W'\,$ thus bounds some
$Q^{(5)}\,$, so that we have $W\cup\hh^3\sim W\amalg_M
W'''\,=W\cup\hh^2\cup ...\cup\hh^2\,$. 
 
\medskip

By definition there is a connected five fold $Q$ cobording $W_1\,$ to 
$W_2\,$. As in [Ki] we may forget about 0- and 5-handles and make
modifications  on $Q$ that replace all 1- and 4- handles by 3- and 2-
handles. We may push the  remaining 2- handles close towards the boundary
$W_1\,$ and the remaining 3-handles towards $W_2\,$. Hence $Q$ is given by 
the composite $Q=Q_1\amalg_{W_3}(-Q_2)\,$ of cobordisms, build up from
handles of the same type.
For both pieces we have that $Q_j$ cobords $W_j\,$ to $W_3$, and that 
it is given by attaching five dimensional two handles
$\bigl(W_j\times I\bigr)\cup\hh^2\cup...\cup\hh^2\,$. Now, a corresponding two
surgery on a four fold $W$ is determined by an attaching curve
$S^1\hookrightarrow W\,$. A neighborhood $D^3\times S^1\,$ is removed 
and an $S^2\times D^2\,$  is glued in along the new boundary component
$S^2\times S^1\,$. It remains to be checked that this is the same as
attaching the pair of (four dimensional) two handles described in the
$\eta$-move above.

\medskip

For a four fold $W$ of the form (\ref{eq-PPP}) we may apply
Lemma~\ref{lm-homot}  to find a (PL-) homotopy of the attaching curve to a
path $\cal C$ in $\breve M\,$. By transversality we may choose
${\cal C}\subset\breve M\,$ to be without selfintersection and the homotopy
in $W$ to be an ambient isotopy of curves. Thus we may assume that
$S^2\times D^2\,$ is glued in along a curve in $\breve M\,$. Since 
all framings of the curve in four dimensions are homeomorphic, we may 
choose it such that the upper and lower hemisphere of the fiber
$D^3\,=\,D^3_+\amalg_{D^2_{eq}}D^3_-\,$ lie in an upper and a lower collar
of $\breve M\,$. I.e., we have $S^1\times D^3_{\pm}\hookrightarrow \breve
M\times [0,\pm 1]\subset \breve M\times [-2,2]\,\subset W\,$, and
$S^1\times D^3\,\cap \breve M\times\{0\}\,=\,S^1\times D^2_{eq}\,$.
Correspondingly, we may decompose the sphere of the newly attached piece
into two hemispheres $S^2=S^2_+\amalg_{S^1_{eq}}S^2_-\,$. As a result we
obtain the decomposition
\bea
\Bigl(\breve M\times[-2,2]-(S^1\times D^3)\Bigr)\amalg_{S^1\times
S^2}D^2\times S^2\,\quad=\hphantom{xxxxxxxxxxxxxxxxxxxxxxxxxxxx}\nonumber\\
\Bigl(\bigl(\breve M\times[-2,0]\bigr)-\bigl(S^1\times
D^3_-\bigr)\Bigr)\amalg D^2\times S^2_-\amalg D^2\times S^2_+\amalg
\Bigl(\bigl(\breve M\times[0,2]\bigr)-\bigl(S^1\times D^3_+\bigr)\Bigr)\,\nonumber  
\eea
The first and fourth piece are clearly homeomorphic to $\breve M\times I\,$
with  $R_{\pm}=S^1\times S^2_{\pm}\,\subset \breve M\,$  lying now in the upper
(lower) 
boundary part  and we may omit one of them (e.g. the last one) since the attachement is
just thickening a boundary piece. The gluing of the second piece is nothing
but a two handle attachement of $\hh^2=D^2\times S^2_-\,$ to $\breve M\times
I\,$ along the  framed knot $R_-\,$. The same is true for the gluing of the
third piece. Only now the attaching data $D^2\times S^1_{eq}\hookrightarrow
D^2\times \partial(S^2_-)\,\subset\partial(\breve M\times I\cup\hh^2\,)$ is
inside of the surgered region. We may push a  circle $\{p\}\times S^1\subset D^2\times
\partial S^2_-\,$ outside of this region by moving $p\in D^2$ to a point
$p'\in S^1\,$, so that it is a meridian $\{p'\}\times \partial
D^2_{eq}\subset\partial R_-\,$ of the attaching knot of the first handle.
Similarly, it follows that for a small interval $J\subset\,D^2$ the
attaching ribbon $J\times\partial D^2_-\,$ can be pushed outside of the
surgered region to an annulus as described in the $\eta\,$ move.  
 
\hfill$\Box$
\s


\subparagraph{1.2) A Central Extension by $\Omega_4$ }\label{pg-12}
\
\s

By construction a  morphism in $\wct$ bounds a morphism in the original 
category $\ct$. Explicitly, it is given by the restrictions
\beq\label{cob-restr}
M=\h^{-1}(1)\qquad{\rm  and}\qquad \psi|_{\{1\}\times \Sigma}\;,
\eeq
which are clearly compatible with the compositions.
In summary, we have a canonical functor,
$$
\wct \stackrel {\cal D}{\lra} \ct \;,
$$
which is full by Rohlins result stated in (\ref{bd-pres}). In this section
we investigate the structure of $\cal D$ in more detail. In particular, we
interpret $\wct$ as a central extension of $\ct$ by the cobordism group $\Omega_4\,$,
which is related to different central extensions by [At], [Ko], and [RM]. 
\s


{\em 1.2.1) The $\Omega_4$-Action and Anomalous TQFT's}\label{pg-121}
\s

To  a
connected four fold $Y$ with boundary function $\h\,$  we can canonically assign a 
natural transformation $\xi_Y\,$ of the identity functor on
the subcategory of connected four cobordisms.
We first apply the construction (\ref{clM})
to the identity $id=\Sigma(\bar g)\times [0,1]\,$ of $\ct$, so that
$(id)^{cl}\,\cong\,S^3\,$. From this we obtain $\xi_Y(\bar g)$ 
by identification of $(id)^{cl}$ with the boundary of $D^4\# Y$, 
i.e., a four ball with $Y$ connected to it. Clearly, the result of
composing a connected morphism $[W,\psi]$
in $\wct$ with $\xi_Y$ on either side is $[W\# Y,\psi]\,$. In particular, 
we have $\xi_Y\circ\xi_Z\,=\,\xi_{Y\#Z}\,$. 

The images of these transformations in $\wct$ only depend on the class of
$Y$ in the cobordism group $\Omega_4\cong{\bf Z}$, and 
thus form a free abelian group generated by $y\,=\,\bigl[\xi_{{\bf
CP}^2}\bigr]\,$. Therefore, we can view the  functor $\cal D$ 
as a quotient map to the orbits of the morphisms under the free action of
$\Omega_4\,\subset Nat_{\wct}(id)\,$. In this picture we can define a  {\em
projective} (or {\em anomalous}) TQFT  as a fiber functor on $\wct$.   
(Examples are the TQFT's constructed in [RT] and [Tu] using quantum groups,
see also Section 5.)

Assuming that $\bar g=(0)$ is assigned to a one dimensional vectorspace,
and the functor on $D^4\in End(0)$ is non-zero, 
this yields a number $\theta\neq 0\,$  for $y\,$, which determines all of
$\xi$.
Clearly, the  TQFT-functor  factors through $\cal D$ into a fiber functor on $\ct$ (an
ordinary or  ``anomalie free''  TQFT) if and only if $\theta=1\,$. 
\s
 

{\em 1.2.2) Another Two-Cocycle}\label{pg-122}
\s

The extension of $\ct$ is in fact non-split. The
precise behavior can be given by  a ``two-cocycle'' $\mu$ of $\ct$,  which
measures the  
non-additivity of the signature for the  composite of two four manifolds
$W_1$ and $W_2$.  
It turns out that $\mu$ only depends on the charts of  $M_i=\h_i^{-1}(1)$. 
We write 

$$
\begin{array}{rl}
\psi_1\,:\,\Sigma(\bar g )&\to M_1\amalg_{\sim} \H^-(\bar g_1^-)\qquad
\psi_1\,:\,\Sigma(\bar g )\to M_2\amalg_{\sim} \H^+(\bar g_2^+)\\
&\quad {\rm where\; }\quad \,\bar
g\,=\,\bar g^+_1\,=\,\bar g^-_2\,
\end{array}
$$
 
with identifications $\sim$ along the common surfaces. We denoted by
$\,\H^{\pm}(\bar g)$ the respective unions 
of handlebodies. In the following discussion we shall disregard the extra
$S^3\,$, which have no effect on the signature.
To give  explicit expressions we introduce the Lagrangian subspaces (in
rational homology)  
$$
\Lambda_i=\,ker\bigl(H_1(\psi_i)\bigr)\,\subset\,H_1\bigl(\Sigma(\bar g)\bigr)
$$
and
$$
V^{\pm}\,=\,ker\bigl(H_1(i_{\mp})\bigr)\;,
$$
where  $i_{\pm}$ is the inclusion of $\Sigma(\bar g)$ into
$\H^{\pm}\bigl(\bar g\bigr)\,$ so that $H_1\bigl(\Sigma(\bar
g)\bigr)\,=\,V^+\oplus V^-\,$.

If $\omega$ is the skew form on $H_1\bigl(\Sigma(\bar
g)\bigr)\,$, let us define a bilinear $\phi'$ form on 
$U'=\Lambda_1\,+\,\Lambda_2\,$
as follows. For $\lambda\in U'\,$ denote $\lambda^{\pm}$ the components in
the subspaces $V^{\pm}\,$. Also, we choose a $\lambda_i\in\Lambda_i$ with
$\lambda=\lambda_1+\lambda_2$. Then
\beq\label{defphi}
\phi(\lambda,\eta)= \omega (\lambda_2,\eta_1)\,-\,\omega(\lambda^+,\eta^-)
\eeq
It is easily checked that $\phi$ is symmetric and does not depend on the
choices of $\lambda_i$ and $\eta_i\,$. Also, it factors into a bilinear
form
$\phi\,:\,U\times U \to {\bf R}\,$ on the quotient space  
$$
U\,=\,\frac {\Lambda_1\,+\,\Lambda_2}{\Lambda_1\cap
V^-\,+\,\Lambda_2\cap V^+}\;.
$$
We denote by 
$
\mu(\psi_1,\psi_2)$ the signature of $\phi\,$. 

\s
Applying a result from [Wa] we find the  anomalie of the signatures:

\bpp
 We have 
\beq\label{addsig}
\sigma(W_1\circ W_2)\,=\,\sigma(W_1)+\sigma(W_2)+\mu(\psi_1,\psi_2)
\eeq
\epp 

{\em Proof :} As in [Wa] we denote by
$Y_+=W_2\,$, $Y_-=W_1\,$, $X_+=M_2\amalg \{\Sigma(\bar
g^+_2)\times[0,1]\}\amalg \H(\bar g^+_2)\amalg \H^-(\bar g)\,$ with
orientation induced by $Y_+\,$, $X_-=M_1\amalg \{\Sigma(\bar
g^-_1)\times[0,1]\}\amalg \H^-(\bar g^-_1)\amalg \H(\bar g)\,$ with
orientation opposite to $Y_-\,$, $X_0=[0,1]\times\Sigma(\bar g)\,$ with
orientation from $Y_-$, and $Z=\Sigma(\bar g)\times {0} \amalg
\Sigma(\bar g)\times {1}\,$ where the orientation on the ${1}$ component is
standard and on the ${0}$ component opposite.(Consult again
Figure~\ref{fig-cob-cat}.)  

With this orientation the bilinear form  on the two component space
$H_1(Z)=V_0\oplus V_1$ is given by $\tilde {\omega}=-\omega\oplus
\omega\,$. From the inclusions of $X_+,\,X_-\,$, and $X_0$ we have
Lagrangian subspaces 
$$
A\,=\,V^-_0\oplus \Lambda_1\qquad\quad C\,=\,V^+_0\oplus \Lambda_2\;
$$
and
$$
B\,=\,\Bigr\{(-x,x)\,:\,x\in H_1(\Sigma(\bar g))\,\Bigl\}\;,
$$
where we consider $\Lambda_i\subset V_1\,$. It is now straight forward to
show that the form in [Wa] reduces to (\ref{defphi}). Its signature is
identified with the signature anomalie of the four folds. Also, the gluing of the
four ball in the composition does not change the signature.
\ep 

\s

Let us record the explicit anomalies in two special cases of interest.

\blm
\ben
\item The anomalie vanishes if  one of the $\Lambda_i$ coincides with a 
standard subspace $V^{\pm}\,$.
\item For the torus $\Sigma((1))$ assume that  $e^{\pm}\in V^{\pm}$ 
$\alpha_i\in\Lambda_i$ are non zero vectors in the one dimensional
Lagrangian subspaces. Then the anomalie $\mu(\psi_1,\psi_2)$ is the
signature of the binary form $\phi$, where
$$
\phi_{ij}=\omega(\alpha_i,e^-)\omega(e^-,e^+)\omega(e^+,\alpha_j) \quad
{\rm with } \quad i\leq j
$$
\een
\elm

{\em Remark:} Our definition of a bounded cobordisms is different from
the description in [Ko] or [Wk] and yields the actual signature of the
linking diagram
on a standard manifold. The cocycle gives also rise to a non-split
central extension of mapping class groups, if we restrict to invertible
morphisms of $\wct\,$. It is related to well known central extensions, 
which are described in terms of canonical framings of $TM\oplus TM\,$ in [At].

\s

\subparagraph{1.3) Ordering and Connectivity}\label{pg-13}
\
\s

We chose the objects of $\wct$ and $\ct$ to be ordered sequences of
surfaces. It is clear that reordering does not change the isomorphism
class. In fact, we have canonical isomorphisms 
\beq\label{perm}
[\pi]\,:\,(g_1,\dots,g_K)\,\isto\,(g_{\pi(1)},\ldots,g_{\pi(K)})
\eeq
for every $\pi\in S_K$, which are homeomorphic to unions of
$\Sigma_{g_j}\times[0,1]\,$. In $\wct$ we bound the $[\pi]^{\rm
cl}$ 's by $D^4$'s.  
\s

The fact that  a cobordism $M$ in $\ct$ is has connected components $M_j$ is  expressed by
the formula
\beq\label{fff}
M\,=\,[\pi]\bigl(M_1\otimes\ldots\otimes M_r\bigr)\;.
\eeq

From this it follows easily that every cobordism $M:\bar g\to \bar g'\,$ in
$\ct$ can be written as a product 
of cobordisms  of the form
\beq\label{elem}
[\pi]\bigl( N\otimes i\!d_{g''_j}\otimes \ldots \otimes i\!d_{g''_{K'}}\bigr)\;,
\eeq
where $N$ is a connected cobordism and $[\pi]$ a permutation as in~\ref{perm}. 

\s

By the connectedness assumption in (\ref{concon})
 a presentation will thus be a collection of presentations of the 
connected components. It is not efficient to try to evaluate the 
presentation of the composition of two cobordisms of
the form~(\ref{fff}) directly. We prefer to expand the product first by
using the most elementary composition rules. These are to write
$M_1\otimes\ldots\otimes M_r$ as the product of the commuting cobordisms 
$id\otimes \ldots id\otimes M_j\otimes\ldots\otimes id\,$, the rule 
$[\alpha]\cdot[\pi]=[\alpha\circ\pi]\,$ for the permutations, and 
$[\pi]M_1\otimes\ldots\otimes M_r\,=\,M_{\pi (1)}\otimes\ldots\otimes
M_{\pi (r)}[\pi]'\,$, if $\pi$ is the  permutation of ordered subsets.  It is not
hard to see, that for the correct expansion and correct bracketing the only
other compositions that will occur are of the form
\beq\label{perm-comp}
[\alpha]\,\cdot\,N \qquad \quad N\,\cdot\,[\alpha]
\eeq
and 
\beq\label{elem-comp}
L\,:=\,N\otimes id \,\cdot\,id\otimes M\;,
\eeq
where the cobordisms $N\,$, $M\,$, and $L\,$ are {\em connected}
and $\alpha$ is a permutation.

%% file: BL2.tex
\section*{2) Enhanced Surgery Presentations}\label{pg-2}
\

Assume we are given a compact  $n$-manifold $X$ and a Morse Function
$f:X\to [0,1]\,$. It is an elementary fact (see [Mi])  that if $f$ has a 
singularity of index $j$ with non-degenerate critical value $a$, then the 
sublevel manifold $X^{\leq a+\epsilon}=f^{-1}([0,a+\epsilon])\,$ is
homeomorphic to $X^{\leq a-\epsilon}$ with a $j$- handle $\hh^j$ attached 
at its boundary. This means, we have an embedding $S^{j-1}\times
D^{n-j}\hookrightarrow f^{-1}(a-\epsilon)\subset \partial X^{\leq
a-\epsilon}\,$, and  we attach $\hh^j=D^j\times D^{n-j}\,$, along the 
respective piece in the boundary $\partial \hh_j\,=\,\partial D^j\times
D^{n-j}\,\cup\,D^j\times \partial D^{n-j}\,$.   

The classical link presentations of [Wc] and [Li] of three folds result
from attaching two handles to a four ball such that
$M=\partial\bigl(D^4\cup\hh^2\cup\ldots\cup\hh^2\bigr)\,$. The theorem of
Kirby  (see [Ki])
states that two such  presentations (with the same signature)
 can be related by sliding the handles across each other.
Both results rely on the fact that the manifold $S^3=\partial D^4$ we
surger on is connected and simply connected. 
Neither is  true for the manifolds $\H\bigl(\bar g^+,l,\bar g^-\bigr)\,$.

For these reasons we shall admit presentations that involve  also surgeries
with one handles. Most of this section is devoted to describe in detail the 
calculus of ``bridged links'', which is the corresponding generalization of the
link calculus for closed manifolds. Some elements of this description were also
used in [FR]. We shall also discuss the calculus of ``bridged ribbon graphs''
analogous to [RT0] or [RT] for generic, planar projections of the link
diagrams. 
\s

 
\subparagraph{2.1) Some Cerf Theory}\label{pg-21}
\
\s

A presentation for a manifold using a Morse function $f$ as above is only
defined if the critical points and critical values of $f$ are non-degenerate.
Since the  space of these functions is disconnected it is not always
possible to deform the presentation of one function into the presentation
given by another function. Nevertheless, if we admit also functions of
codimension one we obtain again a connected space. 
\s

The content of Cerf's theory [Ce] is to determine the connectivity of function
spaces depending on their codimensions,  and the behavior of the singularities  
once a path of function passes through a lower stratum. In this section we
shall review the elements that are relevant for our purposes. We start by
introducing the  function spaces, that will describe the bridged link
presentations and its moves.
\s

{ \em 2.1.1) Excellent and Codimension One Functions}\label{pg-211}
\s

For given $W$ and boundary function $\h$ introduce
 a Riemannian structure such that the gradient flow of $\h$ is parallel to
the boundary piece $\psi\Bigl([0,1]\times \Sigma\Bigr)\,$. As usual we
introduce the space ${\cal F}$ of smooth functions on $W$ which coincide
with  $\h$ near $\partial W\,$. Also, we denote by 
${\cal F}^0$ the ``excellent functions'', which are Morse functions with
distinct critical values, and by ${\cal F}^1\,$ the codimension one functions.
The latter set is the union  of  the set ${\cal F}^1_{\alpha}$ of functions
with distinct 
critical values and only one degenerate critical point, which is a birth
(or death) point, and the set ${\cal F}^1_{\beta}\,$  of Morse functions 
for which  exactly two critical values coincide.

We denote the
stratification of codimension one:

$
\hphantom{xxxxxxxxxxxxxxxxxxxxx}\f_T=\fo\cup\fla\cup\flb
$

\noindent
We also denote by 

$
\hphantom{xxxxxxxxxxxxxxxxxxxxx}p\f\subset\f\qquad\qquad p\fo\subset\fo\qquad {\rm etc.}
$ 

the subspaces of functions for which all  singularities have 
index one or two, all index two singularities have higher values than
the index one singularities, and all singularities are in general position.  

Correspondingly, the codimension one strata in there are denoted

$
\hphantom{xxxxxxxxxxxxxx}\;p\fla,\quad\;p\flb\,,\quad\;p\fo_i\qquad {i=1,2}\qquad\;p\fo_{12}\;.
$

They are the functions with one 1-2 birth point, functions with exactly two critical
values coinciding, excellent functions for which the ascending and
descending manifolds of two index $i$ singularities
intersect in a single trajectory, and functions where the descending
manifold of a 2-singularity and the ascending manifold of a 1-singularity
intersect with an intermediate level three fold in a curve and a sphere that
are tangential to each other in exactly one point ( and transversal in all
other). 

Their union with $p\fo$ is denoted
$p\f_T\,$. 

It is occasionally useful to label the singularities of index
one (three) as to  which components of the surgered manifold are connected 
by the respective one handles.

\begin{defi}
\ben
\item 
The {\em  label}  of a critical point $c$ of index one (three) is the pair
$(\lambda,\mu)$ of connected 
components of $\h^{-1}(0)\,$ which are intersected by the descending
(ascending) manifold of $c$. 
\item For a function $f\in p\fo\,$ we denote by $\breve M$ the level
manifold $f^{-1}(y)\,$, where $\,y\,$ lies above all critical values of  index
one and below all critical values of index two.

\een  
\end{defi}

It is clear that $\breve M$ can be obtained from $\H$ by either connecting
different component to each other or to copies of $S^1\times S^2\,$.
In fact the connection prescription may be entirely determined from the set
of labels of the index one singularities. 

\s

{\em 2.1.2) Paths and Deformations of Paths of Functions}\label{pg-212}
\s

A basic result in [Ce] asserts that the space $\f_T$ is 0-connected and any
path in $\f_T$ can be deformed to a path which is transversal to the lower
dimensional strata $\fl\,$.  By the general  theory of stable manifolds, a
 presentation changes along a path in $\fo$
only by an isotopy.
 A complete
set of  ``moves'' for presentations on the same four fold are  thus given
by the (transversal) passages  through $\fl\,$.

 A generic path of functions in $\f_T$ is conveniently  illustrated by a graphic.
This is a collection of paths $\bigr(f_t(c^i_t),t\bigl)\,$ where
$c^i_t$ is a continuous family of critical points of $f_t$. 
The point where $f_t$ passes through
$\fla$ is given by a beak joining an index $i$ and an index $i+1$ critical
point. Going through $\flb$ corresponds to crossing of two components of the
graphic. Also, we draw a dashed line at points of $t$ where the ascending
and the descending manifolds of two critical points of the same index
intersect each other in a trajectory.

\begin{figure}[ht]
\begin{center}\ 
\epsfbox{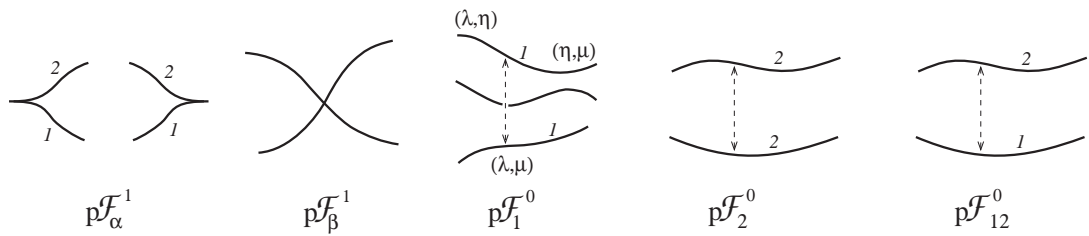}
\end{center}
\caption{Codimesion-One Strata in Graphics}
\label{fig-codone}

\end{figure}

In particular, a path through the codimension one strata of $p\f_T$ is
presented by the pieces of a graphic given in~\ref{fig-codone}.

In the graphic for the intersection of the manifolds of index one critical
points we also indicated the  possible labels.
\medskip

The next more general result in [Ce] concerns  deformations of generic
paths. They may always be chosen such  that the paths pass transversally
through singularities of 
codimensions two. For more general  
results on functions with framings and a useful and humorous summary of
singularities of higher codimensions see [Ig]. The list of elementary
deformations can be  illustrated by graphics as in [Ce], 
[Ki], and [HW, pg26]:

\ben
\item {\em Independent Trajectories}: If the descending manifolds of a  trajectory
$c_t$ of singularities do not intersect the ascending manifolds of another
trajectory $c'_t$ for any $t\,$ then the two trajectories can be moved
independently from each other. (e.g., [HW,pg.65]). For the relevant 
examples see [Ki].
\medskip
  
\item {\em Triangle Lemma} : For $i_1+i_3\leq 3\,$, $inf(i_1,i_3)\leq i_2-1\,$ or 
$i_1=i_2=i_3\leq 2\,$ we have the move:

\begin{center}\ 
\epsfbox{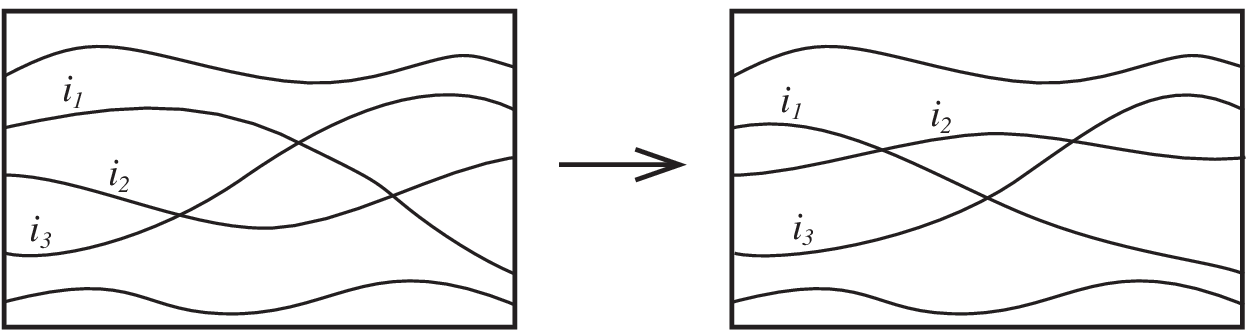}
\end{center}

\bigskip

\bigskip\ 

\item {\em  Beak Lemma}: If $i>0$ we have the move
\begin{center}\ 
\epsfbox{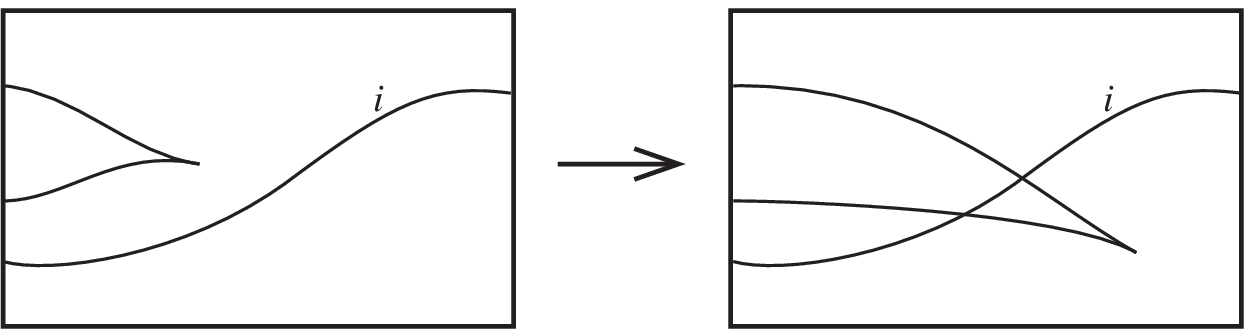}
\end{center}

\item {\em  Dovetail Lemma}:

\begin{center}\ 
\epsfbox{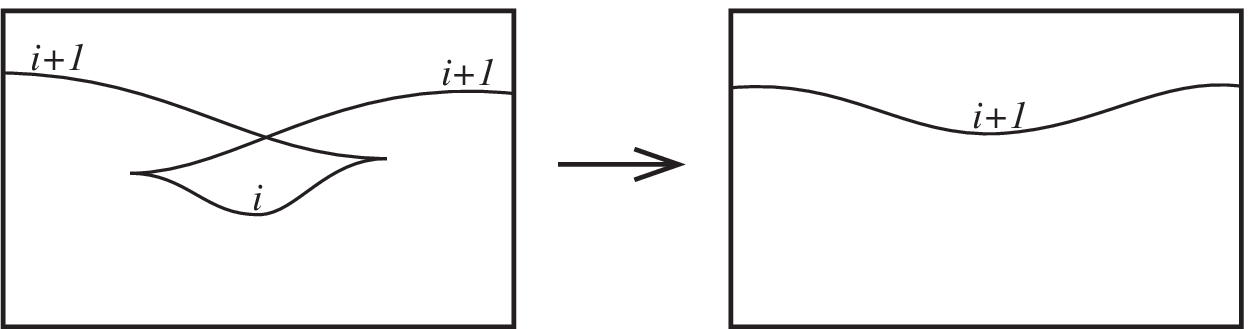}
\end{center}

\een

The remaining moves are obtained from the above by reflections $t\to 1-t$
and $f\to 1- f\,$.

\subparagraph{2.2) Bridged Link Presentations}\label{pg-22}
\
\s

We introduce an enhanced presentation of a three fold by adding both
one and two handles to another compact, oriented three fold, which does 
not have to be connected or simply connected. 
\s 

{\em 2.2.1) Surgery Presentation from Bridged Links}\label{pg-221}\nopagebreak
\s

From the structure of the singularities of a  function $f\in \fo$ we obtain
 a
presentation of the cobordism $\h^{-1}(1)$ by considering the intersections 
($\cong S^0,{\rm or }\,\cong S^1$) of the descending (or unstable) manifolds of
the critical points with 
$\h^{-1}(0)\,$. These yield together with framings of their normal bundles
a unique surgery prescription. To be more precise, let us introduce the
following standard handles:
\s  

{\it 1)A two handle $\hh^2\,=\,D^2_a\times D^2_b$ and two intervals $I,J\cong [-1,1]$
such that $D^2_a=J\times I\,$.
\s

2)A 2-sphere $S^2\subset {\bf R}^3$ with meridian $S^1=\{(x,y,0)\in S^2\}$
and involution $\rho:S^2\to S^2\,:\,(x,y,z)\to(x,y,-z)\,$.
Also denote by $B^3$ the oriented ball bounded by $S^2$ and by $D^2_{\pm}$
the upper and the lower hemisphere.} 
 \s

 On $\H=\H(\bar m^+,l,\bar m^-)$ for arbirary $l$
and a function $f\in p\fo $  a {\em bridged link presentation}
consists of the following data:

\s

{\it 1) An orientation preserving embedding $\phi=\bigl(\phi_1,\phi_1',\ldots,\phi_m'\bigr)$ 
of $2m$ copies of $B^3$.}

\s

The pairs of balls are neighborhoods of the points where the descending
trajectories of the index one singularities intersect $\H$. 
The three fold $\breve M$ just above
all index one singularities is then (for a suitable given Riemann structure)
 naturally
homeomorphic to the manifold obtained by removing the interiors of the
images of the three balls and identifying the pairs of boundaries using the
isomorphisms
$$
\bar\phi_j :=\phi_j'\circ \rho \circ \phi_j^{-1}\;.
$$

Occasionally, we shall prefer the picture of an actual handle attachement,
where we glue in an additional $S^2\times [0,1]$ in between the boundary components. 
\s

The surgery is described further by considering the intersections of the
descending manifolds of the index two singularities with $\breve M\,$. Close
to the singularity, the normal (Morse) form of $f$ determines - up to
homotopy -
a trivialization of the normal bundle. This, in turn, is determined by a
non vanishing section of the normal bundle in any level manifold, i.e., an
embedding of the ribbon: 
$
R:=J\times \partial D^2_b \,\hookrightarrow \breve M\;, 
$
or more precisely of a link of ribbons. Undoing the index one surgery we
obtain the next ingredient of the surgery data:
\s
     
{\it 2) A map $\tau=(\tau_1,\ldots,\tau_s)\,:\,R\amalg\ldots\amalg R\,
\hookrightarrow \,\H$ of $s$ copies of the ribbon $R$ which is an embedding
into the complement of the balls $im(\phi)$ except for a finite number of
intervals 
$J\times \{p\}\subset R$. The right and
left sided limits $\tau_j^{\pm}\,:\,J\times \{p\}\to \H$ are embeddings
into  spheres $\phi_{k(p)}(S^2)$ and $\phi_{k(p)}'(S^2)$ such that}

$$
\bar \phi_k\tau_j^{\pm}=\tau^{\mp}_j\;.
$$
\s

 Since we will be dealing with non-simply connected manifolds (with
no canonical framings of tangent bundles over the one-skeletons) we 
preferred here the language of ribbons over that of framing numbers. 

\begin{figure}[ht]
\begin{center}\ 
\epsfxsize=5.5in
\epsfbox{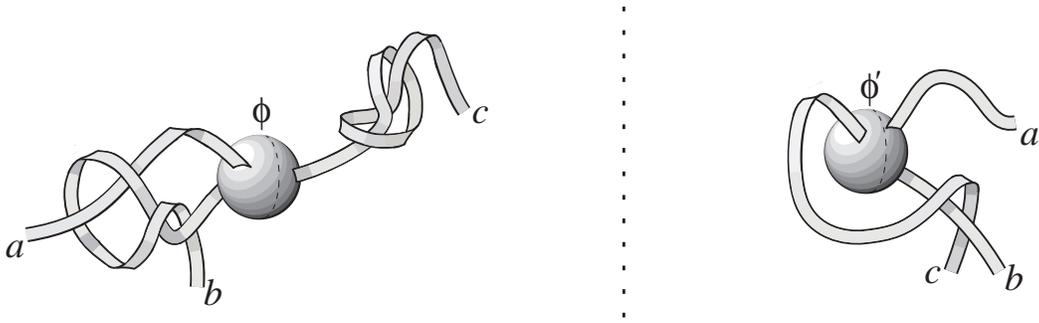}
\end{center}

\caption{Bridged Link}\label{fig-typp}
\end{figure}

A typical picture of the ``bridged link'' is given in Figure~\ref{fig-typp}. The
action of the flip $\rho$ is indicated by the invariant
meridians, drawn here as dashed lines.

It should be kept in mind  that the two spheres may lie in different
components of $\H$ depending on the label of the critical point, and 
are thus not connected by any cancelling ribbon (as in move 4) of the next
section).

The surgery on a component of the resulting ribbon link in $\breve M$ is given as usual by
extending the embedding to the tubes $D^2_a\times \partial D^2_b\supset
J\times 0\times \partial D^2_b\,$. Then we remove the images and glue in 
$\partial D^2_a\times D^2_b$ with reversed orientation. 
\s


{\em 2.2.2) Moves for Bridged Links}\label{pg-222}
\
\s

In this section we compile a list of ``moves'' of bridged link
presentations that do not change the homeomorphy class of the 
presented three fold. Some of them come from deforming the Morse
function through a codimension one stratum, some are obtained by
modification of the bounding four fold. 
\s

\begin{figure}[ht]
\begin{center}\ 
\epsfxsize=5.5in
\epsfbox{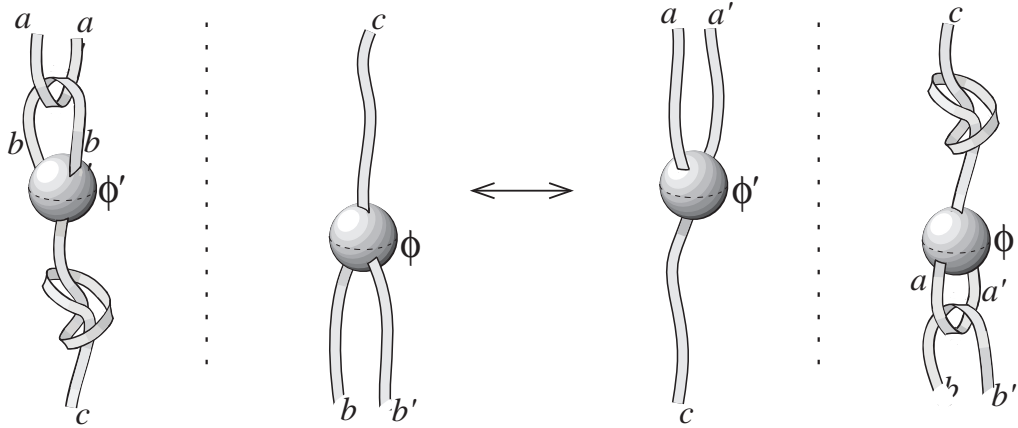}
\end{center}

\caption{Isotopy over Surgery Sphere}\label{fig-isot}
\end{figure}

1)\ub{\em Isotopies} An isotopy is a superposition of isotopies
of the three balls in $\H$ and of the ribbon link in $\breve M$. 
More precisely, after identifying a pair of embedded
spheres and their vicinity via the maps $\phi$ and $\phi'$ 
with the respective standard regions in $\bf R^3\,$, we may consider a pair of regions
$F\subset S^2$ and $F'=\rho(F)\subset {S^2}'\,$. A collar
$(U,F)\cong(F\times I,F\times 0)$ (opposite to the bounded balls) of one region,
can be moved
continuously to the other ending with the isomorphism:
$$
F\times I\,\stackrel{\rho\times (1-t)}{\hbox to 1.5cm{\rightarrowfill}}
\,{F'\times I}\;.
$$
Extending this isotopy (or its inverse) to the identity outside of $F$ and $F'\,$, we
can move pieces of the link diagram through a
gluing sphere as indicated in Figure~\ref{fig-isot}.

We confined generic presentations to have  to 1- and 2-singularities
in general positions; yet,  an isotopy over a surgery sphere also involves situations
where the surgery spheres are not always transversal to the ribbons.
At these points we can deform our path such that we pass only through the
stratum $p\fo_{12}\,$. This means only one ribbon is tangential to the
sphere at a time, and this configuration results from pushing a small 
ribbon loop through the sphere.

\s

\begin{figure}[ht]
\begin{center}\ 
\epsfxsize=5.3in
\epsfbox{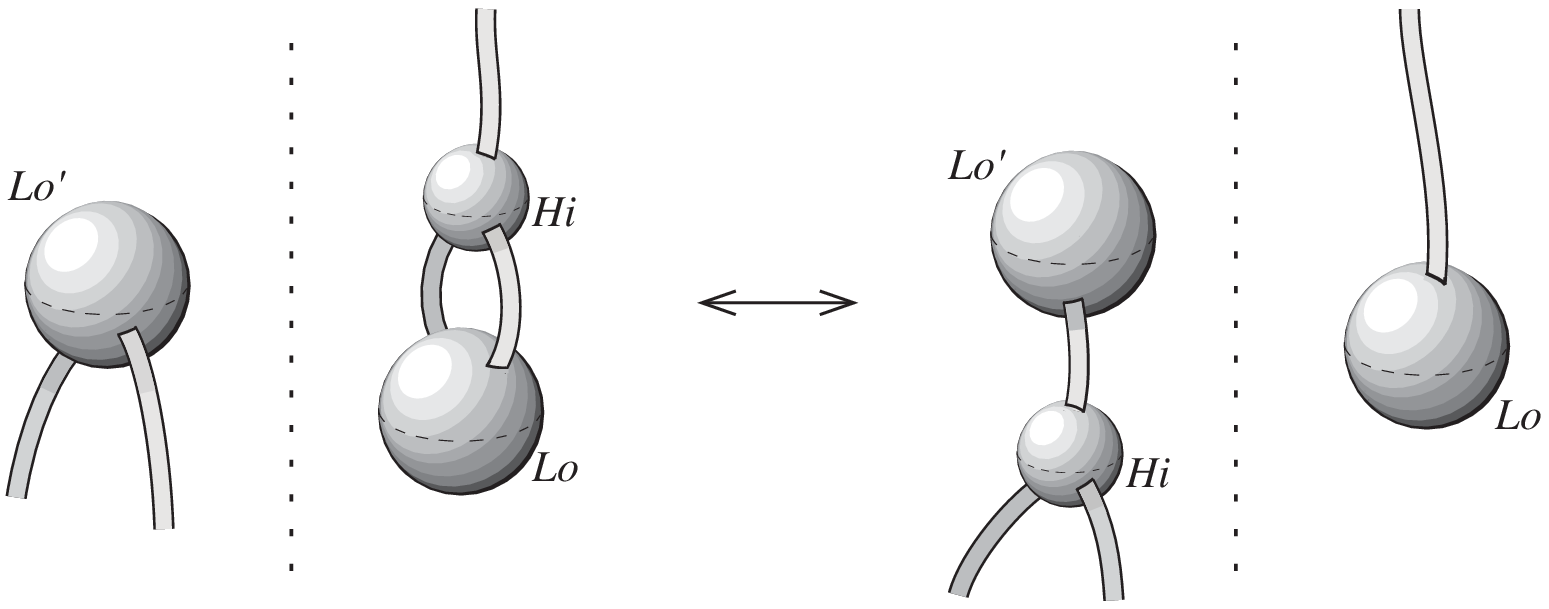}
\end{center}

\caption{One Slide}\label{fig-1slid}
\end{figure}

2) \ub{\em One- (Handle) Slides} They occur when  the Morse function passes
through $\fo_1$ and the  ascending manifold of an
index one singularity $Lo$ is moved transversally through one of the 
descending trajectories
of another index  one singularity, $Hi\,$, with higher value. Considering a
level manifold in between the two singularities it is clear that this 
corresponds
to moving a ball with the attached ribbons through a gluing sphere, as 
was described above for links.   
A  typical situation is suggested in Figure~\ref{fig-1slid}.

\s

\begin{figure}[ht]
\begin{center}\ 
\epsfxsize=5.3in
\epsfbox{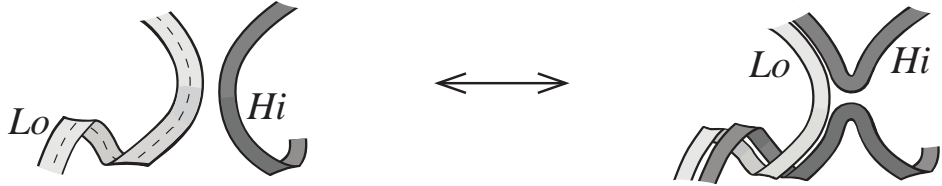}
\end{center}

\caption{Two Slide}
\label{fig-2slid}
\end{figure}

3) \ub{\em Two- (Handle) Slides} This is the Kirby move, denoted ${\cal O}_2$ in
[Ki]. It occurs when the path of  
Morse functions passes through $\fo_2\,$, and the 
intermediate manifolds of
the two index two singularities are moved through each other transversally.
At the point where the descending two fold of the singularity $Hi\,$ 
with higher
value goes through the other critical point $Lo$ we first have the 
corresponding
ribbons tangential along small pieces of their  boundaries. To describe the result
of the slide, the ribbon of $Lo$
is cut along its middle and the closer part is connected to the ribbon 
$Hi$ as
depicted  in Figure~\ref{fig-2slid}. On an intermediate level manifold this
is just an isotopy of the ribbon $Hi$ over the surgered piece along $Lo\,$.

\s

4) \ub{\em 1-2- (Handle) Cancellation} 
 The Cancellation Lemma of Smale (see, e.g. [HW] pg. 172) states that 
when for consecutive  singularities $c_1$ and $c_2$ 
of index $i$ and $i+1$ the intersection
of the descending  manifold of $c_2$ intersects the ascending manifold of
$c_1$ in exactly one trajectory, then $c_1$ and $c_2$ can be canceled 
against each other.
 This is achieved  along a path of codimension one functions passing
through $\fla$ ,i.e., an $i\,-\,i+1\,-$ birth- or death-point.
\s

\begin{figure}[ht]

\begin{center}\ 
\epsfbox{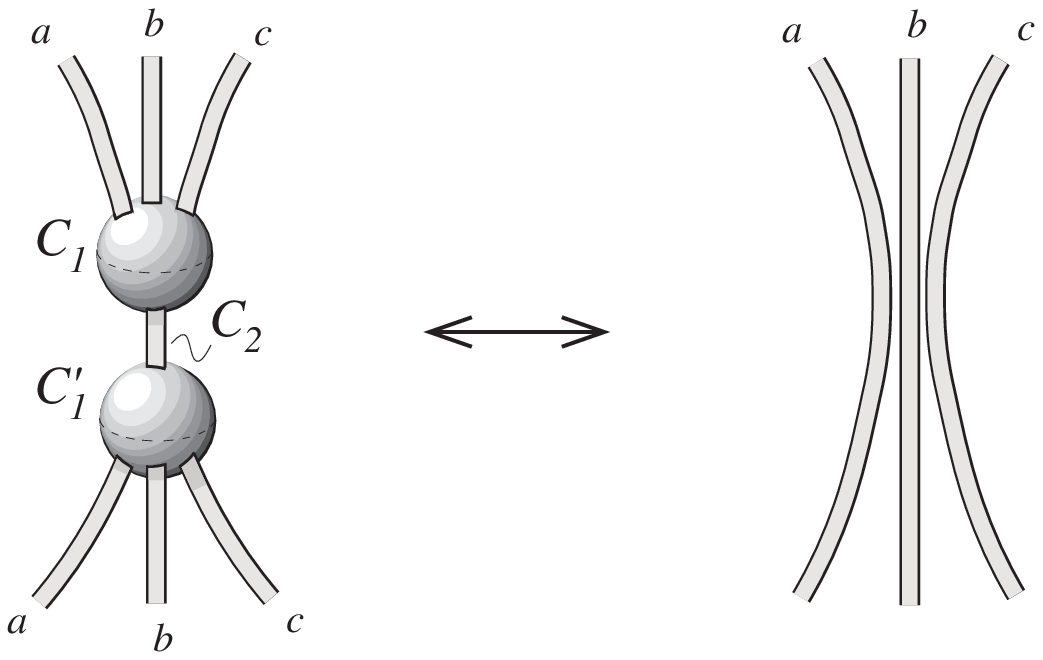}
\end{center}

\caption{1-2-Cancellation}
\label{fig-canc}
\end{figure}

Here we are only interested in the situation where $f$ passes through 
$p\fla\,$, and where the singularities are presented by a pair of balls $c_1$
 and a ribbon $c_2$ intersecting the
gluing sphere of $c_1$ in exactly one interval but no other spheres. In
this
situation we remove in addition to the balls a solid tube $D^2\times I$
around $c_2$ and
collapse the opposite hemispheres of the total removed region using $\rho$
and the framing of $c_2$, see Figure~\ref{fig-canc}.

\begin{figure}[ht]

\begin{center}\ 
\epsfbox{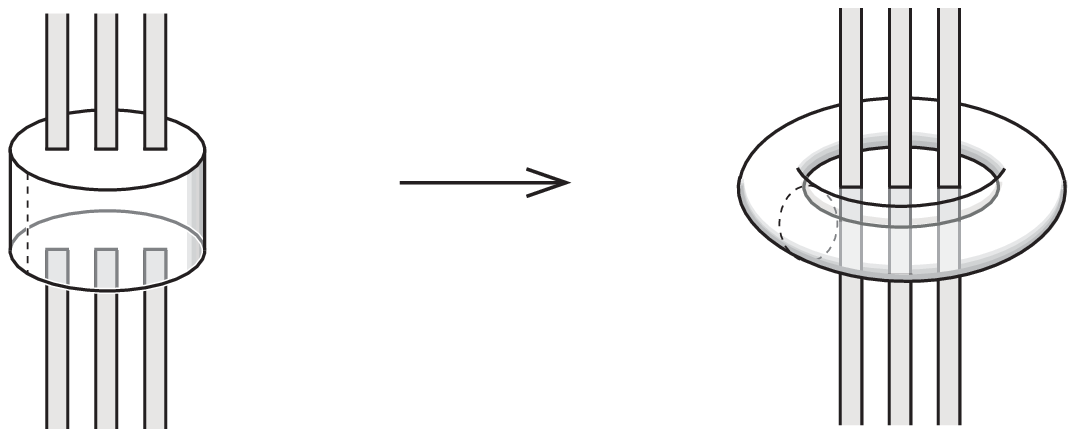}
\end{center}

\caption{Cancellation Surgery}
\label{fig-canc-pr}

\end{figure}

In order to explicitly show that this move preserves  the surgered manifold
we perform the
one surgery,
where we already removed the solid torus of the two surgery. This
corresponds to identifying the remaining hemispheres, which leaves us with 
Figure~\ref{fig-canc-pr}, where the torus is empty. Completing the two
surgery we obtain the right hand side of 
Figure~\ref{fig-canc}.

\s
 
\begin{figure}[ht]

\begin{center}\ 
\epsfbox{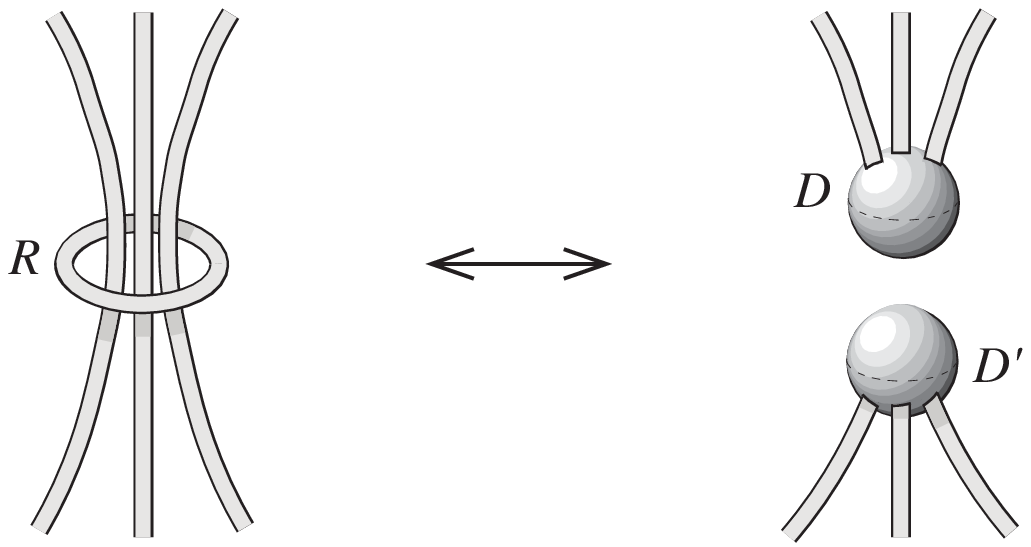}
\end{center}

\caption{Modification}
\label{fig-modi}
\end{figure}

5) \ub{\em  Modification (Handle Trading)} This  move  also changes the
bounding four manifold, exchanging an
index one singularity (or handle) for an index two singularity. Suppose we 
have a ribbon-component $R=S^1\times [1,2]$ which can be extended to an
embedding of a disc $D=S^1\times [2,0]/S^1\times 0$, (which may intersect
other ribbons-components nicely). Then we duplicate $D$ and push the two 
copies away from each other. We complete each of the discs to a sphere, so
that possible ribbons are attached to the outside, and define the maps
$\phi$ and $\phi'$ such that $\bar\phi$ is the correct identification of
the hemispheres $D$, see Figure~\ref{fig-modi}.

\begin{figure}[ht]

\begin{center}\ 
\epsfxsize=5.3in
\epsfbox{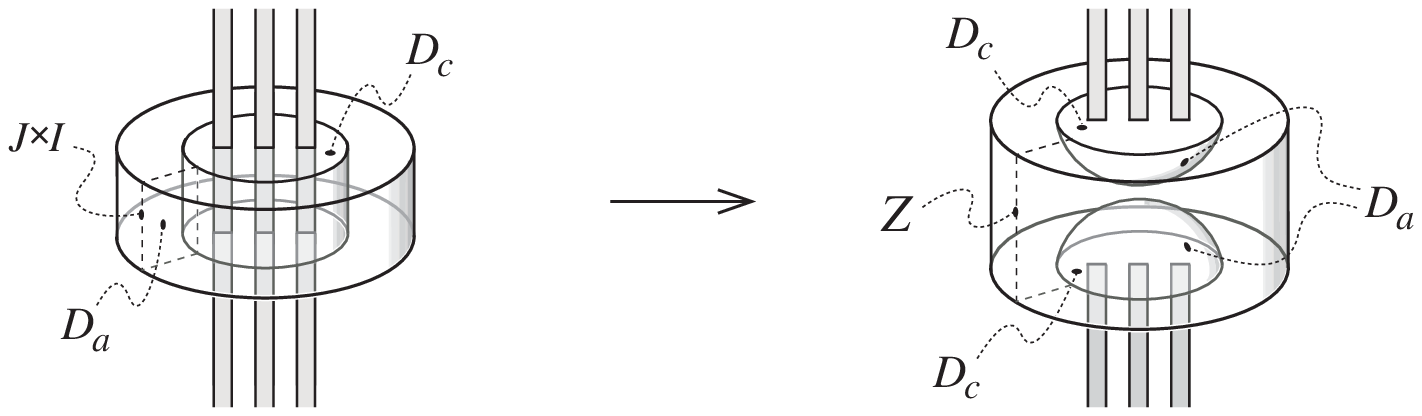}
\end{center}

\caption{Modification Surgery}
\label{fig-modi-pr}
\end{figure}

The procedure of modification is described in [Wa]. To see  that
the identification of the hemispheres is the correct one we shall make the
argument explicit: We parametrize a tubular neighborhood of $R$ by
$J\times I \times \partial D_b^2$ and the thickened, bounded disc by
$J\times D_c^2$, with identification $\partial D^2_b\,=\,\partial D^2_c$.
Instead of removing only the solid torus we also take out the thickened disc and
glue it back later. The opposite torus $\partial D_a\times D_b\,$ consists of
two parts $J\times 0\times D_b^2$ and $Z\times D_b^2$ with
$Z=\bigl(\partial D_a^2-J\times 0\bigr)$, which are identified 
along two copies of $D_b\,$. If we glue  the second part into the empty
region along $Z\times \partial D_b^2\,$, we obtain th region in
Figure~\ref{fig-modi-pr} where the insides of the two spheres $D^2_b\cup
D^2_c$ are missing.

The other part, $J\times 0\times D_b^2$  is first glued to the extra piece  $J\times
D^2_c$ along the cylinder $J\times \partial D^2$ giving $J\times S^2$ and
the included ribbons are $J$ times an interval on $S^2$. 
Shrinking $J$ to a point, it is clear from the picture 
that we obtain exactly the desired one surgery.

\medskip

To describe equivalences of cobordisms we need to include two more moves,
which change the signature of $W$ and the 
number of components  of $\H$:

\medskip

6)\ub{\em Signature} Add an unknotted, ribbon separated from the rest by a
sphere  with framing $\pm 1$. This is the Kirby move ${\cal O}_1$ and
results from 
connecting a projective space ${\bf CP}^2$ ($\overline{{\bf CP}^2}\,$) to $W\,$. We sometimes denote
such a link component by $\FLO\,$ ($\frL\,$).

\medskip

7)\ub{\em 0-1-Cancellation} If for a bridged link presentation an $S^3$
component of $\H$ 
contains nothing but a single surgery sphere we may omit the $S^3$-component,
and the respective pair of surgery spheres. (The content of the presentation in
this  situation is to connect the $S^3$ component to another component which is
of course trivial.)  

\medskip

Finally, let us introduce a move which is a combination of the previous
ones:
\medskip

8)\ub{\em $\eta$-move, $^0\bigcirc\mkern-13mu\bigcirc^0$-move} 
  An $\eta$-move is the addition  or subtraction of  a pair
of ribbons $\,^0\bigcirc\mkern-15mu\subset^{\cdots}_{\cdots}\,$ where one
ribbon extends to a disc which is intersected exactly 
once by the second ribbon $\cal C$ but no other ribbon of the presentation.
As described in the proof of Proposition~\ref{prop-conn} this corresponds
to the two-surgery $\,W\,\to\,(W-S^1\times D^3)\cup D^2\times S^2\,$.

When the ribbon $\cal C$ can also be contracted to an annulus, and the pair
of ribbons can be isolated from the rest,  we call the associated $\eta$-move an
$^0\bigcirc\mkern-13mu\bigcirc^0$-move. It corresponds to the surgery $W\to
W\# S^2\times S^2\,$. Using two slides shows that
$^0\bigcirc\mkern-13mu\bigcirc^0= ^0\bigcirc\mkern-13mu\bigcirc^{2n}\,$.
Thus, we shall tacitly assume that the
$^0\bigcirc\mkern-13mu\bigcirc^0\,$-move also implies addition and removal
of an isolated  $^0\bigcirc\mkern-13mu\bigcirc^1\,=\,\FLO\cup\frL\,$-piece. This corresponds
to connected summing with the non-orientable $S^2$-bundle $S^2\tilde\times S^2\,$.


\subparagraph{2.3) Equivalences of Bridged Links}\label{pg-23}
\
\s


{\em 2.3.1) Completeness of Moves for Bridged Links}\label{pg-231}
\s

A surgery presentation of a three manifold ( e.g., a morphism of $\wct\,$),
may be changed by either changing the bounding four fold $W$ within its
 4+1-cobordism class or by changing the handle decomposition of $W\,$.
 As indicated in the beginning of this chapter it is equivalent to speak
of excellent functions on $W$ instead of handles so that we may use the
results of Cerf's Theory from Section 2.1) and the correspondence between
codimension one strata of $p\f_T\,$ and the handle slides cancellations
given in Section 2.2.2). We know that any two excellent functions in
$p\f_T\,$ can be connected by a path in $\f_T\,$. Yet, this space is too large as it
contains singularities of all indices so that we have to deal with all
types of handles. For moves between bridged links we are only interested in
decompositions of the form (\ref{eq-PPP}), which is equivalent to
considering only functions in $p\f_0\subset p\f_T\,$. We need the following
refinement of the connectivity result from [Ce].
\medskip

\blm\label{lm-12-cerf}
The space $p\f_T\,$ is 0-connected.
\elm

\begin{figure}[ht]

\begin{center}\ 
\epsfbox{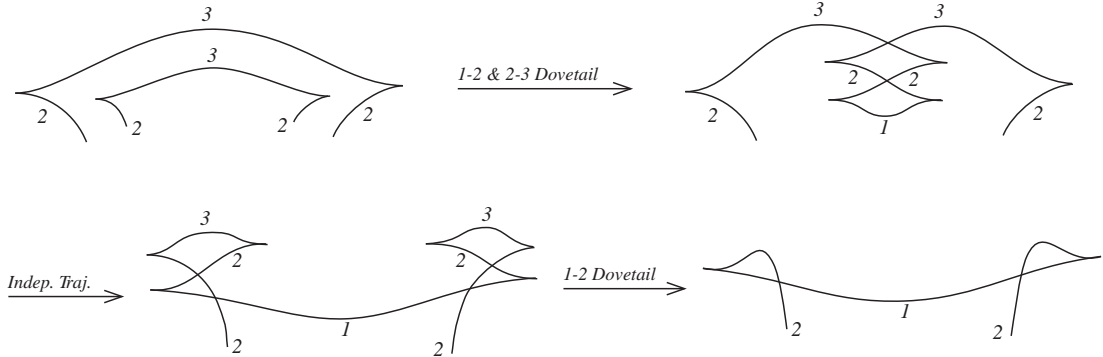}
\end{center}

\caption{Elimination of Three Handles}
\label{fig-cerf-pf}
\end{figure}

{\em Proof :} The elimination of components of a graphic with index zero and
index four  can be literarily taken from [Ki]. The trajectories of index three
singularities start and end in birth - and death - points, and, similarly as for
index one trajectories in [Ki], can be rearranged as disjoint trajectories
as indicated in first
graphic of Figure~\ref{fig-cerf-pf}, using the Beak Lemma. 
Introducing a double Dovetail we find the second graphic, and another
application of the Beak Lemma and the Principle of Independent Trajectories
yields the third picture. Removing the remaining Dovetails we find a graphic 
which contains only index one and two singularities. By the Principle of
Independent Trajectories we can move all index one trajectories below the
index two trajectories. We thus end up with a path in $\f_T\,$.
\hfill$\Box\,$.
\s
 
Put in a different language, Lemma~\ref{lm-12-cerf} states that all
handle decomposition of a four fold of the form~(\ref{eq-PPP}) can be
related to each other by isotopies, 1- or 2-slides, and 1-2-cancellations.
If we combine this result with Proposition~\ref{prop-conn} we arrive at the
following completeness of moves:

\begin{cor}\label{co-compl} If two bridged links on $\H$ present the same manifold with
boundary, then they 
are related by a sequence of the moves {\em 1)- 8)} from Section 2.2.2).
\end{cor}
\medskip

Note, that the modification has been added to the list, even though it
was not used in the preceding arguments. 

In fact it can be expressed as 
an $\eta$-move, where the long ribbon $\cal C$ is chosen as
a connecting strip between the two spheres, followed by a cancellation
along $\cal C\,$.
\s


{\em 2.3.2) Reduction of Moves}\label{pg-232}
\s

It turns out that the complete set of relations between bridged links, as listed in
Corollary~\ref{co-compl}, is highly redundant. It is not hard to imagine
that there are many possibilities to select a minimal subset of moves from
it, which generates all other moves and equivalences.
 Let us assume that the ribbons passing through a general surgery sphere
penetrate at given intervals $I_1\ldots I_k\ldots\,$. 
Then the set of moves that has been proposed in Theorem~\ref{TWO} of the
introduction follows from the following reduction:
\medskip

\bpp\label{pp-quint-mv}

Two bridged link presentations yield the same three manifold if they are related 
by the following moves:
\ben
\item  {\rm Isotopies} Aside from general isotopies of surgery spheres and 
ribbons, which are constant along at the intervals $I_j$, we have isotopies
where ribbon pieces are pushed through surgery spheres of different
components. These can be written as combinations of

a) Opposite actions of the braid group of the two spheres.

b) An untangled, unbraided loop being pushed through. 

\item {\rm One Slides} We may assume that the $Hi$ -sphere is pushed
through $Lo$ first, followed by the attached ribbons. More precisely, before the
move we have ribbons $r_1\ldots r_m$ going through $I_1\ldots I_m$ of $Lo$
and ribbons $r_{m+1}\ldots r_{n+m}$ going through $I_1\ldots I_n$ of $Hi\,$.
After the move we have $r_1\ldots r_{n+m}$ going through $I_1\ldots
I_{n+m}$ of {\rm Lo}. Also, we only consider one slides where the two $Lo\,$
spheres are in different components. 

\item void

\item {\rm Cancellation} We only consider cancellation with no other
ribbons but the cancelling one attached. 

\item {\rm Modification} As stated.

\item {\rm Signature} As stated.

\item void

\item void

\een

\epp

{\em Proof:} An isotopy over a surgery sphere can be deformed into an isotopy
where the strands in a vicinity of the sphere are always radial, except for
passages through the stratum $p\fo_{12}\,$. The latter give rise to the
case {\em 1.b)}. The ambiguity of pushing the strands into the prescribed
positions $I_j\,$ are given by elements of the braid group as described in
{\em 1.a)} and $2\pi$-twists of the strands, which can be reexpressed as the
combination of two loops and one braid.

In order to verify 2)
we deform a one slide  such that first the $Hi$-sphere is pushed through
$Lo$, (the attached ribbons going through $Lo$ will then be loops at $Lo$,)
and then the ribbons at $Hi$. Conjugating everything with isotopies at $Lo$
we can confine ourselves to the situation where the ribbons are ordered as
described.  
 
\begin{figure}[ht]

\begin{center}\ 
\epsfxsize=5.7in
\epsfbox{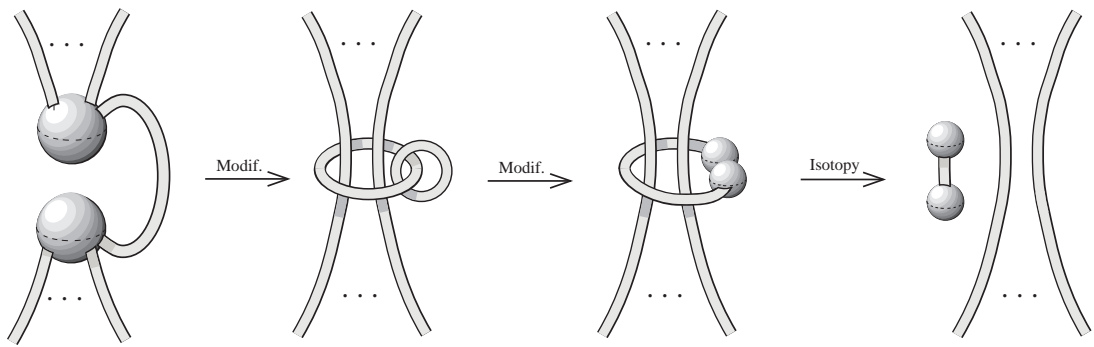}
\end{center}

\caption{Specialization of Cancellation}
\label{fig-canc-red}
\end{figure}

If we have the general situation of a cancellation as in
Figure~\ref{fig-canc} we may push the cancelling ribbon $c_2$ to the side
an use modifications as in Figure~\ref{fig-canc-red} to reduce the
cancellation to an isolated one. The reduction can also be found by using
two slides. Specifically,  we may separate one penetrating ribbon after the other
from the surgery spheres by sliding it over the extra cancelling ribbon.
   
\begin{figure}[ht]

\begin{center}\ 
\epsfxsize=5.7in
\epsfbox{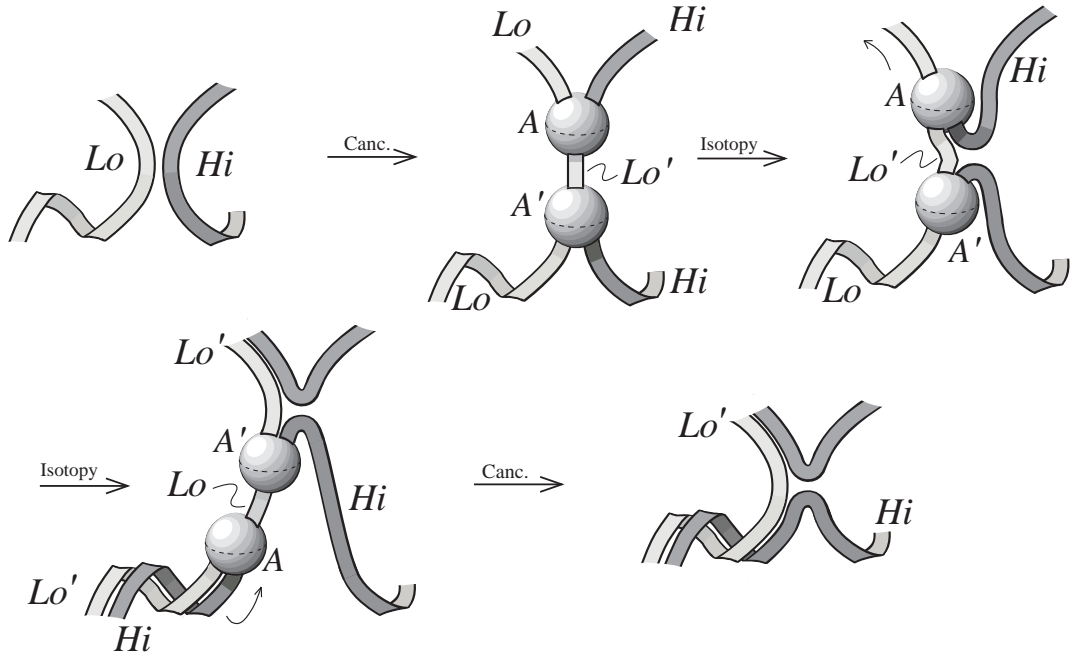}
\end{center}

\caption{Elimination of Two-Slides}
\label{fig-2slide-repl}
\end{figure}

The two-slides can be expressed in terms of cancellations. This is shown in
Figure~\ref{fig-2slide-repl} where we  move the surgery sphere $A$  along 
the $Lo$ ribbon of  the two-slide.
 
Move 6) is omitted if we are only interested in cobordisms of $\wct\,$.
It is easy to see that the removal of $\,^0\bigcirc\mkern-15mu\subset^{\cdots}_{\cdots}\,$
as in move 8) can also be given by a modification at the annulus followed by a
1-2-cancellation along the ribbon $\cal C\,$. 
Finally, if we have a pair of surgery  spheres in the same component we may 
move them close to each other and replace the moves {\em 1a), 1b)} and {\em
2)} by ordinary isotopies through an annulus, which are  conjugated by a modification.  
\hfill$\Box$
\s

The following, alternative  set of generating moves among ordinary links constitues the
proper generalization of Kirby's original calculus to non-simply connected
manifolds:

\bpp\label{pp-red1}
Suppose for two ordinary link presentations on a connected not necessarily
simply connected manifold $M_0\,$ the surgered manifolds are homeomorphic.

\noindent
Then we can relate the
presentations by a sequence of the following moves
\begin{enumerate}
\item Isotopies
\item ${\cal O}_2$-moves 
\item $\eta$-move
\item ${\cal O}_1$-moves\quad.
\end{enumerate}
If we consider only cobordisms in $\wct$ the ${\cal O}_1$-move is again
omitted.
\medskip

\noindent
If $M_0\,$ is also simply connected we may replace the $\eta$-move by the 
$^0\bigcirc\mkern-13mu\bigcirc^0$-move.

\epp

Proposition~\ref{pp-red1} can be directly derived from
Proposition~\ref{pp-quint-mv}, thus yielding a more structured proof of
Kirby's original theorem. The strategy is to consider the formal spaces of
bridged links $\bigl({\cal BL}\bigr)$ and ordinary links $\bigl({\cal
L}\bigr)$ in $M_0\,$ and define the quotient space $\overline{\bigl({\cal BL}\bigr)}$     
 and $\overline{\bigl({\cal L}\bigr)}$, where the relations from the two
previous Propositions have been imposed. We certainly have an inclusion
$\i\,:\,\bigl({\cal L}\bigr)\hookrightarrow \bigl({\cal BL}\bigr)\,$
and it is implied in the proof of Proposition~\ref{pp-quint-mv} that this
factors into a map $\bar{\i}\,:\,\overline{\bigl({\cal
L}\bigr)}\to\overline{\bigl({\cal BL}\bigr)}\,$. If we choose a path
between associated surgery spheres for each bridged link we can define a
map $p\,:\,\bigl({\cal BL}\bigr)\to \bigl({\cal L}\bigr)\,$ by moving
the spheres close to each other along the paths and elimnating them with 
a modification. What needs to be established is that $p\,$ factors into a
map $\bar p\,:\,\overline{\bigl({\cal BL}\bigr)}\to\overline{\bigl({\cal
L}\bigr)}\,$, which is independent of the choice of the recombination 
paths. 
It is then obvious that $\bar p$ and $\bar{\i}\,$ are inverses of each
other so that the two link-calculi are equivalent on connected manifolds.
A detailled discussion will be given in a separate paper, [Ke3].$\smile$
\medskip

\begin{figure}[ht]

\begin{center}\ 
\epsfbox{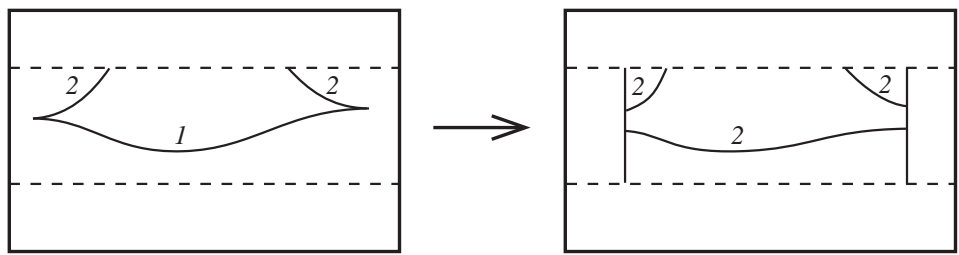}
\end{center}

\caption{Elimination of One-Handles}
\label{fig-graca}
\end{figure}

In this context the interpretation of the $\eta$-move is to compensate
for the recombination ambiguity of the isolated cancellation diagram.
For the same reason the $\eta$-move was introduced in [Ki], when
trajectories of one-handles that started and ended in different birth and
death points had to be eliminated. The graphic of the function is 
changed as indicated in Figure~\ref{fig-graca}, where the vertical lines
indicate modifications of the fourfold.   
If the recombination is, e.g., done
along the trajectories ending at the birth point, the first modification is
an addition of a $^0\bigcirc\mkern-13mu\bigcirc^0$-element, and the second
a subtraction of an $\eta$-configuration, where $\cal C$ runs along the
trajectories of the surgery spheres.     

\medskip

The fact that for simply connected manifolds $M_0\,$ a general $\eta$-move is the
composite of ${\cal O}_2$-moves and a
$^0\bigcirc\mkern-13mu\bigcirc^0$-move has been proven in [Ki]. A slightly
more direct argument starts from a general PL-homotopy of the curve $\cal
C\,$. This is made transversal, such that only at a finite number of times
there may be single crossings of $\cal C$ with another segment $\cal X$, which
can be a piece of $\cal C$ it self or another ribbon. If we move the
annulus aroud $\cal C$ to this crossing it can be substituted by a
two-slide of $\cal X$ over the annulus. For the change of framing of $\cal C$,
but also the reexpression of the $^0\bigcirc\mkern-13mu\bigcirc^0$-move,
using auxiliary $\FLO$'s and ${\cal O}_1$-moves see [Ki].
 
Finally, note that the reduction of the $\eta$-move can also be taken from
its five dimensional interpretation. In particular, it is clear that the
attaching curve $\cal C\,$ for five dimensional 2-handle can be contracted inside
of $M_0\times I\cup\hh^2\cup ... \cup\hh^2\,$.


\subparagraph{2.4) Bridged Ribbon Graphs}\label{pg-24}
\
\s

The construction of the invariants of [RT] from a link presentation is
preceded by a reduction of the presentation to spaces of so called
ribbon graphs, which are generic projections of the links into ${\bf R}^2\,$.
 
In this paragraph we shall describe the analogous space of {\em bridged} ribbon
graphs, which will contain not only ribbons but also special pairs of
coupons. The relations we will impose on this space are the ones that
appear  in Proposition~\ref{pp-quint-mv} and the usual  relations that take
care of isotopies and the ambiguity of choosing a projection.

\begin{figure}[ht]

\begin{center}\ 
\epsfbox{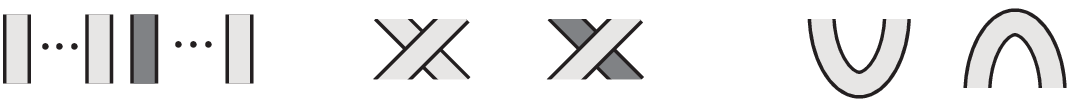}
\end{center}

\caption{Ribbon Segments}
\label{fig-cdr-elem}


\begin{center}\ 
\epsfbox{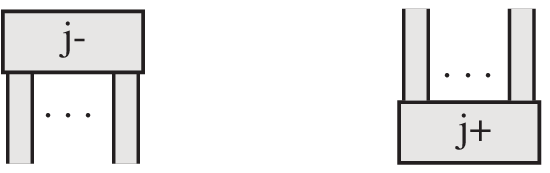}
\end{center}

\caption{Pairs of Coupons}
\label{fig-br-elem}
\end{figure}

We confine ourselves to  presentation in a union of $K$ handlebodies $H_g$
or spheres $S^3\,$. They are given by the union  
of $K$ diagrams of bridged  ribbon graphs in ${\bf R}^2\,$. 
Clearly, the interior of $H_g$ is homeomorphic to
${\bf R}^3$ with $g$ strands, extending to $\infty$, removed. 
Thus in a projection of a ribbon diagram into ${\bf R}^2$ we present the
strands as infinite, parallel strips in vertical direction. The pieces 
of an ordinary  directed ribbon graph as in [RT0] are then given by the
elements depicted  in Figure~\ref{fig-cdr-elem}.

Here the black strands are the ones representing the holes in $H_g$.
In a bridged ribbon graph we have in addition pairs of coupons, one with strands 
going up (labeled $(j,+)$) and one with the same number of strands going down.
The two coupons, as in Figure~\ref{fig-br-elem}, live  in different
${\bf R}^2$'s of the K-component  ribbon-graph, if they present surgery
spheres in different handlebodies.

We have the usual relations for
projections of ribbon graphs, e.g., $Rel_1-Rel_{10}$ of [RT0].
In these relations we also include auxiliary black strands, whenever they
make sense. 
Isotopies of the gluing sphere lead to the moves $Rel_{11}-Rel_{13}\,$.

The last relation also takes care of the ambiguity of identifying a surgery
sphere with a $+$ or a $-$ coupon.
It is  understood that the upside down versions of these pictures are
also included in the set of moves. 
\medskip

Next, we describe the $1a)$ and $1b)$ Isotopies, which are also used to
resolve the ambiguity of moving general strands into the standard positions
$I_1\ldots I_n\,$. The moves are given by the pictures for
$Rel_{14}-Rel_{15}$ and their mirror images.

The reduced one slide is given by  $Rel_{16}$ and reflections. The
modification and the  
reduced cancellation are given by $Rel_{17}-Rel_{18}$. We remind ourselves
that $Rel_{17}$ does not imply the relations $Rel_{14}$ and $Rel_{15}\,$, since
the coupons may be in different components.

Finally, we have the signature move $Rel_{19}$, which is the inclusion of an
isolated, unknotted ribbon $\FLO$ with one twist.  


\nopagebreak
\begin{figure}[ht]

\begin{center}\ 
\epsfysize=1.2in
\epsfbox{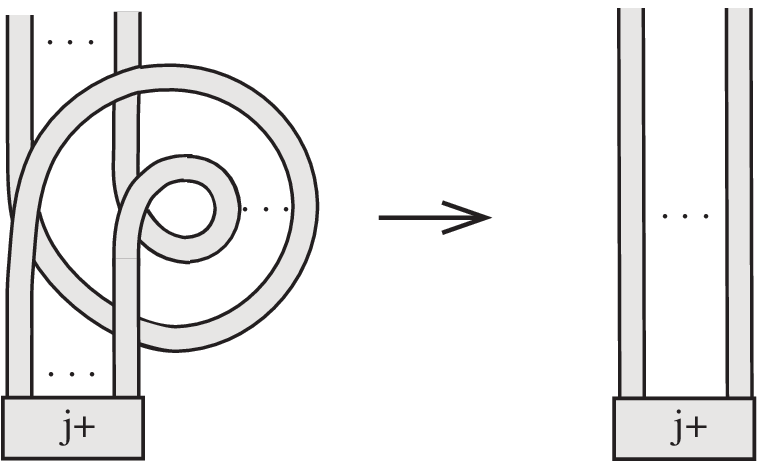}
\end{center}

\caption{$Rel_{11}$}
\label{fig-rel-11}
\end{figure}
\begin{figure}[ht]

\begin{center}\ 
\epsfysize=.9in
\epsfbox{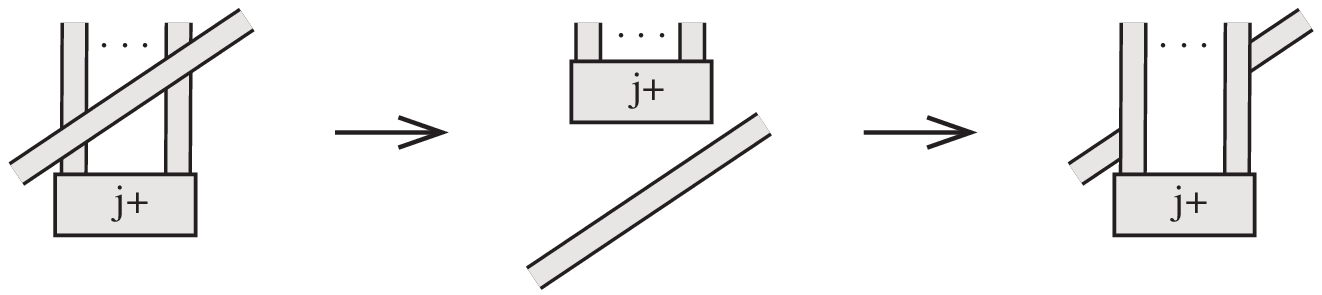}
\end{center}

\caption{$Rel_{12}$}
\label{fig-rel-12}
\end{figure}
\begin{figure}[hb]

\begin{center}\ 
\epsfxsize=5.2in
\epsfbox{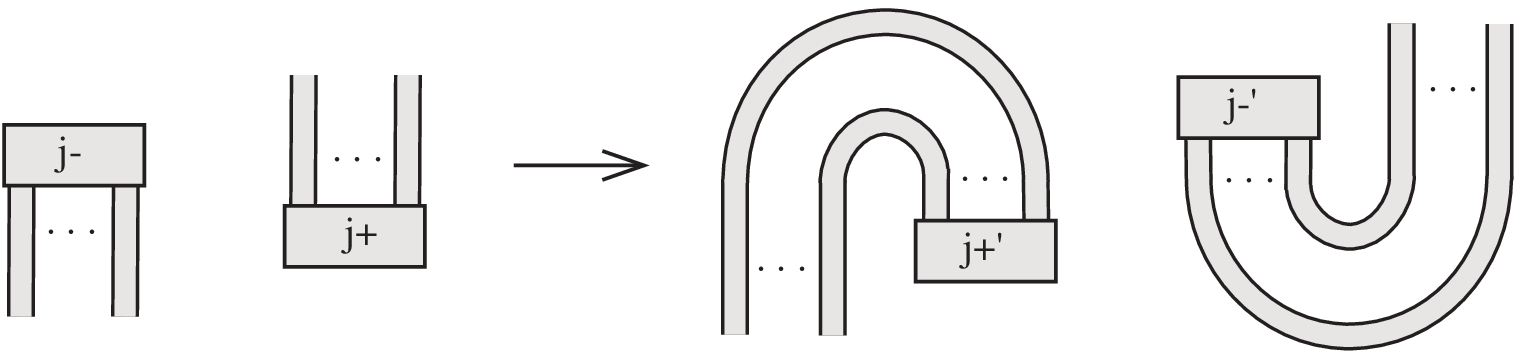}
\end{center}

\caption{$Rel_{13}$}
\label{fig-rel-13}
\end{figure}
\begin{figure}[ht]

\begin{center}\ 
\epsfbox{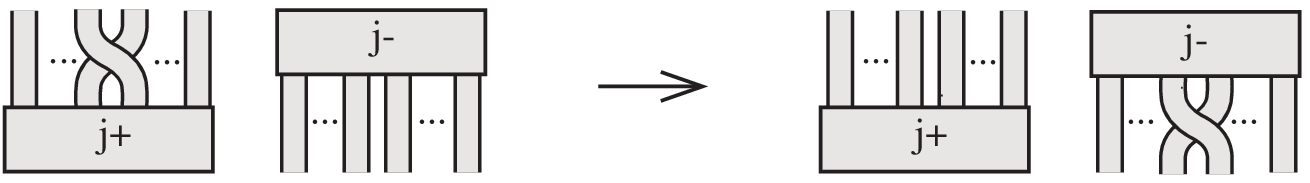}
\end{center}

\caption{$Rel_{14}$}
\label{fig-rel-14}
\end{figure}

\begin{figure}[hb]

\begin{center}\ 
\epsfbox{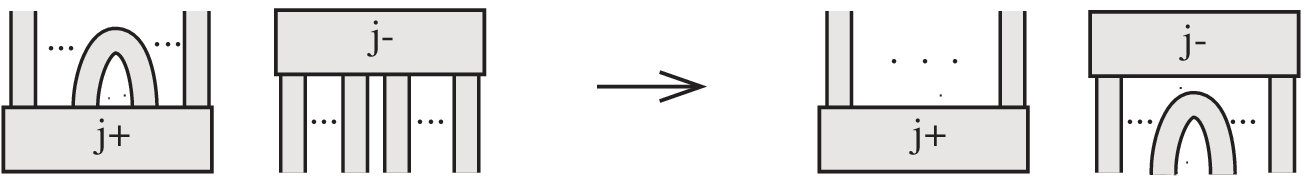}
\end{center}

\caption{$Rel_{15}$}
\label{fig-rel-15}
\end{figure}
\begin{figure}[ht]

\begin{center}\ 
\epsfbox{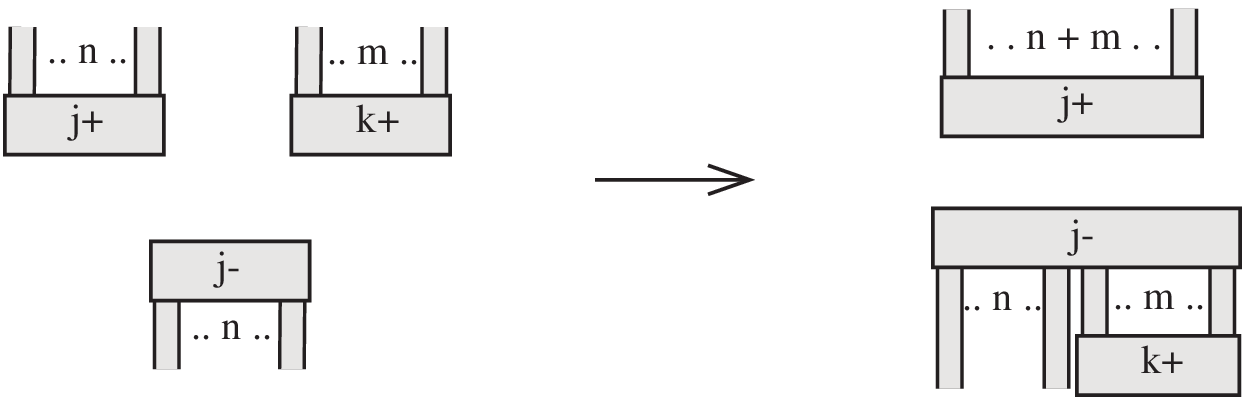}
\end{center}

\caption{$Rel_{16}$}
\label{fig-rel-16}
\end{figure}

\begin{figure}[ht]

\begin{center}\ 
\epsfbox{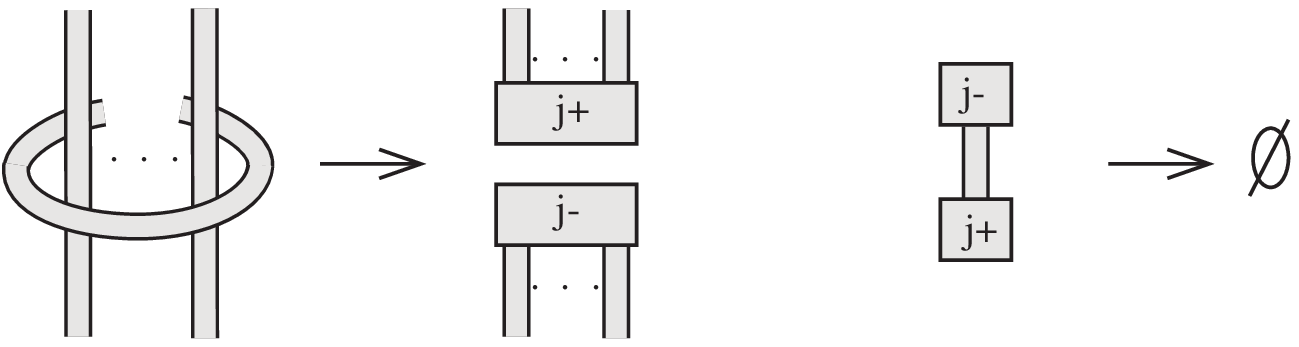}
\end{center}

\caption{$Rel_{17}\qquad
Rel_{18}$}
\label{fig-rel-17} 
\end{figure}

\parbox{6in}{
\bpp \label{pp-BL-moves}
\
\nopagebreak
\ben
\item Two directed, bridged, ribbon graphs are projections of 
bridged link presentations of the same manifold if and only if they are
related by relations $Rel_1-Rel_{19}$.
\item If the coupons $(j,\pm)$ lie in the same component of the presentation the
moves $Rel_{13}, \,Rel_{14}, \,Rel_{15}\,$, and $ Rel_{16}$ follow from the
other moves. 
In particular, if the presentation has only one component
$Rel_1-Rel_{12}$ and $Rel_{17}-Rel_{19}\,$ form a complete set of moves
for presentations of  cobordisms in $\ct$.
\item For  presentations of {\em connected} morphisms in $\wct$ we can omit
$Rel_{19}\,$.   
\een  
\epp}

%% file: BL3.tex
 \section*{3) Standard and Tangle Presentations of Cobordisms}\label{pg-3}
\

\s

The results in Chapter {\em 2} on bridged links are valid for presentations
on any compact, oriented three manifold. In this chapter we wish to
consider the special case of  surgery presentation on the three fold $\H$
discused in Chapter {\em 1}. The boundaries and the charts of the presented
cobordisms are always understood to be the canonical ones of $\H\,$.
We shall use the special form of $\H\,$ to simplify the presentation even
further. 

To this end we first introduce ``standard'' bridged link presentations for
which $\breve M$ and the bridged link diagram in the handlebodies is
normalized. The relevant presentation is thus in $S^3$ from where it can be
reduced  to tangles in the spirit of [RT]. 

A larger part of the discussion is devoted to the subtleties of the
composition rules.  


\subparagraph{3.1) Standard Presentations on $S^3$}\label{pg-31} 
\
\s

\begin{figure}[ht]
\begin{center}
\epsfxsize=5.2in
{ \ \epsfbox{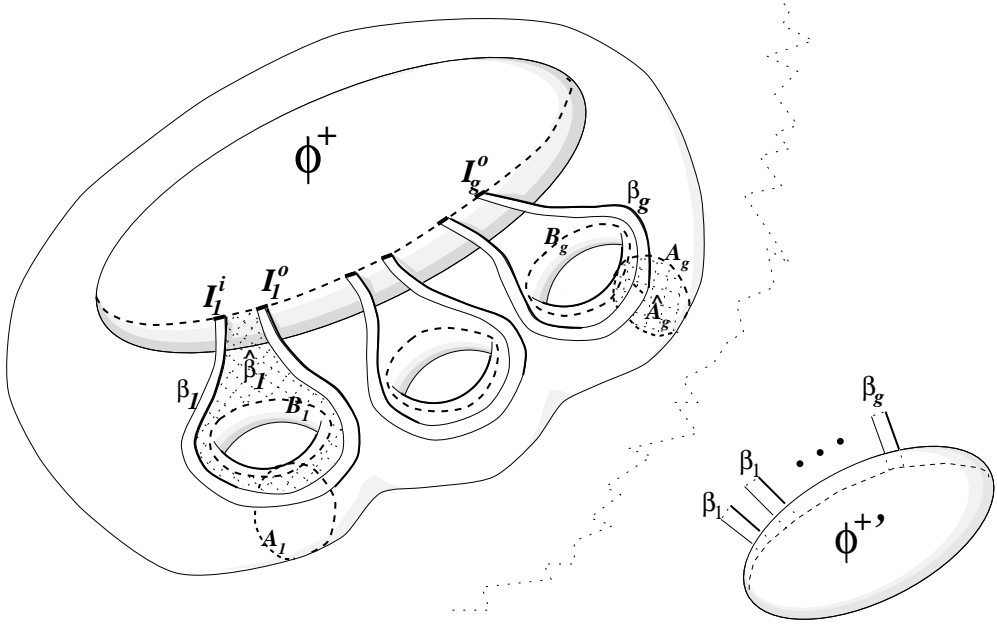}}
\end{center}
\caption{Standard Link in $H^+_g\,$}
\label{fig-stan-+}
\end{figure}

\begin{figure}[ht]
\begin{center}
\epsfxsize=5in
{\ \epsfbox{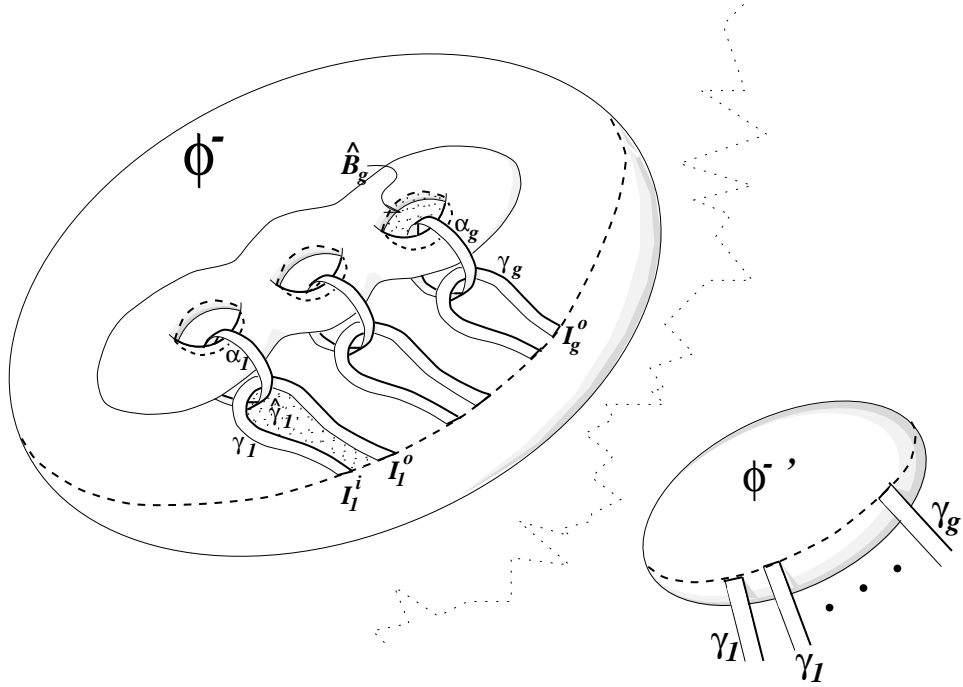}}
\end{center}

\caption{Standard Link in $H^-_g\,$}
\label{fig-stan--}
\end{figure}

In this paragraph we introduce a {\em standard} presentation of a connected 
cobordism, and we show that every cobordism admits one such presentation.
We begin with the surgery descriptions of $\breve M$ on $\H$ with only one
extra $S^3$ and  
the standard link presentation in the $H^{\pm}_g$ components. The latter
are indicated in
Figures~(\ref{fig-stan-+}) and~(\ref{fig-stan--}). 
\s


{\em 3.1.1) Standard Presentation of $\breve M\,$ and the class $\cal
U$}\label{pg-311} 
\s

The one handle structure of a presentation can be normalized as prescribed 
in the following Lemma:

\blm\label{lem-one}  
Suppose $\h^{-1}(1)$, and $W$ (as in Chapter 1) are connected. Then after 
connecting a sufficient number of modifications on $W\,$ we can reduce the
number of one handles in the decomposition of $W\,$ to the minimum 
$|\pi_0(\h^{-1}(0))|-1$.

If $\h^{-1}(0)\,=\,\H(\bar g^+,1,\bar g^-)\,$, and if the extra sphere $S^3$ is
labeled by $\beta\,$, then we can assume that the index one singularities
have labels $(\lambda,\beta)$ where $\lambda$ runs over all other
components. 
\medskip

We denote by $\cal U$ the isotopy classes of such presentations with the
surgery spheres in certain preassigned positions. 
\elm

\s

{\em Proof }: It is always possible to find a Morse function in 
$\fo$ such that all 
critical points are in general position.  
Using the Principle of Independent Trajectories we may deform $f$ such 
that the critical points of index $i$ have
values in  $\,]b_i,b_{i+1}[\,$ , for given numbers $0=b_0<b_1<...<b_5=1\,$.

By the following process we may reduce the index one points to the minimal
number.
Let us consider the first  critical value $c\in\,]b_1,b_2[\,$ with label
$(\alpha,\beta)\,$.  Now, if $\alpha=\beta$ we may change the index of the
singularity from one to two using a modification and push it above $b_2$. 
More generally, if $x\in\,]b_1,b_2[\,$ is a critical value we 
consider the graph, whose vertices are the components of $\H$ and which are
connected by an edge whenever there is a critical value in $]b_1,x[\,$ with
the corresponding label. If the label of $x$ introduces an edge at the
component
of the given graph we may replace it as above, since the descending
trajectories end in the same component  of $f^{-1}(x-\epsilon )\,$.

We end up with a  label graph with a minimal number of edges. That is, if we
remove any edge from the graph the component it lies in will  decomposes into two
components.

The components of $f^{-1}(b_1)\,$ consist of the $D^4$'s, that are created
by the index zero points, and the components of $\h^{-1}(0)\,$. If a $D^4$
appears in a label of a critical point, we know by minimality that the the
other component of the label is different from this $D^4$. Hence we may
apply Smales Lemma and remove the $D^4$-component by  a 0-1-cancellation.

We arrive at the  situation where the labels contain only components of
$\h^{-1}(0)\,$. Hence, $\breve M$ is the direct union of $D^4$ 's and 
connected sums of components of $\h^{-1}(0)\,$.

Replacing $f$ by $1-f$ we use the same
argument to get rid of  singularities of index three, since $g^{-1}(1)\,$ is
connected.

Now, the remaining singularities in $f^{-1}([b_2,b_4])$ do not change the number
of components of the level sets. Hence the singularities of index 0 and 4 
come in pairs, 
which belong to  additional components of $W$. However we assumed $W$ to be
connected. This implies the absence of index 0 and 4 points.

A one-slide corresponds to replacing the labels of the singularities as
indicated in the middle of Figure~\ref{fig-codone}.  
It is now easy to see that the form of the one-singularities described in 
part {\em 2)} can be obtained  by sliding the one handles across each
other, i.e., the label graph can be moved into a star form with $S^3$ in
its center. \hfill $\Box$

\s

{\em Remark :} To obtain a standard presentation we could have also started
with a connected manifold $\h^{-1}(0)$ right away. In fact, the cobordisms
$\,\bigl(\breve W=f^{-1}([0,b_2]),\,f\bigr)\,$ for  $f$'s with a minimal
number of index one singularities are all  
isomorphic.  $\breve M$ is just the connected sum of the common $S^3$ and
each  $H^{\pm}_g$.

Also, there always exists a path of Morse functions connecting two such 
presentations, which have the singularity structure as in
Lemma~\ref{lem-one} on $\breve W\,$. Certainly, there are also paths
of Morse functions for which the connecting one handles are slid across
each other. However, the description of these will be rather complicated since the 
space of one handle configurations is not simply connected. 

Still, it is convenient to view  the 
one handles surgeries to be given by a 
Morse function when we describe  compositions of cobordisms.

\s


{\em 3.1.2) Standard Links in $H^{\pm}_g\,$}\label{pg-312}\nopagebreak
\medskip

We define standard forms of a bridged links inside of the $H^{\pm}_g\,$.
They are depicted in Figures~\ref{fig-stan-+} and ~\ref{fig-stan--}.
\medskip

We start by drawing  Heegaard diagrams on $\Sigma_g\,$  that yield 
splittings of $S^3\,$. 
They consist of curves  $\{A_1\ldots A_g\}$ on $\partial H_g$ where the
$A_i$  are contractible to  the
inside of $H_g^+$, and a Heegaard diagram  $\{B_1\ldots B_g\}$ where each
$B_j$ is contractible to the outside $H^-_g$, and intersects the other diagram 
only once in the curve $A_j\,$.
\medskip

Moreover, we draw on the standard sphere that gives the one-handle
attachement
of $H^{\pm}_g\,$ a line $\,G^g\,$ which is disjoint from the equator 
and contains $\,2g\,$ intervals which we call in the order they are aligned
along $\,G^g\,$ 

\beq\label{eq-intervals}
I^i_1,\,I^o_1,\,I^i_2,\ldots\,,I^o_{g}\,\subset \,G^g\,
\eeq

In the components $H^-_g$ we thicken the curves $A_j$ to ribbons 
$\alpha_j$ in
the surface and push them slightly off into the interior of $H^-_g\,$. 
We depict the surgery ball of the one handle attachement as the complement
of a ball in $\,S^3\,$, which contains $\,H^+_g\,$ so that the surgery
sphere surrounds $\,\Sigma_g\,$. For each $\,\alpha_j\,$ we introduce
another ribbon $\,\gamma_j\,$, which starts at at the interval $\,I^i_j\,$
 and ends at the interval $\,I^o_j\,$. It shall follow a radial direction
away from $\,\Sigma_g\,$ and surround $\,\alpha_j\,$ close to $\,\Sigma_j\,$
as depicted in  Figure~\ref{fig-stan--}. the identification of the
$\,\phi^-\,$-spheres for the 1-surgery is as they appear in the picture,
i.e., the content of the ball on the left is inserted in the ball on the
right part of the figure. 
\medskip

In the $H^+_g$ components  we push the ribbons $\,\beta_j\,$ along the $B_j$ -
curves to 
the inside. Each $\,\beta_j\,$ connected to the surgery sphere $\,\phi^+\,$ as
indicated in Figure~\ref{fig-stan-+} such that it starts at the interval
$\,I^i_j\,$ and ends at the interval $\,I^o_j\,$.

Let us denote by $\cal S$ the isotopy classes of these standard
presentations.
\s


{{\em 3.1.3) The $\sigma$-Move and the Lemma of Connecting Annuli}\label{pg-313}
\s 

Clearly}, in a standard presentation the cobordism is entirely given 
by a bridged link diagram in $S^3$. Yet, in order to ensure the existence
of such presentations and in order to describe a complete set of moves 
as in
Proposition~\ref{pp-red1} we need to include a special version of the 
$\eta$-move combined with two slides and isotopies. We will call this
move a $\sigma$-move. In the context of a combinatorial tangle
presentation of the mapping class group  an analogous move was introduced 
in [MP], where it was called $K_3\,$. This move also turns out to be a
special case of  the Lemma of ``Connecting Annuli'', which we will
discuss in the end of this section. The exceptional cases of this lemma are points of
caution for the composition rule of presentations of general cobordisms.
\s

\noindent
The {\em $\sigma$- move at the j-th handle}  in $H^-_g$ is described as follows: 
\s

\begin{figure}[ht]
\begin{center}
\epsfxsize=5.5in
{ \ \epsfbox{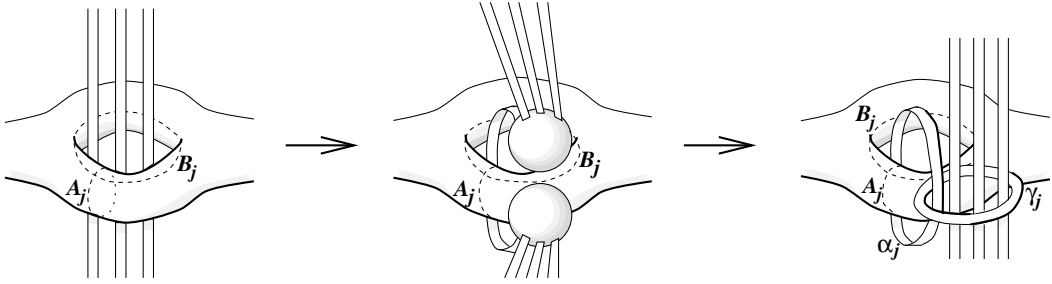}}
\end{center}
\caption{$\sigma$-move in $H^-_g\,$}
\label{fig-push-A}
\end{figure}

We introduce a disc $\hat B_j\,$ in $H^-_g\,$  that is bounded by the curve
$B_j\,$. 
We may assume that all ribbons pass through the disc $\hat B_j\,$ transversally. 
The next step is to undo a cancellation along the disc. We then move the
pair of  surgery spheres around the handle of $H_g$ with meridian $A_j\,$, 
stretching the cancelling ribbon such that it coincides with the ribbon
$\alpha_j\,$. 
The spheres are recombined by a modification as indicated in
Figures~\ref{fig-push-A}. 

One half of the annulus of the modification is
pushed through $\,\phi^-\,$ such that the other half is identical with the
$\,\gamma_j\,$-ribbon of the standard presentation. Furthermore, we push
the strands that were intersecting $\hat B_j\,$  into the
$S^3\,$-component, such
that the disc $\,\hat\gamma_j\,$ bounded by $\,\gamma_j\,$ and the interval
$\,I_j\,$ on $\,G^g\,$ between $\,I^i_j\,$ and $\,I^o_j\,$ is only intersected by
$\,\alpha_j\,$. 
\s

\begin{figure}[ht]
\begin{center}
{ \ \epsfbox{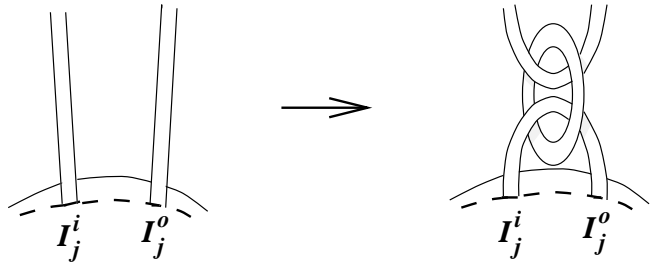}}
\end{center}
\caption{$\sigma$-move in $S^3\,$}
\label{fig-sigma}
\end{figure}

If the $\sigma$-move is applied to a standard presentation it may be
described as a move at the surgery sphere $\phi^-\subset S^3$ as indicated in
Figure~\ref{fig-sigma}. The annulus is what has been the
$\,\alpha_j\,$-band before the move, and the lower loop is the second half
of the newly created $\,\gamma_j\,$-ribbon.

We define the $\sigma$- move for $H^+_g\,$ in the same manner. 
Here the un-cancellation
is done at the discs $\hat A_j\,$ which are bounded by the $A_j\,$'s.  
Similarly, we push the modification annulus to the outside so that 
the annulus $\,\hat\beta_j\,$ bounded by $\,B_j\,$ and $\,\beta_j\cup
I_j\,$ is not intersected by any other ribbons. 

The
move in $S^3$ when  
applied to standard presentation is given  in exactly the same way as for
$H^-_g\,$. 
\medskip

The $\sigma$-move (e.g. in $H^-_g\,$ ) can also be described as an
$\eta$-move, where we introduce an $\alpha_j$-ribbon linked to the annulus, 
which will be extended to the $\,\gamma_j\,$-ribbon, and two slides of each 
ribbon passing through $\hat B_j\,$ over $\alpha_j\,$.

It is also a special case of the following Lemma for
``Connecting Annuli''. 
It gives rules for replacing  a piece of a link  as in Figure~\ref{fig-modi} with
only two strands passing through $R\,$, such that we are left with a link
with fewer components.

\begin{figure}[ht]
\begin{center}
\epsfxsize=5.2in
{ \ \epsfbox{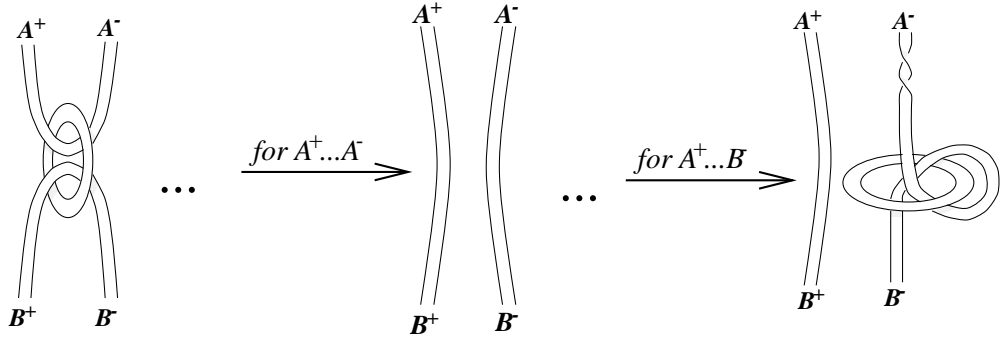}}
\end{center}
\caption{Connecting Annulus}
\label{fig-conn-lm}

\end{figure}

\blm\label{lem-con}
Suppose a disc bounded by a surgery ribbon is penetrated by two pieces of a
ribbon diagram as indicated on the left of Figure~\ref{fig-conn-lm}. We
assume that the ribbons are oriented and that the 
orientation is defined by which face is upward in the plat graph. 
Moreover, we distinguish the cases where, if we follow the ribbon connected
to $A^+$ in the  diagram, we will return to the annulus at the point
$A^-\,$, $B^-\,$, or $B^+\,$. 

\ben
\item[$A^+...A^-\,$] The diagram may be substituted by two straight ribbons
joining  $A^{\pm}$ to  $B^{\pm}$ respectively, as indicated in the middle of
Figure~\ref{fig-conn-lm}. 
\item[$A^+...B^-\,$]The diagram may be substituted by one straight ribbon
joining  $A^+$ to  $B^+$ and one tangle joining $A^-$ to $B^-$, as
indicated in the right of  Figure~\ref{fig-conn-lm}. 
\item[$A^+...B^+\,$] The ribbon piece running from $A^+$ to $B^+$ can be
replaced by a ribbon for which the corresponding closed ribbons (as in the
middle picture of \ref{fig-conn-lm}) are isotopic. 
\een
\elm

{\em Proof :} In the case $A^+...A^-\,$ we can slide the ribbon with labels
$A$ over the ribbon with labels $B$. Concluding with an $\eta$-move this
give us the move labeled $i)$ of Figure \ref{fig-conn-lm}.

\begin{figure}[ht]
\begin{center}
\epsfxsize=5in
{ \ \epsfbox{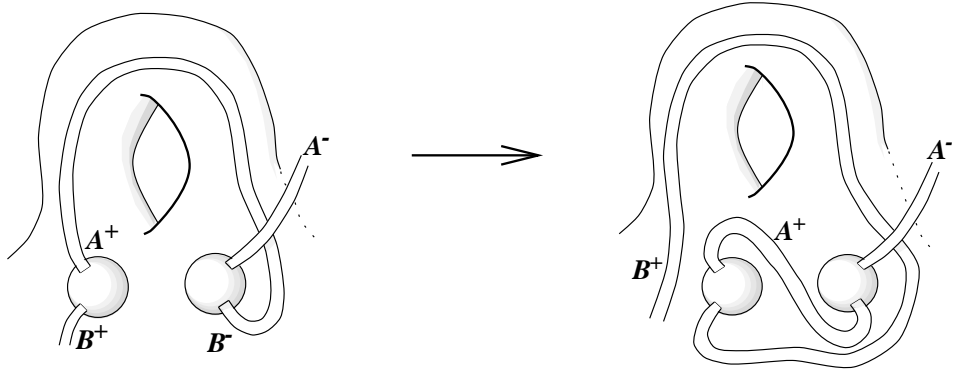}}
\end{center}
\caption{Recombination of Annulus}
\label{fig-conn-pff}
\end{figure}

For the case $A^+...B^-\,$ we start by adding to the
diagram the handle in which the $A^+...B^-\,$ lives and do a modification.
The surgery sphere connected to $A^+$ and $B^+$ is then pushed through the
handle as described in Figure~\ref{fig-conn-pff}, so that it arrives in the correct
position at the other sphere. Undoing the modification we obtain the
desired move $ii)$ in Figure~\ref{fig-conn-pff}. 

For the last case we remark that after a modification we
obtain two components of a  bridged link, to which we can apply isotopies
independently. They are different from those with fixed end points
$A^{\pm},B^{\pm}\,$ if $\pi_1(M)\,$ is non abelian. 
\hfill$\Box\,$   

\s


{\em 3.1.4) Existence of Standard Presentations and a Projection on $\,\cal
S$}\label{pg-314}\nopagebreak 
\medskip

It is clear from Lemma~\ref{lem-one} that we have an inclusion  
\beq\label{incl-pres}
\iota:\,{\cal S}\,\hookrightarrow\,{\cal U}\;.
\eeq
 
In this section we wish to show that $\iota$ is onto, i.e., induces an
isomorphism, if we  mod out the moves for link diagrams. More specifically,
we shall construct a standard representative of the inverse of $\iota$ on
link diagrams. 
 
\bpp\label{pp-ex-pro}
Every cobordism admits a standard presentation. More precisely:

There is a canonical map 
$$
\Psi\,:\,{\cal U}\to{\cal S}
$$ 
that sends a general bridged link  presentation in $\H$
to a standard presentation of the same cobordism in $\wct$.

The composition 
$$
\Psi\circ\iota\,:\,{\cal S}\,\to\,{\cal S}
$$
is the composition of $\sigma$-moves at the surgery spheres in $S^3\,$,
one for each handle.
\epp 

In place of a  {\em proof} we give the definition of $\Psi\,$:
\s

The first piece of $\Psi\,$ is one  $\sigma$-move at each handle of each 
$H^{\pm}_g\,$ - component. 

As a result of this operation we obtain a bridged link which looks
inside of a handlebody $\,H^{\pm}_g\,$ as in Figures~\ref{fig-stan-+}
or~\ref{fig-stan--} with additional components $\,\li''\,$. By prescription
of the $\sigma$-moves $\,\li''\,$ is disjoint from the discs ( or annuli )
$\,\hat\gamma_j\,$, $\,\hat\beta_j\,$, $\,\hat A_j\,$, and $\,\hat B_j\,$
that have been defined in Section 3.1.1). In fact, for sufficiently small,
closed
collars $\,V^{\pm}_g\subset H^{\pm}_g\,$ of the surfaces $\,\Sigma_g\cup\hat
A_1\cup\ldots\cup \hat A_g\cup\hat\beta_1\cup\ldots\cup\hat\beta_g\,$  in
$H^+_g\,$  and 
$\,\Sigma_g\cup\hat B_1\cup\ldots\cup\hat
B_g\cup\hat\gamma_1\cup\ldots\cup\hat\gamma_g\,$  
(containing all the $\,\alpha_j\,$) we may assume that $\,\li''\subset
H^{\pm}_g-V^{\pm}_g\,\cong\,\bigl(S^2-\amalg_{j=1}^g
I_j\bigr)\times ]0,1]\,$. 

Moreover, we may choose a collar  $\,S^2\times[1,1+\epsilon]\,$ of
$\phi^{\pm}\,$ in the $S^3\,$-component 
such that the links different from the $\,\beta_j\,$ and $\,\gamma_j\,$ are
all  in $K^{\pm}\,\cong\,\bigl(S^2-\amalg_{j=1}^g
I_j\bigr)\times ]1,1+\epsilon]\,$. From a homeomorphism
$]0,1+\epsilon]\,\isto\,]1,1+\epsilon]\,$ we obtain a homeomorphism
$\,\delta\,:\,\bigl(H^{\pm}_g-V^{\pm}_g\bigr)\cup K^{\pm}_g\,\isto\,
K^{\pm}_g\,$. It is clear that two links for which  the parts $\li''\,$
are replaced by $\delta(\li'')\,$ are isotopically equivalent.

Clearly, the result has the form of  a standard presentation and we may
give the  definition of $\Psi\,$ as the composition of the $\sigma$-moves
and the maps $\,\delta\,$. \hfill$\Box$
  
\bigskip
\ 

\bigskip


{\em 3.1.5) Moves in a Standard Presentation}\label{pg-315}
\s

What remains to describe the class of connected cobordisms by bridged link
diagrams in $S^3$ is a  set of moves. The assertion of the next proposition 
is that the $\sigma$-move is the only additional move for standard
presentations besides the proper specializations of the Kirby moves to the
$S^3$ component:

\bpp\label{pfff}

 Suppose we have two standard presentations on $\H\,$ of the same
connected cobordism in $\ct$.
Then they are related  by a sequence of the five moves 
\ben
\item Isotopies in $S^3$ (as in 1) of Section 2.2.2)) with fixed surgery
spheres.
\item Two slides (as in 3) of Section 2.2.2)) where the $Hi$-ribbon may be
attached to a surgery sphere. 
\item The signature move as in 6) of 2.2.2).
\item  The $\sigma$-move for each handle. 
\een
In the second move we may confine ourselves to $\kappa$-moves.
We may replace the third move by the
$\,^0\bigcirc\mkern-15mu\bigcirc^0\,$-move, 
if we only wish to preserve the class of the cobordism in $\wct\,$.  
\epp

{\em Proof :} Our task is to show that the set of moves listed in
Proposition~\ref{pp-red1} applied to presentations in $\cal S$ can be
expressed by the above moves. The basic idea is to conjugate other
sequences of moves continuously with the projection $\Psi$.

An isotopy of a presentation in $\cal U\,$ can be decomposed into isotopies
where a singularity of the ribbon diagram is moved through a specific
disc, $\hat B_j\,$ (or $\hat A_j\,$), and isotopies that are constant 
across all $\hat
B_j\,$. If we conjugate the first type of move with the map $\Psi$ we may
push the entire isotopy into $S^3$, where it is given by a passage of a
singularity through the respective loop (or annulus) from Figure~\ref{fig-sigma}. 

In the same way we can express a two slide and an $\eta$-move by  a two slide
and an $\eta$-move which are constant close to  the discs $\hat B_j$ conjugated
by $\Psi\,$. As in the proof of Proposition~\ref{pp-ex-pro} we may also
push the rest of a move to the outside of the handlebody.
Thus, a complete set of moves of $\cal U$ is given by $\Psi^{\pm 1}$ and the
types of moves from Proposition~\ref{pp-red1}, which are constant on 
the handlebodies. 

Since  $S^3-\bigcup\phi\,$
is simply connected we may use the arguments in the proof of
 Proposition~\ref{pp-red1} and replace the
$\eta$-move by two slides and $\,^0\bigcirc\mkern-15mu\bigcirc^0\,$-moves. 
This yields the presentation of $\wct\,$.

By the same Proposition, we can replace the
$\,^0\bigcirc\mkern-15mu\bigcirc^0\,$-moves by ${\cal O}_1\,$-moves if we 
consider presentations of $\ct\,$.  
The reduction of two handle slides to $\kappa$-moves
follows similarly from [FR].  

Since we started from presentations in $\cal S\,$, and since
$\Psi$ is in this case a  combination of the $\sigma$-move in $S^3\,$, we
conclude that the list given in Proposition~\ref{pfff} is complete.

\hfill$\Box$
\s


\subparagraph{3.2) Tangle Presentation of Cobordisms}\label{pg-32}
\
\medskip

As in  [RT] or Section 2.4) we seek to describe in this section a planar 
presentation of cobordisms. We shall use the reduction to standard presentations 
in $S^3$ to derive via a suitable projection a planar presentation of 
cobordisms in terms of ribbon  graphs in ${\bf R}\times [0,1]\,$, where the 
ribbons are allowed to end in the boundary of the strip. Thus the cobordisms
are presented by {\em admissible tangles}, i.e., configurations 
of such ribbons, which  fulfill a certain orientability condition.

The formulation of a complete set of moves is easily given using results 
from Section 2.4) and 3.1.5)

The composition rules contain a few subtleties: We shall give a Cerf theoretic
derivation of a decoration rule for boundaries with many components. The rule has
also been stated in a more rudimentary form in [Tu]. In our construction of 
presentations we also need to 
include a second rule related to the $\sigma$-move, which applies also to
connected boundaries.   In particular, we discuss the obstructions given
by the Lemma of Connecting Annuli to na\"\i ve compositions of tangle
diagrams.   
\medskip


{\em 3.2.1) From Standard Presentations to Admissible Tangles}\label{pg-321}
\medskip

The choice of the projection of a standard bridged link in $S^3$ to an
admissible tangle 
depends on  a few more conventions regarding the positions of the surgery spheres and
the links in a fixed $S^3$.

\begin{figure}[ht]
\begin{center}
{ \ \epsfbox{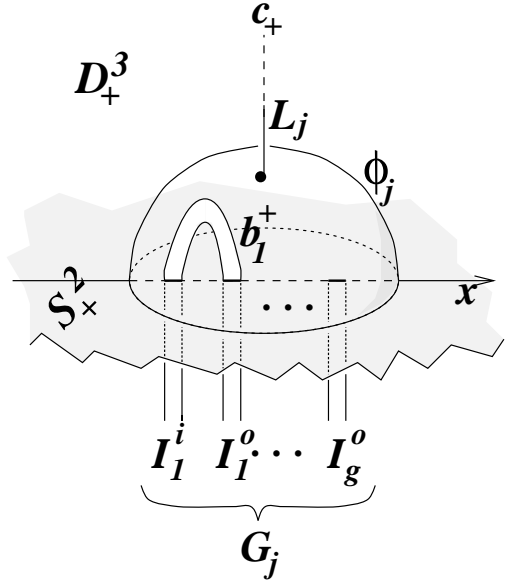}}
\end{center}
\caption{Alignment of Surgery Spheres in $\,S^3\,$}\label{fig-isto-II's}
\end{figure}

To start with we fix two spheres $S^2_{\pm}\,\subset S^3\,$, which separate
$S^3$ into three pieces homeomorphic to $[-1,1]\times S^2\,$ and
$D^3_{\pm}\,$. Inside the standard $S^2$ we fix a point $\infty$ and a
homeomorphism $S^2-\infty\isto {\bf R}^2\,$.  For any pair $\bar g_{\pm}$
we  fix a standard alignment of the $K_++K_-$ surgery spheres of a
standard presentation along the respective copy of ${\bf
R}^2$. More precisely, we shall fix a sequence of intervals $G_{j}^+,
j=1,\ldots,K_+$ on the $x$-axis of ${\bf R}^2\cong S^2_+-\infty\,$.
The standard position of the $j$-th surgery sphere $\phi_j$ is then specified by the
property that $\phi_{j}^+$ lies in the closure of $D^3_+$, 
such that the special line on $\,\phi_j^+\,$, containing the intervals
$\,I^{i/o}_s,\,s=1,...,g_j\,$, coincides with $\,G_j\,$. Also we require that the 
order in which the intervals appear on the $x$-axis is the same as
in~(\ref{eq-intervals}), see Figure~\ref{fig-isto-II's}.

Similarly, we define positions for the opposite spheres on $S^2_-$.

As in the remark of Section 3.1.1) we shall consider only  Morse function 
$f\big |_{\breve W}\,$ of
the surgery presentation, such that the surgery-spheres are always in the
defined positions, i.e, we will consider presentations on the 
manifold where the handlebodies are already connected to the $\,S^3\,$ by 
one-surgeries.
\medskip

Next we introduce interior points $c_{\pm}\in D^3_{\pm}\,$ and (unbraided) lines
$L^{\pm}_j$ in $D^3_{\pm}\,$ connecting $c_{\pm}$ to the surgery spheres $\phi^{\pm}_
j\,$ . The complement, $Q$, of the surgery balls, the points $c_{\pm}$,
and the lines $L_{j}^{\pm}\,$ is clearly $\cong S^2\times]0,1[\,$. 

Moreover, we can
introduce a homotopy of
embeddings $f_t:Q\hookrightarrow Q$, such that $f_t$  is the identity along 
$\phi^{\pm}_j\cap S^2_{\pm}\,$, $f_0=id\,$, and $f_1$ maps $Q$ onto the 
$[-1,1]\times S^2\,$-piece of $S^3\,$. Thus, for any  ribbon diagram of a standard
presentation and isotopy thereof with all ribbons in
$Q-[-1,1]\times\infty\,$, we have a canonical deformation to  equivalent
diagrams and isotopies in the $[-1,1]\times {\bf R}^2\,$. 
\medskip

\begin{figure}[ht]
\epsfxsize=4in
\begin{center}
{ \ \epsfbox{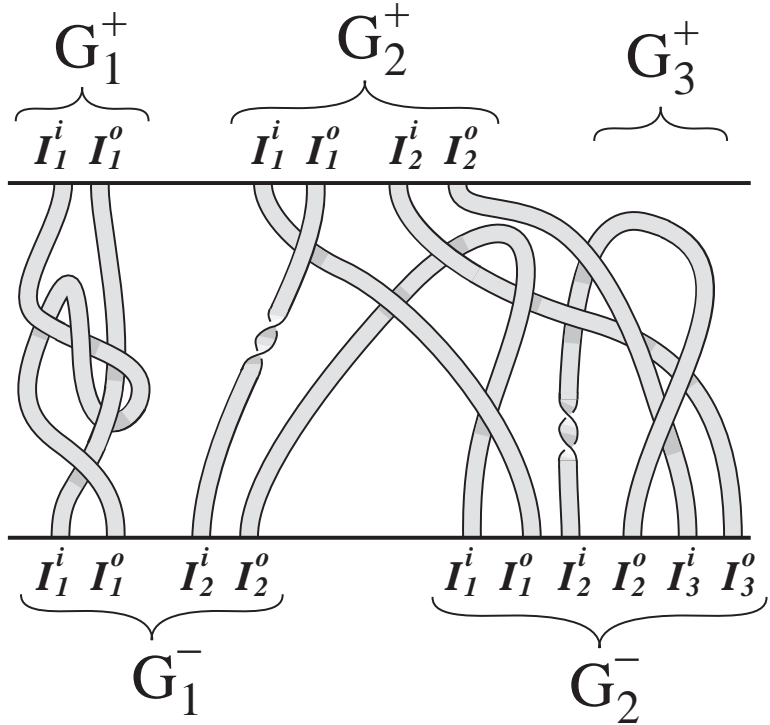}}
\end{center}

\caption{Plat Ribbon-Tangle}\label{fig-ex-admi}
\end{figure}

The diagram in $[-1,1]\times{\bf R}^2\,$ is then projected into a strip
$[-1,1]\times X\,$ where $X\cong{\bf R}\subset{\bf R}^2\,$ is the $x$-axis,
so that the $G_j$ are arranged along the lines $X_{\pm}:=\{\pm 1\}\times X\,$.
If the diagram is in a general position the projection yields a ribbon
graph (in the sense of [RT0]). It has the properties that to any of the
$2\sum_j g^{\pm}_j$ intervals of $ X_{\pm}\,$ a ribbon is  attached,
and, furthermore, the closed ribbons that result by inserting the strips
$b_k^{\pm}\,$, as 
indicated in Figure~\ref{fig-isto-II's},  are all orientable (i.e., $\cong
I\times S^1\,$). We call a ribbon graph in $[-1.1]\times X$ with these
properties an {\em admissible tangle }. An example with $\bar g^+=(1,2,0)$ and 
$\bar g^-=(2,3)\,$ is shown in Figure~\ref{fig-ex-admi}. Inserting the
bands $b_j^{\pm}$ yields four closed, orientable  ribbons.
\medskip

\begin{figure}[ht]

\begin{center}
\epsfxsize=4.2in
{ \ \epsfbox{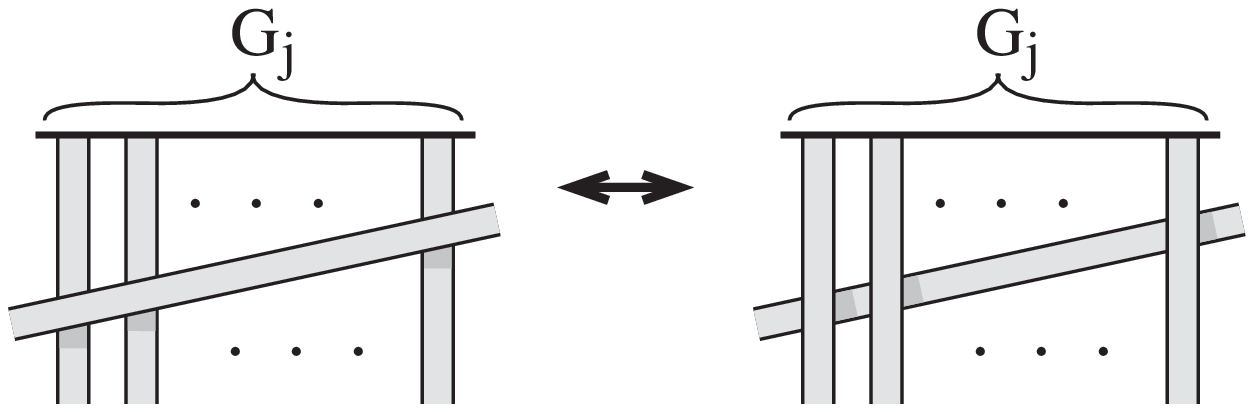}}
\end{center}
\caption{$\tau$-move}
\label{fig-tau-move}
\end{figure}

We can always deform  a given standard presentation into a position where
all ribbons lie in $Q$. Yet, a general isotopy of standard diagrams can be
merely chosen transversal to the complement of $Q$. I.e., the possible
singularities occur when an individual ribbon passes transversality though a
line $L_j$ or through $[-1,1]\times\infty\,$. For a line $L_j$ 
this gives rise to the 
additional {\em $\tau$-move at the group $G_j\,$} 
in the set of admissible tangles, which is given by
Figure~\ref{fig-tau-move} and its reflections. Here, a strand is moved
through the ribbons emerging from a group $G_j$  of $2g_j\,$ intervals on
$X$.  If we move a ribbon through $]-1,1[\times\infty\,$ the corresponding
move of admissible tangles is to push a ribbon through the strands of all
groups at once, and can thus be written as a combination of $\tau$-moves.
\medskip

\s


{\em 3.2.2) Moves for Admissible Tangles}\nopagebreak\label{pg-322}
\medskip

The possible moves for standard presentations are given in
Proposition~\ref{pfff}. It is straight forward  to derive from this
an equivalent  set of moves for the admissible tangles. A general isotopy of a 
ribbon diagram can always be written as a composition of the moves in [RT0]
and the $\tau$-move. 
As usual for presentations of connected morphisms in $\wct\,$ we may forget the
signature move 6). 
Since we give our presentation on a connected and simply
connected manifold we can also use [FR] to replace the two slides by the
$\kappa$-move, or 
by the remarks of Section 2.4) to use $Rel_{17}$ and $Rel_{18}$ for an
enhanced  bridged ribbon presentation. 
Given the projection rule of Figure~\ref{fig-isto-II's}, the form of the
$\sigma$-move in the plat tangle presentation is obvious from the one given
in Figure~\ref{fig-sigma}. 
Combining these observations with the results from Proposition~\ref{pp-BL-moves}   and
Proposition~\ref{pfff},
 we may now formulate the presentation  of cobordisms in terms of
admissible tangles:

\bpp
There is an isomorphism between  morphisms in $\wct$ and
the set of admissible, planar, ribbon tangles, \mbox{$\widetilde{{\cal
T}\!\!g\,}$}, divided by
the following relations
\begin{enumerate}
\item Isotopies: $Rel_1-Rel_{10}\,$, and $\tau$-move, 
\item $\sigma$-move, 
\end{enumerate}
and, alternatively, for links in $S^3$
\begin{enumerate}
\item[3.] $\kappa$-move,
\item[4.]$\,^0\bigcirc\mkern-15mu\bigcirc^0\,$-move,  
\end{enumerate}
or, for bridged links in $S^3\,$,
\begin{enumerate}
\item[3.] Isotopies $Rel_{11}\,$  and $ Rel_{12}\,$, 
\item[4.] Modification $Rel_{17}\,$, 
\item[5.] Cancellation $Rel_{18}\,$.  
\end{enumerate}
\epp
\s

{\em 3.2.3) Compositions of Admissible Tangles}\label{pg-323}
\medskip

The composition of two cobordisms can be presented  by a tangle that is
build up from the original tangle diagrams. In a more fancy language this 
means that we can endow the set of  admissible tangles with a composition
structure of a category and extend the presentation to a functor.

As outlined in the end of Section 1.3) it is sufficient to give the rules for the
compositions of tangles corresponding to products of cobordisms  as  in equations 
(\ref{perm-comp}) and (\ref{elem-comp}). We start with the first type, which
is easier.
\medskip

\begin{figure}[ht]

\begin{center}
\epsfxsize=5in
{ \ \epsfbox{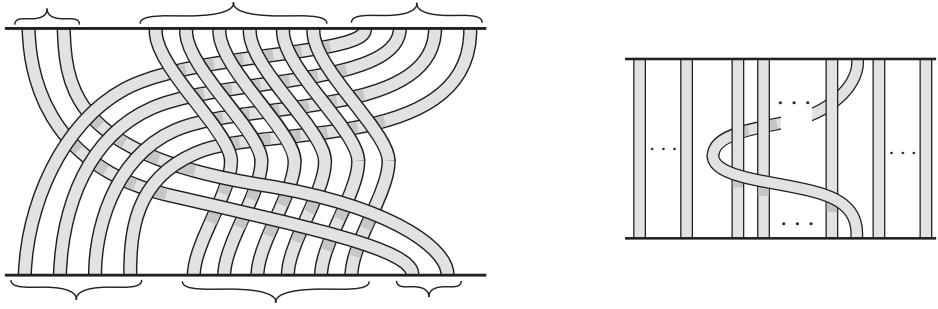}}
\end{center}
\caption{Permutation of Components,\qquad\quad Pure Braid\qquad\qquad}\label{fig-perm-gg}
\end{figure}

If we consider a standard presentation with surgery spheres $\phi^+_j$ in
positions  $G_1,\ldots, G_K\,$ we may use an isotopy to bring the spheres
into positions $G_{\pi(1)\,},\ldots\,,G_{\pi(K)\,}\,$ for  a given
permutation $\pi\in S_K\,$, so that they still describe the same cobordism.
The effect of the composition in (\ref{perm-comp}) is to permute the handlebodies
 $H_g$, so that  the bridged link with the $\phi_j$'s moved into new positions is
in fact a standard presentation of the composite. On the level of tangle
diagrams the isotopy will be given by a braid $b\in B_K\,$  of groups of
strands, such that its class in $S_K$ is $\pi$. The example $\pi=(1,3)\in
S_3\,$ is depicted in the left of Figure~\ref{fig-perm-gg}. 
To see that the definition does not depend on the choice of the braid
element $b$ we observe that a standard generator of the pure braid group 
(see right of Figure~\ref{fig-perm-gg}) can be eliminated using
$\tau$-moves. 

\medskip

\begin{figure}[ht]
\begin{center}
\epsfxsize=5.5in
{ \ \epsfbox{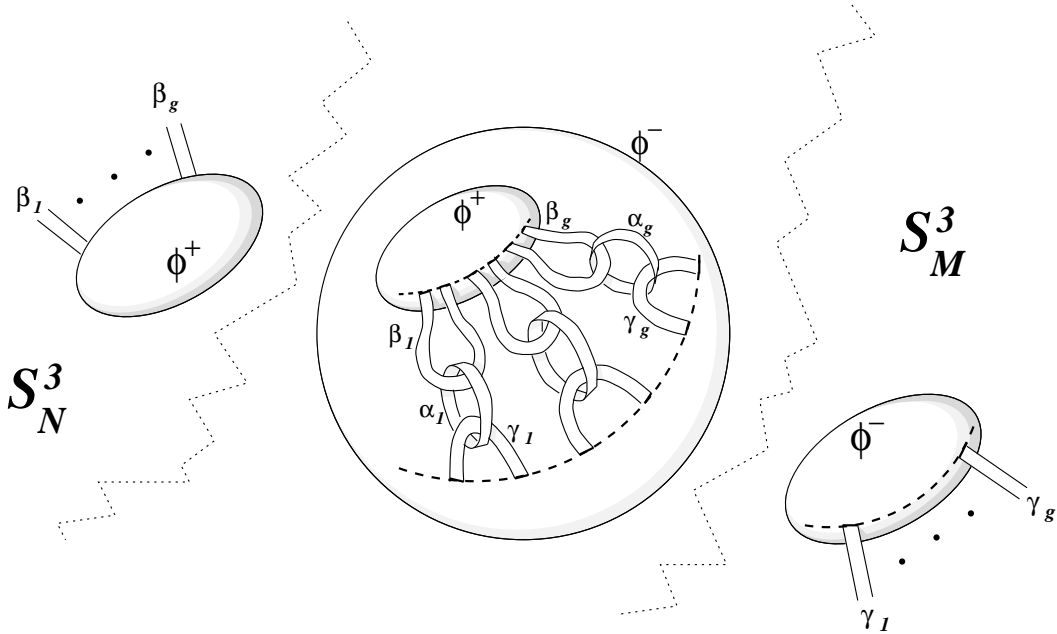}}
\end{center}
\caption{Composition of $H^{\pm}_g\,$'s}\label{fig-cobcomp-stan}
\end{figure}

The composition of two cobordisms is given by gluing corresponding
components of the standard surfaces together. The bridged link presentation of the
composite three fold close to a gluing surface
$\Sigma_g$  is obtained from the original presentations  by  taking the union  of
the bridged links from  $H^+_g$ and $H^-_g$ yielding a bridged link in an
$S^3$-component. 
In a standard presentation the ribbon diagram in this $\,S^3_j$-component
for the $j$-th surface is 
depicted in $b)$ of Figure~\ref{fig-cobcomp-stan}. It is connected to the
standard $S^3$ of either factor of  the cobordism product by the surgery spheres
$\phi^{\pm}\,$. Their partners in the standard $S^3_{N/M}\,$-components of
the two cobordisms are sketched in a) and c) of the same figure.
\medskip

The presentation is immediately simplified by carrying out the one-surgery
along $\,\phi^-\,$ explicitly,
which can also be thought of as a $0-1$-cancellation. In terms of the
elements of  Figure~\ref{fig-cobcomp-stan} this means that we insert the
content of the ball bounded by $\,\phi^-\subset\S^3_j\,$ into the region
bounded by $\,\phi^-\subset\S^3_M\,$. As a result the $\,S^3_j\,$-component
disappears and the vicinity of the $\,\phi^-\,$-sphere in the
$\,S^3_M\,$-components is relaces by the link digram around $\,\phi^+\subset
S^3_N\,$.  
\medskip

 We are left with a a bridged link presentation on the union  of the
standard $S^3$'s, that are part of the presentation of either cobordism,  
and the handle bodies $\,H^{\pm}_g\,$, which are not glued, and on which
the presentation at the boundaries is of standard  form. 
The presentation of the composite itself is however not standard since we have
two instead (of one)  $S^3$-components, and the pairs of surgery spheres
$\,\phi_j\,$, that survives the simplification described above, connect
them to each other.

To obtain a standard presentation we need to be more specific about the
ordering of the surfaces and the  connectivity of the cobordisms. We assume
that we are in the situation of equation~(\ref{elem-comp}).
\medskip

We denote the morphisms $M:\bar p\,\to\,\bar g^2\otimes\bar g^3\,$ and
$N:\bar g^1\otimes \bar g^2\,\to\,\bar q\,$, with $\bar g^2\neq
\emptyset\,$, so that $L:\bar g^1\otimes\bar p\to\,\bar q\otimes \bar
g^3\,$. On the standard sphere $S^2$ we mark the groups $\bar G^1_j\,$,$\bar G^2_j\,$,
and $\bar G^3_j\,$ in the given order. We move the surgery spheres in the
corresponding positions of the spheres $S_M:=S^2_{+,M}\,$ and $S_N:=S^2_{-,N}\,$ in
the copies of $S^3\,$, which define the tangle presentations of $M$ and
$N$. Thus, the surgery spheres in $S_N\,$ which belong to the first group
and
the spheres in $S_M$, which belong to the third group  connect to
corresponding spheres in the $H^{\pm}_g$'s. The spheres in the second
groups of $S_N$ and $S_M$ connect to each other as described above. 
\medskip

We modify the presentation by introducing a 1-2-birth-point in $M$ given by a
pair of surgery spheres $\phi$ and $\phi'$ connected by a ribbon $R$.
We choose one of the spheres $\phi_{j_o}\,$ in the  second group $\bar
G^2$ and push the sphere $\phi'$ through it into the $S^3$ component of the
other cobordism. We may arrange it that $\phi$ sits inside the $D^3_-$ of
$N$ and $\phi'$ inside $D^3_+$ of $M$. We then expand the spheres $\phi$ and
$\phi'$ until they coincide with $S_N$ and $S_M$ respectively. Since
$\pi_0\bigl(Diff(S^2)^+\bigr)$ is trivial we may deform the isomorphisms
$\phi\to S_N$ and $\phi'\to S_M\,$  such that
they are compatible with the respective identifications with 
standard spheres (and the arrangement of the other surgery spheres).
\medskip

\begin{figure}[ht]
\begin{center}
\epsfxsize=5.5in
{ \ \epsfbox{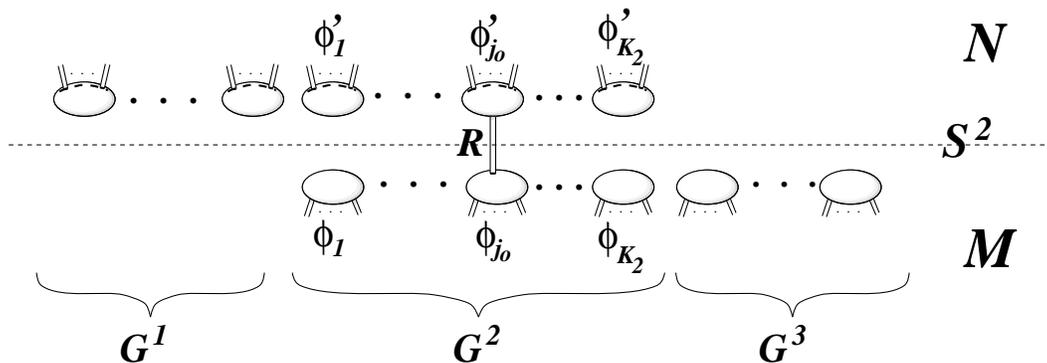}}
\end{center}
\caption{Composition of $S^3\,$'s}\label{fig-mid-pres}
\end{figure}

We  now carry out the index one surgery explicitly by gluing
the $D^3_+\cup S^2\times I\,$-piece of the $N$-presentation to the 
$S^2\times I\cup D^3_-\,$ piece of the $M$ presentation along $S_N\,$ and
$S_M\,$. The result is a bridged link diagram in a single $S^3\,$, with a
natural decomposition $D^3_+\cup S^2\times I\cup S^2\times I \cup D^3_-\,$.
Along the sphere in the middle the presentation has the form as in Figure
\ref{fig-mid-pres} .
\medskip

The respective surgery spheres of the second group are in the correct
position for the modification move, or, in the case of the sphere
$\phi_{j_o}\,$ through which we chose to move $\phi'\,$, a cancellation along
the ribbon $R$. The spheres of the first and third group are pushed into
corresponding positions on $S^2_{-,M}$ and  $S^2_{+,N}$ respectively.
\medskip

\begin{figure}[ht]
\begin{center}
{ \ \epsfbox{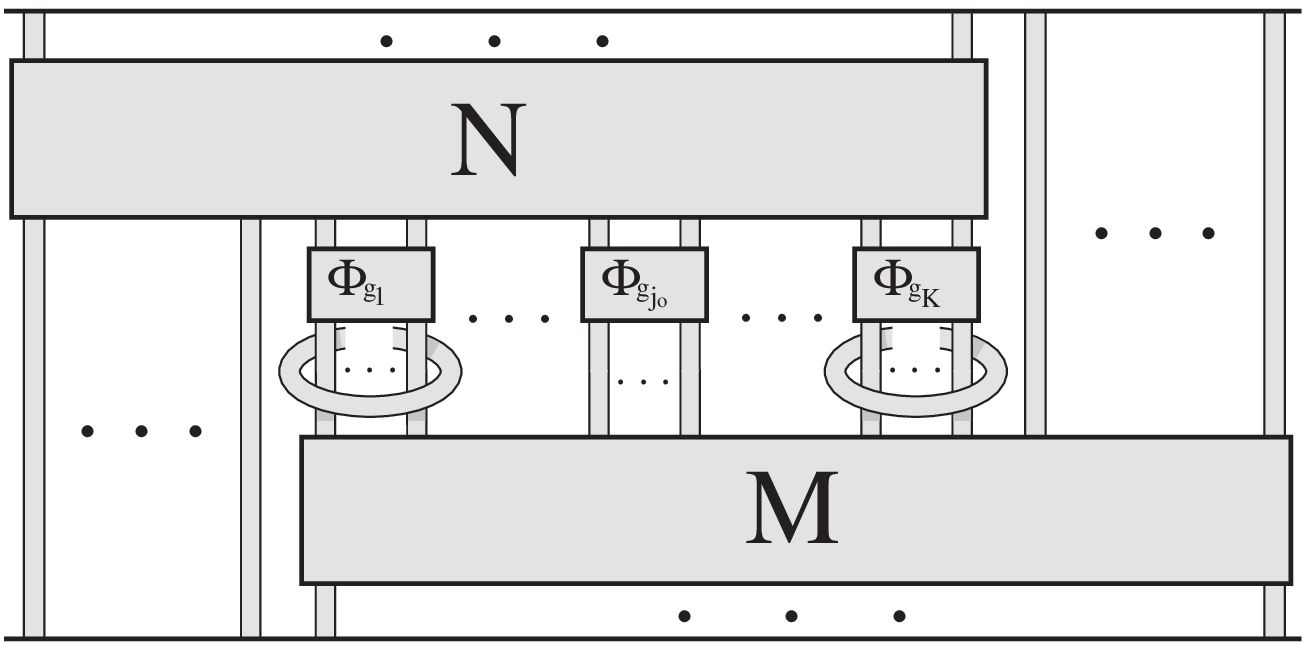}}
\end{center}
\caption{Composition of Tangles}
\label{fig-tang-comp}
\end{figure}

Resizing the $S^2\times I\amalg_{S^2}S^2\times I$ part  we obtain a
standard presentation, and for a generic projection a presentation of the
composite in terms of admissible tangles. It is now clear how the
composition rule of tangles should look like in order to give a functorial
presentation 
of the cobordism category $\wct$. Its definition is summarized in
Figure~\ref{fig-tang-comp}. 

\begin{figure}[ht]
\begin{center}
\epsfxsize=3.6in
{ \ \epsfbox{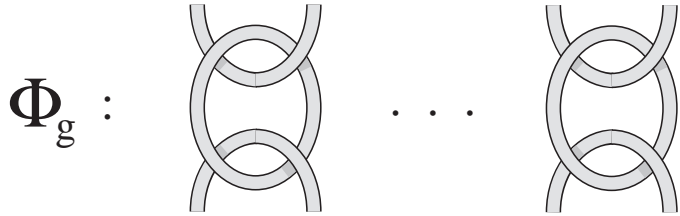}}
\end{center}
\caption{$\Phi_g\,$-Elements}
\label{fig-tang-coco}
\end{figure}

The boxes $M$ and $N$ stand for the tangle
presentations of the respective cobordisms. The content of the boxes
$\Phi_{g_j}\,$ is given in Figure~\ref{fig-tang-coco}. Note that 
$\Phi_{g_{j_o}}\,$ is the only group that is not decorated with an annulus.
\medskip

\s


{\em 3.2.4) Na\"\i ve Compositions, Connecting Annuli, and Closed Tangles}\label{pg-324}
\medskip

In situations, where the ribbons of a presentation that start from an
interval $I$ at the
lower end of the presentation do not return there in an interval
different from the partner interval of $I$ we may apply the 
Lemma~\ref{lem-con} of Connecting Annuli and replace the boxes $\Phi_g\,$
(or the respective pieces 
in Figure~\ref{fig-tang-coco})  by vertical strands.
\medskip

\begin{figure}[ht]
\begin{center}
{ \ \epsfbox{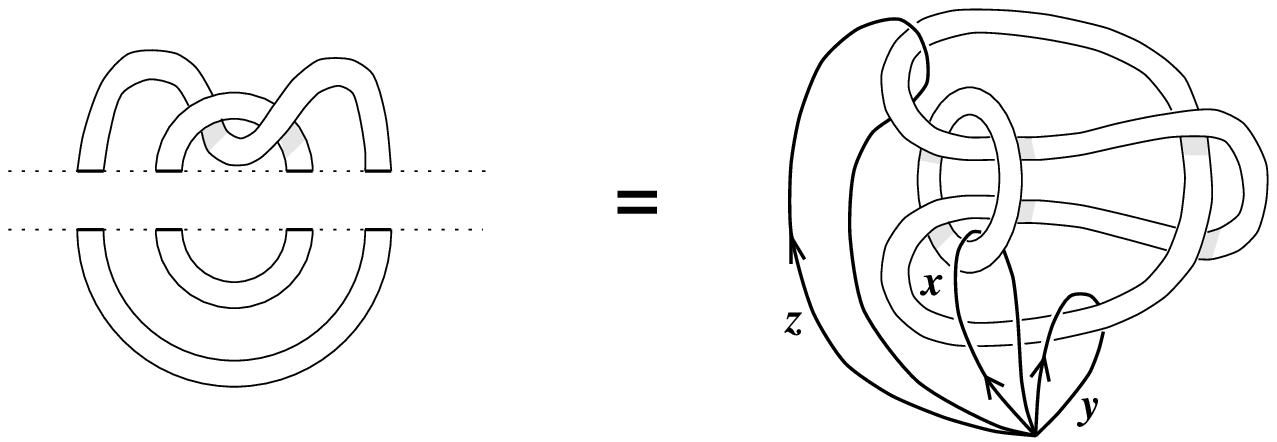}}
\end{center}
\caption{Example of $\emptyset\to\Sigma_2\to\emptyset$ Composition: $\,M_{\cal W}\,$}
\label{fig-link-ex}
\end{figure}

Let us give a simple  example with one boundary component, 
where this cancellation is not possible:
\medskip

We choose cobordisms $\Sigma_2\to\emptyset$ and $\emptyset\to\Sigma_2\,$
given by the tangle presentations on the left side of Figure~\ref{fig-link-ex}.

The na\"\i ve composition of these tangles  yields the
$\,^0\bigcirc\mkern-15mu\bigcirc^0\,$ link that can be
removed all together, i.e., we obtain $S^3$. The correct 
composition leaves us with the Whitehead link $\cal W$ on the right side
of Figure~\ref{fig-link-ex}. 
The three manifold $M_{\cal W}$ 
presented by it is however nontrivial. In particular, we have
\footnote{The link $\cal W$ can be presented as the closure of the braid
$\sigma_1^{-2}\sigma_2^{-2}\sigma_1^2\sigma_2\,$. With this we may compute,
as in [Bi] or [Rf], $\pi_1\bigl(S^3-{\cal L}\bigr)\,$. It is the free group in $x$,
$y$, and $z$ with relations $x^{-1}zx=zyz^{-1}\,$ and
$xyx^{-1}=yzy^{-1}\,$. The loops along the ribbons are $zy^{-1}$ and
$xyx^{-1}y^{-1}x^{-1}yxy\,$. Dividing by these relations gives the asserted
group.}) 
$\pi_1\bigl(M_{\cal W}\bigr)={\bf Z}(x)\,\oplus\,{\bf Z}/2(y)\,$.
For a description of
$M_{\cal W}\,$, when the components of the Whitehead link have framings
$\pm 1$ consult [Rf]. 

\s

The composition rules can be simplified at  the expense of starting with a
more specialized class of tangles. Let us define the space of {\em closed
tangles}, $c\!\widetilde{{\cal T}\!\!g\,}$, by
the  property  that the ribbon $R_j^{\pm}$ starting at an interval
$I_j^{i}\subset X_{\pm}\,$ ends in the partner interval $I_j^{o}\subset
X_{\pm}\,$. An even smaller space is given by the {\em special, closed
tangles}, $s\!c\!\widetilde{{\cal T}\!\!g\,}\subset c\!\widetilde{{\cal
T}\!\!g\,}$, for which each ribbon  $R_j^{\pm}\,$  with the additional segment
$b_j^{\pm}\,$ inserted
bounds a disks $\,D_j\,$. We require that $\,D_j\,$ is penetrated by only
one  ribbon, which has to be different from the 
$R_j^{\pm}\,$'s.

It is not hard to show that every closed tangle is
equivalent  to a special closed tangle if we  admit the moves of Section
3.2.2). The category of (specialized) closed tangles is endowed with
the {\em na\"\i ve composition rule}. This means we omit the insertions
$\Phi_g\,$, but keep the decoration rule.
\medskip

There is a {\em full } functor $\cc\,:\,\widetilde{{\cal T}\!\!g\,}\to s\!c\!\widetilde{{\cal T}\!\!g\,}\,:\,t\,\mapsto \,\cc(t)\,$ from the total space of
admissible tangles
to the space of special, closed tangles. In order to define it we
observe that the tangle in Figure~\ref{fig-tang-coco} can be written as a
(na\"\i ve) composition of its upper and its lower half
$$
\Phi_g\,=\,(S^+)^g\circ (S^-)^g\;.
$$
For $\tilde g^{\pm}\,=\sum_j g_j^{\pm}\,$ we then set
\beq\label{eq-cl}
\cc(t)\,:=\,(S^-)^{\tilde g^+}\circ t \circ  (S^+)^{\tilde g^-}\;.
\eeq
The composition $\circ$ used here is just placing the tangles on top of
each other. 
\medskip

Mainly by using $\sigma$-moves it is easy to see that $\cc$ and
the inclusion factor into isomorphisms, once we impose the relations on the
tangle categories. Summarily, we have
\beq\label{eq-tang-iso}
\wct\;\isto\;\overline{\widetilde{{\cal T}\!\!g\,}}\;\stackrel{\bar\cc}
{\isto}\;\overline{\widetilde{s\!c\!{\cal T}\!\!g\,}}\;\stackrel{\bar\i}
{\isto}\;\overline{c\!\widetilde{{\cal T}\!\!g\,}}
\eeq

%% file: BL4.tex
\section*{4) Applications and Implications}\label{pg-4}
\
\s

The applications we shall be concerned in the first part of this chapter
are of 
purely topological nature. They regard special classes of manifolds for which we 
derive presentations using a tangle presentations of the mapping class
group, derived from the results of the previous chapter.

In the second part  we extend the  Reshetikhin-Turaev invariant to 
presentations with bridged links and give a simple proof of invariance.
We shall discuss  the algebraic implications of the Lemma~\ref{lem-con}
for Connecting Annuli. From this we will find a canonical natural
transformation for a BTC and a canonical, central element of a
quasitriangular Hopf algebra, which projects onto selfconjugate objects.

We conclude with some remarks on how the presentations of cobordisms of two
folds with boundary are obtained. 
In particular we explain how the operation of a ``glue tensor
product'' acts on the presentations.
\s


\subparagraph{4.1) Invertible Cobordisms and Presentations of Mapping
Tori}\label{pg-41} 
\
\s

An interesting family of morphisms in $\wct$ are the invertible cobordisms
$Aut(\Sigma)\,$ for a
connected surface $\Sigma\,$. They are given by $M=\Sigma\times I\,$,
equipped with possibly non canonical charts, as for  $\psi=\psi'\amalg id\,$ with
$\psi'\in\diff(\Sigma)^+\,$ .
For example it is not hard to find the standard presentations for 
different Dehn twists, and construct a 
presentation of the mapping class group of $\Sigma$ in terms of tangles.
This presentation is identical to a the one constructed by [MP].

However, in [MP] the use of
Cerf theory in 3+1 dimensions was avoided by referring to the explicit
presentation of the mapping 
class group  of [Wj] in terms of generators and relations.

From the tangle form of $\pi_0\bigl(\diff(\Sigma)^+\bigr)\,$  we derive
surgery presentations of Heegaard splittings and 
general mapping tori over $S^1$, with fiber $\Sigma$. In a few examples 
we compare those to known presentation. Namely, lens spaces and
``planar presentations'' of trivial bundles of the form 
$\partial\bigl(N^{(4-j)}\times D^j\bigr)\,$.  

\s


{\em 4.1.1) Tangle Presentation of the Mapping Class Groups}\label{pg-411}
\s

We start this section with a  derivation of the presentation of the mapping
class group in the category of admissible tangles. To begin with, we
remark that the trivial tangle, $I\,:\,(g)\to(g)\,$, given by  $2g$
vertical ribbons, 
represents the identity cobordism in $\wct\,$. This follows directly from
Theorem~\ref{THREE} and identity the $I\cdot t \,=\,t\cdot I\,=\,t\,$, which is
by the composition rules holds obviously for any tangle $t\,$. 
As an instructive exercise let us derive this explicitly from the
topological situation: 

In the 
corresponding standard presentation we obtain two concentric surgery spheres $\phi^+$
and $\phi^-$ in $S^3$ with $2g$ straight, radial lines joining them. As in 
the derivation of the composition rule in Section 3.2.3) we carry out the
one-surgery for the 
$\,\phi^-\,$-sphere explicitely, i.e., we do a  0-1-cancellation.

The resulting  link-presentation of the cobordisms consists of the left 
side of Figure~\ref{fig-stan--} and the left side of
Figure~\ref{fig-stan-+}, where we consider the two spheres as partners.

\begin{figure}[ht]
\begin{center}
\epsfxsize=4in
{ \ \epsfbox{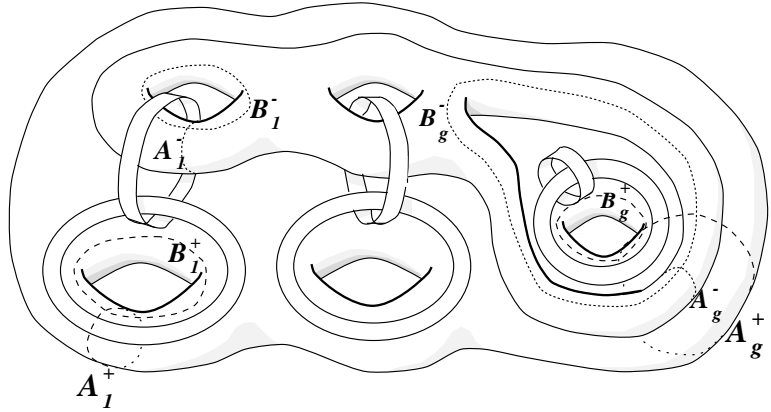}}
\end{center}

\caption{Identity on $\Sigma_g\,$}\label{fig-pres-id}
\end{figure}

We perform  the one-surgery along these spheres, too, so that we have a
presentation on the one-componest manifold $H_g^+\#H_g^-\,$, without
one-surgery data. 
The detailed result is given in Figure~\ref{fig-pres-id}.

Now, we
can use an isotopy which slides the one handles with meridians $A_j$ of the
inner handlebody 
over the $\beta_j$-ribbons  so that the meridians of the inner and the outer
handlebody are aligned. Using $\eta$-moves we may then remove all surgery 
ribbons. The resulting cobordism are just two nested copies if the same handle
body, i.e., $\Sigma\times I$ with canonical charts at the
boundaries. This is the identity cobordism.

\medskip

\begin{figure}[ht]
\begin{center}
\epsfxsize=4.8in
{ \ \epsfbox{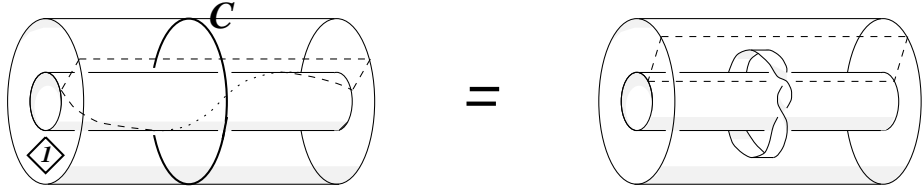}}
\end{center}

\caption{Dehn Twist as Surgery}
\label{fig-dehn-surg}
\end{figure}

The mapping class group of a surface $\Sigma$ is generated by
Dehn twists $\delta_C$ along a sufficient number of Jordan curves $C$
on $\Sigma$. The corresponding cobordism $<C>\,:=\,\Bigl(\Sigma\times I, id,
\delta_C\Bigr)\,$ is equivalent to one of the form
$<C>\,=\,\Bigl(\widetilde{\Sigma\times I}, id, \id\Bigr)\,$. Here 
$\widetilde{\Sigma\times I}$ is the same manifold with canonical boundary
maps, but a surgery done inside. The surgery presentation is given by
pushing the curve  $C$ inside 
the thickened surface $\Sigma\times I\,$ and inserting a ribbon with
framing number 1 with respect  the canonical framing of $T(S^3)_1\,$. In
Figure~\ref{fig-dehn-surg} the cylindrical neighborhood of $C$ is shown;
the equality follows from an ordinary isotopy.
\medskip

The symbol $\FLN$ is short hand for $|n|$ isolated 
ribbons with framing number $sgn(n)=\pm 1\,$. (It may be omitted if we are
only interested in the cobordism classes in $\ct$). The admissible tangles
giving the presentation of the Dehn twists can be obtained by composing the
$<C>$ with the standard presentation of the identity.
We obtain an additional link component in  $H^{\pm}_g$, which can be pushed into the
$S^3$ component with the map $\Psi\,$ (or only $\sigma$-moves for the
affected handles).
\medskip

\begin{figure}[ht]
\begin{center}
\epsfxsize=5in
{\ \epsfbox{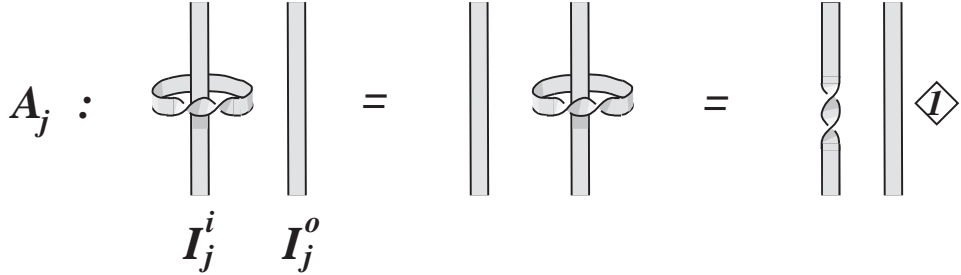}}
\end{center}

\caption{$A_j\,$ Dehn Twist}
\label{fig-Dehn-A}
\end{figure}

\begin{figure}[ht]

\begin{center}
\epsfxsize=5in
{\ \epsfbox{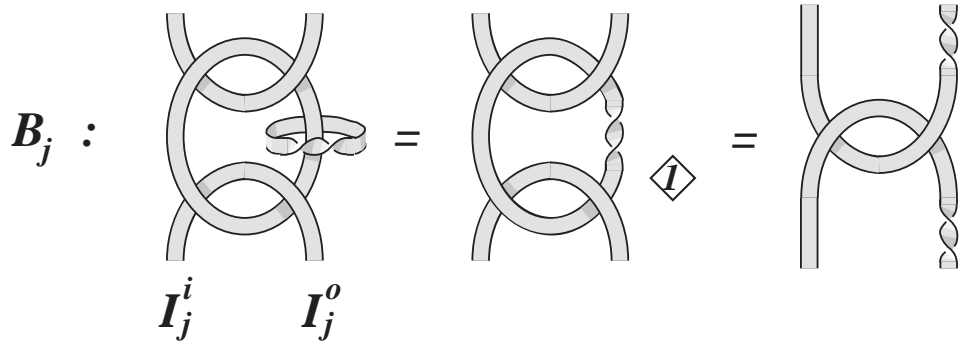}}
\end{center}

\caption{$B_j\,$ Dehn Twist}\label{fig-Dehn-B}
\end{figure}

Specifically, we obtain Figure~\ref{fig-Dehn-A} by
inserting a 1-framed ribbon along $A_j$ in $H^+_g\,$, 
and moving it along $\beta_j\,$ and through $\,\phi^+\,$ into
the $\,S^3\,$ component of the presentation. A Dehn twist along the curve
$\,B_j\,$ is given by placing a small one-framed ribbon around
$\,\alpha_j\,$ in $\,H^-_g\,$. A $\sigma$-move at the $j$-th handle
leaves us with Figure~\ref{fig-Dehn-B}.

\begin{figure}[ht]
\begin{center}
\epsfxsize=5in
{\ \epsfbox{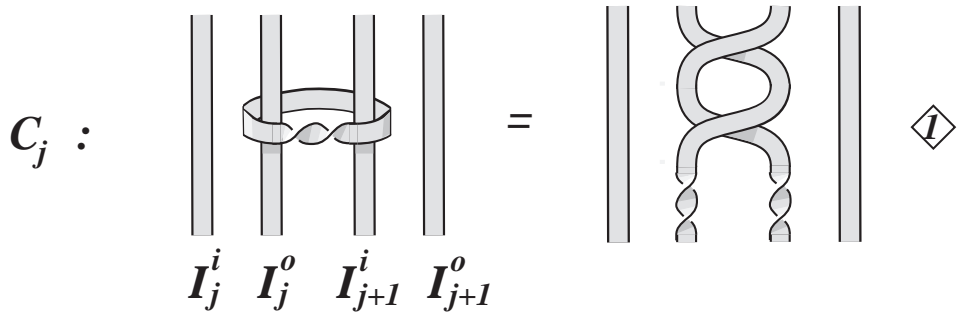}}
\end{center}

\caption{$C_j\,$ Dehn Twist}
\label{fig-Dehn-C}
\end{figure}

The curve $C_j$
is the one intersecting $B_j$ and $B_{j+1}$ exactly once and no other $A$
or $B$ curve. For the tangle presentation in
Figure~\ref{fig-Dehn-C} we insert the respective ribbon in $H^+_g\,$ 
move it towards $\,\phi^+\,$ along the pieces of $\,\beta_j\,$ and
$\,\beta_{j+1}\,$ that emerge `radially' from the intervals $\,I^o_j\,$ and
 $\,I^i_{j+1}\,$. 

Finally, the
curve $D_j$ is the one opposite to $A_j$, i.e., it intersects only $B_j$.
If the respective ribbon is inserted in $\,H^+_g\,$ ot can be pushed onto
a curve $\cal E\,$ on the sphere $\,\phi^+\,$. Now, $\cal E$ separates
$\,\phi^+\,$ into two hemispheres one containing the intervals
$\,I^i_1,I^o_1,...,I^i_j\,$ and the other containing the intervals
$\,I^o_j,I^i_{j+1},...,I^o_g\,$. Thus in the $\,S^3\,$ the respective
ribbon can be moved around the surgery sphere into the $S^2\times I\,$-piece
of the presentation in two ways. The results are depicted in
Figure~\ref{fig-Dehn-C}. Clearly the two possibility differ by exactly one
$\tau$-move since in one instance $\cal E$ was moved through the special
line $L$ from Figure~\ref{fig-isto-II's} in another it was not.

\begin{figure}[ht]
\begin{center}
\epsfxsize=5in
{\ \epsfbox{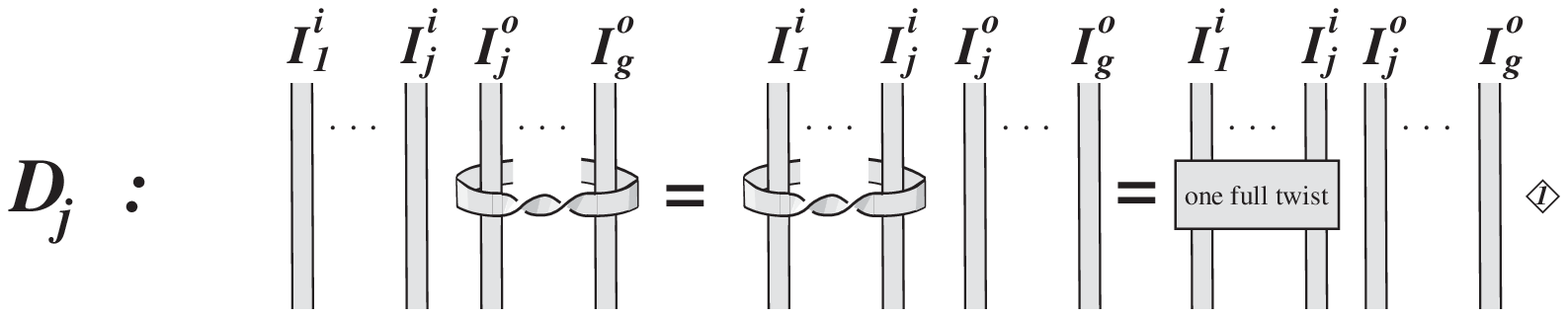}}
\end{center}

\caption{$D_j\,$ Dehn Twist}\label{fig-Dehn-D}
\end{figure}

The equalities in all of the pictures follow from $\kappa$-moves.

\s

Let us consider the subring, $\cal T\,$, of  cobordisms that are generated by 
the $A_j$'s, the $B_j$'s, and tangles presenting  generators of the pure braid 
group $P_{2g}\,$, as depicted in Figure~\ref{fig-perm-gg}.
For compositions it follows by  induction that a strand of a tangle in
$\cal T$ starting at an interval in $X_-$
will end either at its partner interval on $X_-$ or at its own copy in
$X_+$.
Thus, by Lemma~\ref{lem-con}  we may use the na\"\i ve composition rule 
for all elements in $\cal T$, so that $\cal T$ is in fact a group.  
Clearly, all
tangles of the mapping class group presentation lie in $\cal T$.
In fact the converse is also true.

\bpp \label{lem-T=mcg}The group $\cal T$ is isomorphic to $\pi_0\bigl(
\diff(\Sigma)^+\bigr)\,$ via the presentation in terms of  admissible tangles.
\epp

This remark follows immediately if we use that all invertible cobordisms
are of the form $\bigl(\Sigma\times  I, id, \psi\bigr)\,$. However, the generators
of $P_{2g}\,$ can also be produced directly using special Dehn twists. For
example, a Dehn twist along the curve $A_j*B_j*A_j^{-1}*B_j^{-1}\,$ yields a
full twist of the strands at $I_j^i$ and $I_j^o$. (Here * is the
composition of paths as in $\pi_1\,$.) Similarly, we obtain from a twist
along a curve along $D_i*D_j^{-1}\,$ (not intersecting the $A$'s and $B$'s
anywhere else) the full twist of the strands $I^o_i,\ldots, I_j^i\,$.

It is often more useful to replace the generators $A_j\,$ by the generators
\beq\label{eq-new-gen}
S^+_j\,=\,B_jA_jB_j\;,\qquad\quad (\,S_j^-:=(S_j^+)^{-1}\,)
\eeq
which is identical to the tangle $S^+\,$ introduced in~(\ref{eq-cl}).

Let us remark as a word of warning that, e.g., for $g=2\,$ the simple braid
of the strands $I^o_1$ and $I^i_2$ presents a non-invertible cobordism
$F\,$. This follows immediately from $F\circ A_1=F\circ A_2\,$. In
particular the na\"\i ve composition  is not applicable. Note also, 
that the simple braid on a pair, $I_j^i\,,I_j^o\,$, is equivalent to $S_j^+\circ S_j^+\,$.
\medskip

The results in [MP], specifically  Propositions 5.2 and 6.1, follow directly from
Proposition~\ref{lem-T=mcg}. When comparing our presentation to the  one in
[MP] we have to keep in mind that 
we have to include the
$\tau-move\,$, since we consider surfaces without  punctures.$\smile$
\medskip

A more detailed analysis of the presentation of
$\diff\bigl(\Sigma\bigr)\,$ as a product of the pure braid group $P_{2g}$, 
the group ${\bf Z}^g\,$, generated by $S_j^{\pm}\,$, and the group ${\bf
Z}^g\,$ generated by $B_j\,$, should also give an alternative proof of
the results in [Wj].   
\s


{\em 4.1.2) Presentations of Manifolds from ${\cal
T}=\pi_0\bigl(\diff(\Sigma)\bigr)\,$}\label{pg-412}  
\
\s

From the presentation of cobordisms we may derive link presentations in
$S^3$  of closed manifolds. We start with the easier example.
\s

{\em Heegaard-Splittings:}
\s

\begin{figure}[ht]
\begin{center}
\epsfxsize=5in
{\ \epsfbox{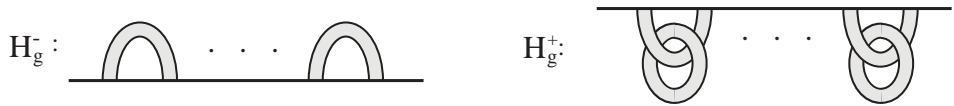}}
\end{center}

\caption{Presentations of $H^{\pm}_g\,:\,\emptyset \leftrightarrow \Sigma\,$}
\label{fig-H-pres} 
\end{figure}

Any  three fold $M$ can be presented  by a Heegaard-splitting, i.e.,
we glue $H^+_g$ to $H^-_g\,$, where the boundary identification is given by
 an element in $\psi\in \diff(\Sigma)^+\,$,
with associated cobordism $<\psi>\,$.
Now, the manifolds $H^{\pm}_g$ may also be considered as cobordisms with 
tangle presentation as in Figure~\ref{fig-H-pres}. Hence  $M$ may be written
as the composite of the three cobordisms $H^-\circ <\psi>\circ H^+\,$. 

For the
composite $<\psi>\circ H^+_g\,$ we simply close the tangle presenting
$<\psi>\,$ at the bottom with the ribbons $b_j^-\,$. It is clear from the
description of $\cal T$ that the strand emerging from $I^i_j$ at $X_+\,$
ends at $I^o_j\,$. We may  apply $\eta$-moves to these $g$ different strands
after composing with $H^-_g$. Thus the link presentation is obtained from
the presentation of $<\psi>\circ H^+_g$ by
omitting the strands ending in $X_+\,$, ( and adding annuli if the closure
of a strand goes through several intevals). 
\medskip

The corresponding  presentation of a lens space $L\bigl(\psi \bigr)\,$ with
$\psi=\prod_jS_1^+T_1^{n_j}\,$ is, for example, easily identified with the
familiar chain of unknots with framing numbers $\{ n_j\}\,$.
\s

{\em Mapping Tori:}
\medskip

\begin{figure}[ht]

\begin{center}
\epsfxsize=5in
{\ \epsfbox{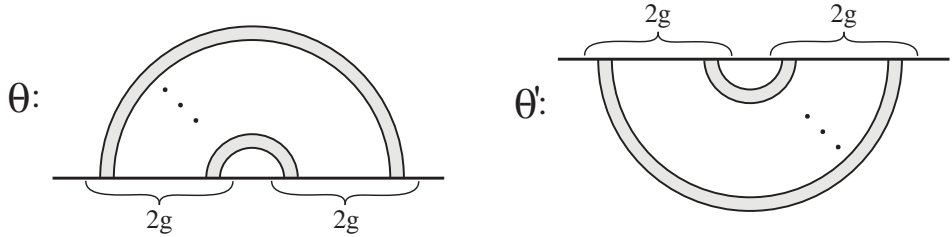}}
\end{center}

\caption{Pairings $\Sigma\times I\,: \,\Sigma\amalg-\Sigma\leftrightarrow
\emptyset\,$}
\label{fig-cob-rig} 
\end{figure}

A more interesting case is provided by bundles  over $S^1\,$ with a connected
surface $\Sigma$ as fiber. They are classified by conjugacy classes in
$\pi_0\bigl(\diff(\Sigma)\bigr)\,$ and can be given as the mapping torus of
a representative.  The main ingredient for the surgery description is to fix
a pair of ``rigidity morphisms'' of $\wct$:
\medskip

For a connected surface of genus $g$ let us introduce cobordisms
$\theta\,:\,\Sigma\amalg\Sigma\to\emptyset\,$ and
$\theta'\,:\,\emptyset\to\Sigma\amalg\Sigma$ 
as indicated in Figure~\ref{fig-cob-rig}.

It is easy to see that \beq\label{rigid}
(\theta\otimes id)\circ (id\otimes \theta')\,= \,(id\otimes\theta)\circ (\theta'\otimes
id)\,=\,id\;,
\eeq
using the composition rules for admissible tangles. In general, a cobordism
$\Sigma\amalg\Sigma\to\emptyset\,$ may be constructed from $\Sigma\times I$
with canonical boundary charts by composing one of the components of the
chart  with an orientation reversing map $\psi\in\diff(\Sigma)^-\,$. 
Using the relations (\ref{rigid})
and the general form of an invertible cobordism, it is easy to show that
$\theta$ is also of this form. 
For the closed composition we have the following identity:

\blm  The composition $\theta\circ\theta'$ is homeomorphic to
$\Sigma\times S^1\,$. 
\elm
\medskip

{\em Proof:} The direct proof is easy using the previous remarks on the general
structure of cobordisms $\emptyset\to\Sigma\amalg\Sigma$ . In fact, it
follows that $\theta$ and $\theta'$ are of the form 
$\bigl(\pm\Sigma\times I, id, \psi\big)\,$ from which the assertion follows
immediately. Nevertheless, we wish to give a more complicated proof which 
reveals another surgery presentation of $\Sigma\times S^1\,$ starting
directly from the link diagram of $\theta\circ\theta'\,$. 

\begin{figure}[ht]
\hspace*{.5in}\epsfbox{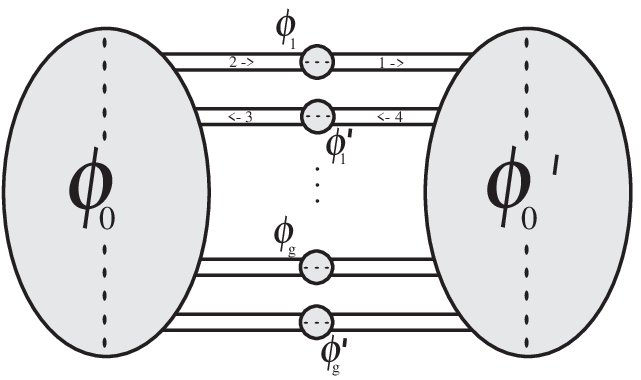}\qquad\qquad\epsfbox{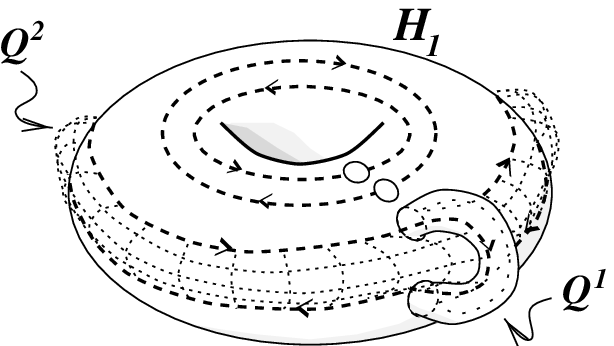}
\caption{Presentation of $\Sigma\times S^1\,$}
\label{fig-ppp}
\end{figure}

For the composite we
may use the na\"\i ve composition rule  for one component and insert $\Phi_g$ and
the extra ribbon $R$ for the other component. The elements in $\Phi_g$ may not
be replaced by the identity (this would yield a presentation of the
connected sum of $g-1$ copies of $S^1\times S^2$). However, we may apply 
an (un-) modification, which introduces $g$ pairs of surgery spheres.
We do the same with $R$, so that we end up with $g+1$ index one surgeries
and $g$ index two surgeries. The resulting surgery presentation is now planar
and is shown on the left of Figure~\ref{fig-ppp}. Following the first
ribbon we pass through the pieces $1, \phi_0', \phi_0, 2, \phi_1 , \phi_1' , 3,
\phi_0 , \phi_0' , 4 , \phi_1' , \phi_1 , 1\,$.
\medskip

It is a general principle that if a surgery graph is planar the surgered
four manifold $W$ can be given as a product $N\times I$ and the planar
graph can be used to give a surgery description of the three fold $N\,$. 
More specifically, write $D^4=D^3\times I\,$, so that we have a standard
piece $S^2\times I=\partial(D^3)\times I\,\subset S^3=\partial(D^4)\,$.
Now, for a planar diagram the attaching curves for $j$-handles
$S^j\times D^{3-j}\hookrightarrow S^3\,$ can be brought into the form
$i\times id_I\,:\,S^j\times D^{2-j}\times I\hookrightarrow S^2\times I\,$.
It is clear that instead of attaching a four dimensional $j$-handle
$\hh^j=D^{j+1}\times (D^{2-j}\times I)\,$ to $D^4=D^3\times I\,$ we may as well
attach a three dimensional $j$-handle $Q^j\,=\,D^{j+1}\times D^{2-j}$ to $D^3$ and
form the product with $I\,$.
\medskip

In our example, attaching a one handle to $D^3$ along the discs
$\phi_0,\phi_0'\subset S^2\,$ yields the solid torus $H_1$ on the right of 
Figure~\ref{fig-ppp}, where the ribbon pieces $1$ and $2$ (or $3$ and $4$)
are glued to loops in $\partial H_1\,$. If we attach another $Q^1$ along
discs $\phi_j,\phi_j'$ we obtain the depicted handlebody $H_2$ of genus 2, and the
corresponding attaching (Jordan) curve for $Q^2\,$. It is now clear from the
picture that the manifold   $H_2\cup Q^2$ is homeomorphic to the solid
torus $H_1$ with another torus removed from the inside. It can also be
described as $S^1\times A_1$, where $A_1$ is an annulus. The total surgery
will result in the manifold $ S^1\times A_g\,$, where $A_g$  is the disc
$D^2$ with $g$ small discs removed from the inside. 
\medskip

The surgered four-manifold is therefore $W\,=\,S^1\times A_g \times I\,$. 
But with $A_g\times I=H_g\,$ we find $M\,=\,\partial W\,=\,S^1\times \Sigma\,$.  
$\smile$ \hfill$\Box$
\s

{\em Remark:} Clearly, the class of manifolds that have a planar bridged
link presentation and as above allows a reduction of dimensions is much larger
then the class of manifolds with planar link diagrams (which are just
connected sums of $S^1\times S^2\,$'s).
\s

\begin{figure}[ht]

\begin{center}
\epsfxsize=5.5in
{\ \epsfbox{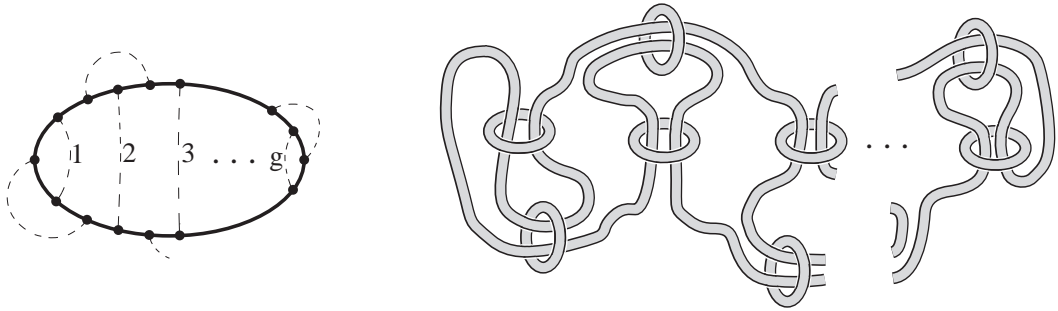}}
\end{center}

\caption{Other Presentation of $\Sigma\times S^1\,$}\label{fig-2low-surg}
\end{figure}

{\em Remark:}\footnote{This presentation has been communicated to me by
Robion Kirby.} It is in fact possible to reduce the presentation of
$\Sigma\times S^1$ by two dimensions using $W=\Sigma\times D^2$.
A surgery presentation of $\Sigma$ is given by attaching to $D^2\,$
$2g$ one handles in the way indicated in the left of
Figure~\ref{fig-2low-surg} and closing the manifold with a two-handle. Here
the  attaching curves are pairs of points $\bullet$ in the boundary $S^1$
joined by a dashed line. The attaching curve for the two handle passes
through every piece of $S^1$ exactly once. The corresponding link diagram
for $\Sigma\times 
S^1$ is indicated on the right of Figure~\ref{fig-2low-surg}. For $g=1$ this
presentation is identical to the one above; for $g>1$ the equivalence of the
link presentations is left as an exercise to the reader.$\smile$   
\s
\s 
 
\begin{figure}[ht]

\begin{center}
\epsfxsize=3.5in
{\ \epsfbox{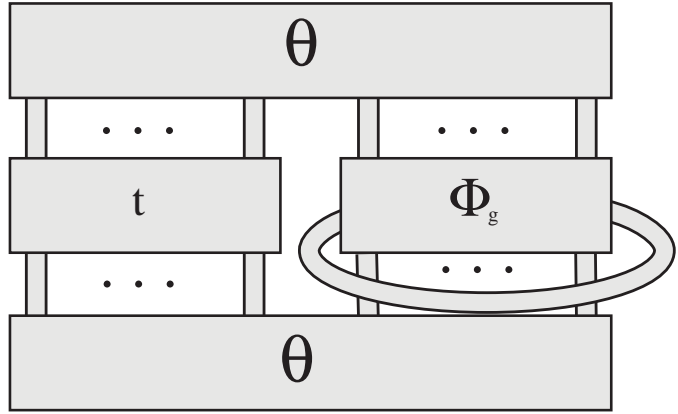}}
\end{center}

\caption{Tangle Presentation of Mapping Torus}
\label{fig-pres-MT}
\end{figure}

It is clear now that the mapping torus of $\psi\in \diff(\Sigma)^+\,$  is
given by the composition, $\theta'\circ \bigl(<\psi>\otimes
id\bigr)\circ\theta\,$.
We give a bloc diagram of the  corresponding link presentation in 
Figure~\ref{fig-pres-MT}. In the box $t$ the tangle corresponding to $\psi$
is inserted. This concludes our discussion of surgery presentations of
mapping tori. 
 
\subparagraph{4.2) On the Reshetikhin Turaev Invariant}\label{pg-42}
\
\s

The construction of the invariant of closed three dimensional manifolds 
as in [RT]
is based on close relations between tangle categories and general, abelian  BTC's.
The new elements of bridged links and admissible tangles we have
encountered so far also have natural counterparts in semisimple BTC's and
finite-dimensional nicely-quasitrangular Hopf algebras.  

The two aspects we wish to address here are an extension of the [RT] -
invariant for bridged ribbon presentations as described in Section 2.4)
and a discussion of the Lemma~\ref{lem-con} for Connecting Annuli; in
particular the 
structure of the transformation associated to the tangle in
Figure~\ref{fig-conn-lm}. 
\s


{\em 4.2.1) Invariants  of Three Manifolds from Bridged
Ribbons}\nopagebreak\label{pg-421}
\medskip

Throughout this section we use the same notation as in [Tu]. Starting point
of the construction is an abelian, strict, semisimple, balanced BTC, with only a
finite set $\cal I$
of isomorphism classes of irreducible objects, each of which contains only
one element. We denote the braid element by $c_{i,j}:i\otimes j\to j\otimes
i\,$ and the balancing $\theta_j\in End(j)={\bf C}\,$, so that
$\theta_{i\otimes j}\,=\,c_{j, i}c_{i, j}\theta_i\otimes
\theta_j\,$.
Rigidity provides us with a pair of morphisms $1\to X\otimes X^{\vee}\,$
and $X^{\vee}\otimes X\to 1\,$. We define the maps 
\beq\label{tep}
\rho_X\;:\;\bigl(\theta_X\otimes
1\bigr)\,c_{X^{\vee},X}\;:\;Hom(1,X^{\vee}\otimes X)\,\isto\,Hom(1,X\otimes
X^{\vee})\,,
\eeq
with $\rho_{X^{\vee}}\rho_X=id\,$ and the corresponding ones on
$Hom(X^{\vee}\otimes X,1)\,$. Applying these to the rigidity morphisms we
produce corresponding morphisms for the opposite product,
$\otimes'$. These morphisms are associated to maxima and minima  in a directed,   colored ribbon graph. The morphisms $1\to 1$ associated to an annulus, with a
morphism $f:X\to X$ inserted, defines a canonical, generally cyclic, and
$\otimes$-factorizable trace
$$
tr_X\,:\,End(X)\,\to\,{\bf C}\;.
$$
As usual we define the $S$-matrix and the $q$-dimensions:
\beq\label{def-S}
S_{i,j}\,=\,tr_{i\otimes j}(c_{j, i}c_{i, j})\qquad\qquad
dim(j)=S_{1,j}=\,tr_j(1)
\eeq
\medskip

In the construction of a three manifold invariant it is usually required 
that the $S$-matrix is invertible. (This is sometimes called the
``modularity'' axiom.) In the bridged link formalism it suffices to start
for a seemingly weaker condition. All we require is that there is a vector
$\hat d$ with $S\hat d\,=\,1\,$, i.e.,
\beq\label{eq-d-hat}
\sum_j \hat d(j)S_{j,i}\,=\,\delta_{i,1}\;.
\eeq
It will turn out that~(\ref{eq-d-hat}) also implies $\hat d(1)\,\neq
\,0\,$, but we shall add this property  here to the list of assumptions.
\medskip

Let us now give a definition of a functional $BR\to\{ BR\}\,$ on the bridged links in unions
$U=\amalg_{\omega =1}^rS^3\,$, from which we wish to define an invariant of
the three manifolds they are presenting. The prescription to compute
$\{\;\}\,$ is as follows:
\medskip

For each object $X$ of the semisimple category we introduce bases
$$
\{f^X_{\alpha}\}_{j\in\Lambda_X}\subset Hom\bigl(X,1\bigr)\qquad\; {\rm and}\qquad\qquad \{e^X_{\alpha}\}_{\alpha\in\Lambda_X}\subset
Hom\bigl(1,X\bigr)\,
$$ 
\qquad\qquad with $\hphantom{xxxxxxxxxxxxxxxxx}
f^X_{\alpha}\circ e^X_{\beta}\,=\,\delta_{\alpha\beta}1_1\,$.
\medskip
 
A coloration is now not only a labelling of the directed ribbons with
elements $l\in\cal I$ but in addition a labelling of a pair of coupons $(j,\pm)$
with an (the same) element of $\Lambda_X$. Here, $X=l_1\otimes\ldots
\otimes l_n\,$,
where $l_k$ are the labels of the ribbons entering the coupon $(j,+)$.
\medskip
 
In the next step we associate to a plat bridged ribbon graph $BR$ with
coloration $C$  a composition of elementary morphisms.
As in [RT0] we insert braid and rigidity morphisms in place of the tangle
elements depicted in Figure~\ref{fig-cdr-elem}. Furthermore, we assign
$f^X_{\alpha}$ to the 
coupon $(j,+)$ with color $\alpha$, and $e^X_{\alpha}$ to $(j,-)$
respectively. For a closed, bridged ribbon we obtain a morphism $1\to 1$ in every
component, $S^3_{\omega}\,$, and thus a number
$F(BR,C,\omega)$. From this we define the number
\beq\label{eq-def-[]}
\{ BR\}\,:=\,\sum_{C}\prod_{\omega=1}^rF(BR,C,\omega)\prod_Ld\bigl(C(L)\bigr)\;,
\eeq
where $\hphantom{xxxxxxxxxxxxxx}d(j)\;:=\;\frac {\hat d(j)} {\hat
d(1)}\;.$
\medskip

The product runs over all components $L$ of the link diagram in
$\breve M$ and $C(L)\in \cal I$ is the coloration associated to $L$ by $C$.
\s

Let us also define two classical  invariants of the bounding four fold
$W_{\cal L}\,$. The homology of $W_{\cal L}\,$ only depends on
the bridged link $\cal L$. It is given by the cellular complex
$(C,\partial)\,$, where $C_j$ has as a basis
the handles, so $dim(C_2)$ is the number of ribbons in the presentation,
$dim(C_1)$ is the number of pairs of coupons, and $dim(C_0)$ is the number $r$
of components of $U$.
The boundary operation $\partial:C_2\to C_1$ is given by counting
the number of times a ribbon passes through a surgery sphere including
signs for directions. $\partial:C_1\to C_0\,$ assigns to a basis element
the difference of the
components of $U$ in which the two surgery spheres lie. 
\medskip 
 
In particular we have for the Euler number of $W_{\cal L}\,$ the formula 
$$
\chi({\cal L})\,=\, dim(C_2)-dim(C_1)+dim(C_0)\;.
$$
The signature $\sigma({\cal L})$ of $W_{\cal L}\,$ is
 given by the signature of the linking matrix of $\cal
L\,$ restricted to the kernel of $\partial$. 

\s

It is clear that the functionals $\{.\}\,$, $\sigma$, and $\chi\,$ are
isotopy-invariants of the bridged link. The fact that $\{.\}\,$ is
independent of the choice of the projection follows as in [RT0] from
$Rel_1\,-\,Rel_{11}\,$. Naturality of the balancing and $\theta_1\,=\,1$
implies $Rel_{12}\,$. Since all functionals are obviously invariant undet
0-1-cancellations, we may confine ourselves to the situation where $r=1\,$,
i.e., the presentation is on only one $S^3\,$.
In order to determine how they change under the remaining moves
$2),\,3)\,$, and $4)\,$ of Theorem~\ref{TWO} let us explain some notations 
for bridged links: 
\medskip

For a given link $\li\,$ let $\li\cup\clr\,$ 
be the link with an additional, isolated cancellation diagram as in
$Rel_{18}$ of Figure~\ref{fig-rel-17}. Similarly we define
$\,\li\cup\frL\,$ as the diagram where an isolated, -1-framed unknot is
added. If a pair of coupons is in a position as on the right in
$Rel_{17}$ of Figure~\ref{fig-rel-17} we denote by $\li\cup{\eta}$ the 
link where the respective piece is replaced by the diagram on the left hand
side of $Rel_{17}$. We have the following ``transformation rules'':

\blm\label{lm-trsf} 
\
$$
\begin{array}{rclrclrcl}
\{\li\cup \clr\}\,&=&\,\{\li\}&\chi\bigl(\li\cup
\clr\bigr)\,&=&\,\chi\bigl(\li\bigr)&\sigma\bigl(\li\cup
\clr\bigr)\,&=&\,\sigma\bigl(\li\bigr)\nonumber\\
\{\li\cup{\eta}\}\,&=&\,\frac 1 {\hat
d(1)}\{\li\}&\chi\bigl(\li\cup{\eta}\bigr)\,&=&\,
\chi\bigl(\li\bigr)+2&\sigma\bigl(\li\cup{\eta}\bigr)\,&=&
\,\sigma\bigl(\li\bigr)\nonumber\\
\{\li\cup \frL\}\,&=&\,\{\li\}\{\frL\}&\chi\bigl(\li\cup
\frL\bigr)\,&=&\,\chi\bigl(\li\bigr)+1&\sigma\bigl(\li\cup
\frL\bigr)\,&=&\,\sigma\bigl(\li\bigr)-1\nonumber
\end{array}
$$
\elm

{\em Proof :} The relations for the ${\cal O}_1\,$- and cancellation moves
follow from the multiplicativity of $\{.\}\,$ and the additivity of 
$\sigma$ and $\chi\,$. E.g., we have $\{\li_1\cup\li_2\}\,$ if  $\li_1$ and
$\li_2$ are two disjoint links that can be separated by a 2-sphere in
$S^3\,$. For the verification of $\{\clr\}=1\,$ we remark that $1\in\cal
I$ is the only  coloration of the cancelling ribbon that contribues to
$F(\clr,C)\,$ and $d(1)e^1\circ f^1=1\,$ by construction.
\medskip

In order to check the modification move we decompose $id\in
End(l_1\otimes\ldots\otimes l_n)\,$ into sums of composites of 
projections onto and injections of irreducibles $k\in\cal I\,$.
They are presented by pairs of coupons as in the cancellation configuration
with an additional ribbon $R$ of color $C(R)=k$ joining them.
We compose this presentation of $id$ with the modification configuration on
the left side of Figure~\ref{fig-rel-17}, and  slide the modification
annulus  $A$ between the coupons  
over $R$. If we sum now over the colorations $C(A)$ we easily find from 
equation~(\ref{eq-d-hat}) that we are left with only $C(R)=1\,$ as a 
possible channel and an extra factor $\frac 1 {\hat d(1)}\,$.
\hfill$\Box$
\s

It is obvious from Proposition~\ref{pp-BL-moves} and
Lemma~\ref{lm-trsf} that there is only one way to construct an invariant
from the given data that is multiplicative with respect to connected
summing and for which $\,\tau(S^3)=1$ :

\begin{cor}
If ${\cal D}^2=\hat d(1)^{-1}\,$ then the expession 
\beq\label{invar}
\tau(M)\,=\,{\cal
D}^{-\chi(\li)-\sigma(\li)}\{\frL\}^{\sigma(\li)}\bigl\{\li\bigr\} 
\eeq
only depends on the three manifold that is presented by the bridged link
$\li\,$.
\end{cor}

Note that the Corollary implies the identity $\{\frL\}\{\FLO\}={\cal
D}^{-2}\,$, which allows us to substitute one of the coefficient with an
expression in $\{\FLO\}\,$. As outlined in the end of this section it is not
hard to extend  $\tau$  to a Topological Field Theory for any BTC satisfying
(\ref{eq-d-hat}), using the same arguments as above.
In other words, we can construct a 
$\otimes$-functor $\tau\,:\,\overline{c\!\widetilde{{\cal T}\!\!g\,}}\,\to
\,Vect({\bf C})\,$, which is specializes to the invariant for cobordisms
between empty surfaces. The vector spaces that is associated, e.g., to the
torus $T\,=\,S^1\times S^1\,$ is canonically identified with
$\,\tau(T)\,=\,Hom(F,1)\,\cong\,{\bf C}^{\cal I}\,$, where
$F\,:=\,\bigoplus_{j\in\cal  I}j^{\vee}\otimes j\,$. For the 
cobordism from $T$ to itself, represented by the tangle $S^+\,$ as defined
in (\ref{eq-cl}), we can give the map $\tau(S^+)\,$ in terms of the matrix
$S_{i,j}\,$ from equation (\ref{def-S}). Clearly, $S^+$ is invertible in 
$\overline{c\!\widetilde{{\cal T}\!\!g\,}}\,$ with inverse $S^-\,$ and
$\tau(S^-)\,=\,\bigl(S_{i,j^{\vee}}\bigr)\,$ for a suitably scaled basis. This implies
a purely algebraic statement, namely that a category which satisfies
(\ref{eq-d-hat}) also has the `stronger' modularity property. This
implication can also be proven formally, without reference to the topological
situation:

\blm\label{lm-mod-eq}
For a semisimple BTC the following are equivalent:
\ben
\item \beq\label{mod-new} 1\in im(S)\eeq

\item \beq\label{xxx}\sum_{j\in{\cal I}}dim(j)S_{j,i}\,=\,{\cal
D}^2\delta_{i,1}\eeq\label{eq-Sd=1} 

\item $$\mbox{$S$ is invertible.}$$
 
\een
\elm

{\em Proof :} We prove this lemma by showing that {\em 1)} implies {\em
3)}, and that {\em 3)} implies {\em 2)}. 
\ben
\item [$ 1)\Rightarrow 3)$:] Let $\hat d(j)$ be as in (\ref{eq-d-hat}). The ``Verlinde-formula''

$$
dim(k)^{-1}S_{i,k}S_{j,k}\,=\,\frac {S_{i\otimes j,k}}{S_{1,k}}\,=\,\sum_pN_{ij,p}S_{p,k}\,,
$$
is a simple consequence of the general cyclicity of $tr_X\,$. ( $dim(k)\neq
0\,$ is part of the semisimplicity condition). 
We multiply this with $\hat d(k)$ and sum over $k\,$. Using the symmetry
of $S$ we arrive at the equation:
\beq\label{sysc}
SYS\,=\,C\;,
\eeq
where $C$ is the conjugation matrix and $Y$ is the  diagonal matrix, whose
entries are the numbers $\hat d(k)dim(k)^{-1}\,$. Thus $S$ is invertible. Note also that
invertibility of  $Y\,$ implies $\hat d(1)\neq 0\,$.
\item [$ 3)\Rightarrow 2)$:] We have {\em 1} and can use again
equation~(\ref{sysc}) and 
$dim(j)=dim(j^{\vee})$ to show:
$$
\sum_j dim(j)
S_{j,i}\,=\,\bigl(SCS\bigr)_{1,i}\,=\,(Y^{-1})_{1,i}\,=\delta_{1,i}\hat
d(1)^{-1}\;.$$
\hfill$\Box$ 
\een
\medskip

The fact that $S\,$ is invertible and that the vector $d(j)\,$ as defined in
(\ref{eq-def-[]}) is also a solution of (\ref{eq-Sd=1}) implies
\beq\label{eq-d=d}
d(j)\,{\cal D}^2\hat d(j)\,=\,dim(j)
\eeq

Moreover, if we specialize (\ref{eq-Sd=1}) to $\,i=1\,$ we obtain the
expression
\beq\label{eq-def-D}
{\cal D}^2\,=\,\sum_{j\in{\cal I}}dim(j)^2
\eeq
as it was originally defined in [Tu]. It also follows immediately that in
the notation of [Tu],
\beq\label{eq-def-delta}
\Delta\,=\,\{\frL\}\,=\,\sum_{j\in{\cal I}}\theta_j^{-1}dim(j)^2\;.
\eeq
Thus the invariant defined in (\ref{invar}) coincides with the construction from [Tu],
where the input data have been the quantum dimensions $dim(j)\,$ and the statistical phases
$\theta_j\,$.  
\medskip

In conclusion, let us note that our proof of invariance avoided the use of the Kirby
$\kappa$-move, and the result from [FR] that it generates all two-handle
slides of links.  Other than in [RT] the 
verification of invariance is straight forward and does not involve any
algebraic computation, if we start from the quatities $\hat d(j)\,$ and
$\{\frL\}\,$. The computations needed to compare it to the original
definition are also much simpler than in [RT] where most parts of the 
$SL(2,{\bf Z})$-representation intrinsic to a BTC hace to be constructed
by algebraic computation rather than topological arguments. Thus the
bridged link presentations and the condition (\ref{eq-d-hat}) are  a much more
natural and convenient starting point for the construction of three manifold invariants
than the Kirby-Fenn-Rourke-Calculus and the modularity condition of Turaev.
$\smile$
\s

The RT-construction of invariants can be generalized to non-semisimple
BTC's, as for instance the representation categories  $mod-H\,$ of a
quantum group $H\,$ at
roots of unity. For finite dimensional, quasitriangular Hopf algebras $H\,$
the first definition was given by M. Hennings, and, independently,  the
generalization to arbitrary BTC's with enough limits by V. Lyubashenko. Their
properties and rules of computation have been studied by L. Kauffman, T.
Ohtsuki, and D. Radford, see [HKLR]. Let us briefly outline the  algorithm 
by which the invariant for  $mod-H\,$ can be  computed:
\medskip

Along a component of a link elements of $H\,$ are inserted and moved as
follows: If ${\cal R}\in H^{\otimes 2}\,$ is the R-matrix 
we insert the elements ${\cal R}^{(1)}_j\,$ in one component of an overcrossing 
and the elements ${\cal R}^{(2)}_j\,$ in the other, and replace the
overcrossing by a singular crossing which does not distinguish over and
under. The elements can be moved through extrema using the antipode and 
successive elements along a component can be multiplied. It is clear that we
reduce a plat link to a link with only singular crossings and  only one
element (which may depend on a summation index) at a prescribed spot 
for each component. The invariant is obtained by evaluating the 
integral $\mu\in H^*\,$ on each of these elements and sum over all
products. 
\medskip

An integral of a Hopf-algebra is (up to normalization) uniquely defined
by the property $\bigl(\mu\otimes 1)\Delta(y)= 1\mu(y)\,$. The topological
translation of this condition is the ${\cal O}_2\,$-move, since doubling a 
piece of a link component corresponds in the above algorithm to taking the 
coproduct of all of the elements along this piece. In the bridged link
formalism we have a priori no requirement on $\mu\in H^*\,$, since the
two handle slides are not among the selected set of generating moves. Still, we 
need to extend the algorithm  to 1-handle
attachements, such that the computed invariant is consistent with the modification move. This
dictates the following prescription: For each pair of coupons we connect
the $\,n\,$ incoming strands to the ones going out at the other coupon
along any path (with singular crossings) and insert the elements
$\Lambda_j^{(k)}\,$ in the $k$-the component, where
$\Delta^{n-1}(\Lambda)\,=\,\sum_j\Lambda_j^{(1)}\otimes\ldots\otimes\Lambda_j^{(n)}$
and 
\beq\label{eq-int-coint}
\Lambda\;=\;\bigl(\mu\otimes 1\bigr)\Bigl({\cal R}^t{\cal R}\Bigr)\qquad\in\,H\;.
\eeq
Now, the move $Rel_{12}\,$ from Figure~\ref{fig-rel-12} imposes the
cointegral constraint $y\Lambda\,=\,\Lambda y\,=\,\epsilon(y)\Lambda\,$,
where $y$ can be any element in $H\,$, if $H\,$ is nicely quasitriangular,
as, e.g., for doubles. Also, for doubles it is easy to infer from
(\ref{eq-int-coint}) that $\mu\,$ has to be a right integral whenever $\Lambda$
is a cointegral. The cancellation move imposes the normalization
$\mu(\Lambda)=1\,$. In fact, it is a well known result from Hopf algebra
theory that $\mu(\Lambda)\neq 0\,$ no matter if $H\,$ is semisimple or not.
Finally it is a fact that for doubles $S(\Lambda)=\Lambda\,$, as required
by move $Rel_{13}\,$ of Figure~\ref{fig-rel-13}.
\medskip

For a semisimple BTC the categorial integral of the braided Hopf algebra
$H:=F^{\vee}\,$ is $\mu\,=\,\sum_j
dim(j)F\Bigl(\ _{j^{\vee}}\ccap\,_j\Bigr)\;\in\,Hom(1,F)\,$ and the
cointegral
$\Lambda\,$ is the projection onto invariance. 

The fact that
$\Lambda^2=0\,$ for non-semisimple BTC's (or equivalently that \mbox{$\tau(S^1\times
S^2)=0\,$}) makes it impossible to extend the invariants to TQFT's that
observe the $\otimes$-rule.
\s  
  
We conclude this section with a 
table of notions and conditions that we have found to be related:
\s

\s

\begin{tabular}{l|r}\hline
SURGERY & BTC / HOPF ALGEBRA\\ \hline\hline
One Handle & Projection on Invariance / Cointegral \\ \hline
Two Handle & Canonical Trace / Integral \\ \hline
Framing & Balancing/ Ribbon Graph Element\\ \hline
Isotopies & Braid and Rigidity Relations / Quasitriangular\\ \hline
0-1-Cancellation & $\chi(\li)$-Normalization with $\cal D$ or ``Quantum-rank''\\ \hline
1-2-Cancellation & Non-degenerate Pairing of Invariance\\
&and Coinvariance ($\Leftrightarrow$ Semisimplicity) or \\ 
& Contraction of Integral and Cointegral\\ \hline
${\cal O}_1$ or $\#{\bf CP}^2$-move& $\sigma(\li)$-Normalization with
$\Delta$ or Moduli\\ \hline 
Modification- , $\eta$-, or $\,\cup \hh^2_{(5)}$-move &  $1\in im(S)$ / Double\\ \hline
\end{tabular}
\s


{\em 4.2.2) The Connecting Annuli and Selfconjugate Objects}\label{pg-422}
\medskip

In the case of topological surgery presentations we have to be careful
about replacing a piece of a link as in Figure~\ref{fig-conn-lm} by
straight strands as in the middle picture. However, for the evaluation of the
RT-invariant this modification can be made, if we properly modify the
summation $\sum_C....\prod_{L}dim(C(L))$ in the definition (\ref{invar}) of
$\tau\,$.
\s

In the simplest, non-trivial  case of Lemma~\ref{lem-con} with a link $\cal
L$, where $A^+$ is 
connected to $B^+$, we have only one link component $L_o$ going through the
annulus. If we applied the na\"\i ve move as for $A^+...A^-$ we would get 
a link ${\cal L}'\,$, where $L_o$ is replaced by  two components $L^+$ and
$L^-\,$. The total number of components stays the same. 
For a given coloration with $C(L_o)=j$ we may replace the annulus
by a pair of coupons, for which we insert dual bases of the (co-)invariance
of $j\otimes j^{\vee}\,$. A possible choice is, such that the canonical
morphisms (associated to the maxima and minima of a ribbon graph) are
$f:1\to j\otimes j^{\vee}$ and $dim(j)e:  j\otimes j^{\vee}\to 
1\,$.  

Thus we can compute the invariant from the link ${\cal L}'$ by confining
the summation to colorations with $C(L^+)=C(L^-)\,$ in the relative
orientation induced by $L_o$, and omitting the $dim(j)^2$-contribution to
the product $\prod_{L}dim(C(L))\,$ coming from $L^{\pm}$.
\medskip

The more interesting case is given by the situation $A^+...B^-$. Here, we
have again
only one link component $L_o$ passing through the annulus, but this time
with linking number 2 instead of 0. This entails that the diagram is zero
for colorations for which $j=C(L_o)\,$ is not selfconjugate.
\medskip

For a more precise statement we observe that for a selfconjugate,
irreducible object 
$k$ the map $\rho_k$ as defined in (\ref{tep}) is an involution on
$Hom(1, k\otimes k)\cong{\bf C}\,$, i.e., $\rho_k=\pm 1\,$.

The invariant may now be computed from ${\cal L}'$ by substituting the
$dim(j)$ contribution in the product $\prod_{L}dim(C(L))\,$ by ${\cal
D}^{-1}\rho_j\,$, and 
confining the summation to colorations with $C(L_o)=C(L_o)^{\vee}\,$.

In fact the tangle in Figure~\ref{fig-conn-lm} defines a natural
transformation of the identity of the BTC, with endomorphisms
$\xi(X):X\to X\,$ uniquely determined by its values on the simple objects:
\beq\label{nat-sc}
\xi(j)=\left\{ \begin{array}{ll}
{\cal D}^{-1}dim(j)^{-1}\rho_j & \mbox{if $j=j^{\vee}$}\\
0& \mbox{elsewise}
\end{array}\right.\;.
\eeq
Hence we may also look at the evaluation of ${\cal L}'$ with the morphism
$\xi$ inserted along $L_o\,$.
\medskip

The natural transformation $\xi$ in the representation category of a
quantum double $H\,$ is given by the remarkable, central element,
\beq\label{R}
\rho\,=\,\sum_j\Lambda'f_j\hat u\Lambda''e_j\;,
\eeq
which projects onto the selconjugate subrepresentation when it is applied
to a general representation of $H\,$. $\smile\;$
 
In (\ref{R}) $\Lambda$ is the cointegral of the double as discussed in the previous subsection.
For the form of $p$ and our  conventions for coproduct, dual bases
$\{e_j\}$ and $\{f_j\}\,$, and $\hat u$ see, e.g., [Ke1].

\s

The constructions of Sections 4.2.1) and 4.2.2) are not confined to 
closed manifolds but apply to cobordisms as well. Let us conclude this
section with a summary of 
the construction of an anomalous TQFT from the tangle
presentation developed in Chapter~3:
\medskip

To a (connected) cobordism $M\,$ in $\wct\,$ we associate a corresponding
admissible tangle,~$t\,$, and to  this the special, closed tangle
$t_M:=\cc(t)\,$ as in (\ref{eq-cl}). For a given coloration we obtain as in [RT0] 
a morphism:
\medskip

$$
F(t_M,C)\,:\,j_1^{\vee}\otimes j_1\otimes\ldots\otimes j_{\tilde g^-}\,\to\,
\,k_1^{\vee}\otimes k_1\otimes\ldots\otimes k_{\tilde g^+}
$$

Summing over colorations of  the internal, closed ribbons we obtain, 
after application of $Hom(-,1)\,$, the linear map:
$$
\{t,\ub j,\ub k\}:Hom\bigl(j_1^{\vee}\otimes j_1\otimes\ldots\otimes j_{\tilde
g^-},\,1\bigr)\,\to\, 
\,Hom\bigl(k_1^{\vee}\otimes k_1\otimes\ldots\otimes k_{\tilde g^+},\,1\bigr)
$$
with
$$
\{t,\ub j,\ub k\}\:=\sum_{\mbox{\tiny internal} C}\prod _{\mbox{\tiny internal}
l}dim(l)\,F(t_M,C)\;.
$$

Here $\ub j$ and $\ub k$ are short hand for the $\tilde g^{\pm}\,$ colors
at  the boundaries $X_{`\pm}\,$. We introduce the canonical injections and 
projections between the total invariance and the product of the invariances
of the groups $G_s\,$:
$$
Hom(j_1^{\vee}\otimes j_1\otimes \ldots\otimes j_{\tilde
g},\,1)\quad\lower 2ex \hbox{$\stackrel{\stackrel {p}{\longrightarrow}}{\stackrel
{\longleftarrow}{i}}$}\quad
\bigotimes_{s=1}^{K}Hom\bigl(j_{s,1}^{\vee}\otimes j_{s,1}\otimes\ldots\otimes
j_{s,g_s},\,1\bigr)\;, 
$$
if we have $K$ components of genera $\{g_s\}$.
Also, we introduce the ``Euler number'' of the tangle:
$$
\chi^+(t)\,=\, (\mbox{\rm\small\# of internal components of $\cc(t)$})\,\,+\,\,\tilde g^+\,\,-\,\,K^+\,\,
+2\,,
$$
where $\tilde g^+$ and $K^+$ are the total number of pairs of ribbons and
the number of components at $X_+\,$, respectively.

For the  object $F:=\oplus_{j\in\cal I}j^{\vee}\otimes j\,$ we have a
canonical decomposition:
$$
\tau\bigl((g)\bigr)\,:=\,Hom(F^{\otimes g},1)\,=\,\bigoplus_{\ub j}Hom(j_1^{\vee}\otimes j_1\otimes
\ldots\otimes j_g,\,1\bigr)
$$
On this we define the linear map associated to the cobordism by the TQFT
$$
\tilde \tau(M)\,:\,
\tau(g^-)\,=\,\bigotimes_{s=1}^{K^-}\tau\bigl((g^-_s)\bigr)\;\longrightarrow\;\tau(g^+)\,=\,\bigotimes_{s=1}^{K^+}\tau\bigl((g^+_s)\bigr)
$$
by the sum of the bloc entries:
$$
\tilde \tau(M)\,=\,{\cal D}^{-\chi^+(t)}\,\bigoplus_{\{\ub j\}\{\ub
k\}}\,\prod_{\,j\in\ub j}dim(j)\,p_+\circ\{t,\ub j,\ub
k\}\circ i_- 
$$
Note, that we consider only the ribbons at $X_+$ in the additional product
of $q$-dimensions. Also, the signature, $\sigma(\li)\,$,
does not appear in the formula, since we only want to have an anomalous functor 
$\tilde \tau\,$ defined on $\wct$. Compatibility of $\tilde \tau\,$ with
the composition rules of Section 3.2.3) follow from (\ref{xxx}). 
Since we projected on the invariance of the individual groups the
$\tau$-move is given by a pure braid with one $1$-  colored strand and
therefore does not change the morphism. Invariance of the functor under 
all other moves follows easily from the results for closed manifolds.


\subparagraph{4.3) Punctured Cobordisms and Glue-$\otimes$}\label{pg-43}
\
\s

The definition of $\ct$ has a natural generalization, which includes 
punctured surfaces. To be more precise we define for any $N\in {\bf Z}^+\,$ 
a cobordism category $\ct(N)\,$ as follows:
\s

The objects are compact, oriented two folds $\Sigma$ with orientation
preserving homeomorphisms 
$$
\zeta\,:\, S(N):=\amalg_{i=1}^{N}S^1\,\isto\,\partial\Sigma\;,
$$
parametrizing the boundaries by a oriented, standard manifold $S(N)\,$. A morphism 
from a surface $\Sigma_1$ to $\Sigma_2$ is given by a three fold $M\,$,
and coordinate maps from its boundary to the composed surface:
$$
\psi\,:\,\partial M\,\isto\,-\Sigma_1\amalg_{\zeta_1}\bigl\{
S(N)\times[0,1]\bigr\}\amalg_{\zeta_2}\Sigma_2\;.
$$
We also impose a similar notion of equivalence as for the closed case.
The composition of morphisms is then given as usual by  gluing the two
three folds together along a common boundary piece $\Sigma_2\,$. We obtain a
three  fold with boundary 
$$
-\Sigma_1\amalg_{\zeta_1}\bigl\{
S(N)\times[0,1]\bigr\}\amalg_{\zeta_2^{-1}\circ\zeta_2'}\bigl\{
S(N)\times[0,1]\bigr\}\amalg_{\zeta_3'}\Sigma_3
$$
For a suitable redefinition of the $\zeta$- coordinate maps we may replace
the two middle parts by one $S(N)\times [0,1]\,$.
\s

The outlined axioms entail representations 
$$\pi_0\bigl(\diff(\Sigma,\partial\Sigma
)^+\bigr)\,\cong\,Aut_{\ct(N)}\bigl(\Sigma\bigr)\,$$  
 and
$$
\pi_1\bigl(\partial
\Sigma\bigr)\,\hookrightarrow\,Nat_{\ct(N)}\bigl(id,id\bigr)\;.
$$  
Also, $\ct(1)$ has a
natural structure of a braided tensor category and contains a canonical,
braided Hopf algebra (see, e.g., [Ke1] and [Ke2]). 
\s

 In [KL] we develop the analogous tangle presentation of
 $\ct(N)\,$. It is obtained from
the one of 
$\ct=\ct(0)$ by considering a surgery presentation of the corresponding
cobordism   in $\ct$  where the cylinders  $S^1\times [0,1]\,$, are filled with 
tubes $D^2\times [0,1]\,$.
The latter are tubular neighborhoods of $N$  strands joining the opposite 
surfaces of the cobordism, and we may assume that they are disjoint from
the other surgery ribbons. A standard presentation of a connected element of
$\ct(N)\,$ is thus given by a standard bridged link diagram in $\H\,$
with $N$ additional ribbons that start and end in  opposite boundary
components  of the $\ct$-cobordism, and have a prescribed standard form
inside the handlebodies $H^{\pm}_g\,$.
The moves of the presentation are obtained by treating the additional
ribbons like singularities with highest values., i.e., they may go through
surgery spheres and can be slid over other two handles. 
\s

A far more interesting aspect of this family of cobordism categories is
another type of ``glue operation'', which may be understood as a
second independent composition among relative cobordisms. We start in the
definition with the choice of some 
orientation reversing
involution $\rho':S^1\isto S^1\,$. From the manifolds $S(N)$ and $S(M)$
we select $K$ components and construct a functor
$$
\ct(N)\times \ct(M)\,\longrightarrow\,\ct(M+N-2K)
$$
as follows:
\s

For two surfaces $\Sigma\in Ob\bigl(\ct(N)\bigr)$ and $\Sigma'\in
Ob\bigr(\ct(M)\bigl)\,$ the product is defined as the sewed surface,
$$
\Sigma\amalg_{\sim}\Sigma'
$$
where we glue respective components $C\subset\Sigma$ and
$C'\subset\Sigma'\,$
of the boundaries together using the 
identifications
$$
\zeta'\circ\rho\circ\zeta^{-1}\,:\,C\isto C'\;.
$$
Moreover, we define the product of cobordisms $M$ and $M'$ in these
categories as a quotient space:
$$
M\amalg_{\sim}M'\quad,
$$
where we use the identification along the cylindrical boundary pieces 
$S^1\times[0,1]\cong T\subset \partial M\,$ given by
$$
\psi'\circ(\rho\times id)\circ\psi^{-1}\,:\,T\isto T'\;.
$$
Clearly, the definitions for objects and morphisms are compatible.
\s

In [KL] we describe a very natural way of organizing these two type
gluings
over surfaces in terms of {\em double categories}, whose ingredients
are always two types of compositions.  Specifically, the pastings over
surfaces are viewed as {\em vertical compositions}, and the gluings
over the cylindrical pieces as {\em horizontal compositions}.
\s

\begin{figure}[ht]
\begin{center}
\epsfxsize=5in
{\ \epsfbox{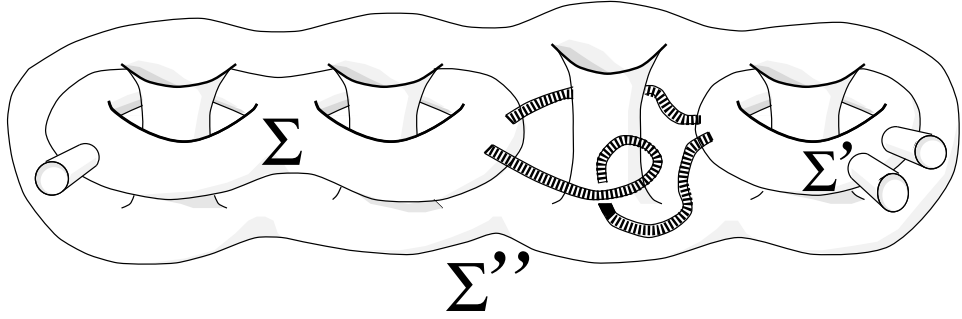}}
\end{center}

\caption{Non-canonical: $\aleph\,:\,\Sigma\otimes_{glue}\Sigma'\,\to\,
\Sigma''\,$}\label{fig-ex-aleph}

\end{figure}

To describe the (horizontal) glue product in terms of standard presentations,
we observe that $\Sigma\otimes_{glue}\Sigma'$ is homeomorphic to a standard 
surface $\Sigma''$ of genus $g''=g+g'+K-1\,$, with $N+M-2K$ holes.
This allows us to define a  cobordism $\aleph\,:\,\cong
\Sigma''\times[0,1]\,$, which may be  seen as a morphism between objects of
the different categories $\ct(N)\times\ct(M)$ and $\ct(N+M-2K)\,$:
$$
\aleph\,:\,\Sigma\amalg_{\sim}\Sigma' \,``\isto\mbox{''}\,  \Sigma''\;.
$$
We present $\aleph$  in a way analogous to the generalized standard
presentation, which we outlined above. In this picture we start from a cobordism
$$
\hat \aleph\,:\,\hat\Sigma\amalg\hat\Sigma'\,\to\,\hat\Sigma''\;,
$$
where $\hat \Sigma^*$ denotes the corresponding, closed manifold. Starting
and ending at the original punctures we include ribbons: $N-K$ of them going from
$\hat \Sigma$ to $\hat \Sigma''\,$, $M-K$ going from $\hat \Sigma'$ to $\hat
\Sigma''$, and $K$ going from $\hat \Sigma$ to $\hat \Sigma'$. An example
for $g=2,\,g'=1,N=3,\,M=4,\,K=2\,$ is given in Figure~\ref{fig-ex-aleph}.

For  two cobordisms $M:\Sigma_1\to\Sigma_2\,$ and
$M':\Sigma_1'\to\Sigma_2'\,$ let us consider the composite:
$$
M\tilde\otimes M'\,:=\,\aleph_2\circ\bigl(M\otimes M'\bigr)\circ\aleph_1^{-1}\,:\,\Sigma_1''\to\Sigma_2''\;.
$$

Let us also briefly describe the basic topological transformation
 used to extract  from the glue operation a compatible composition 
law for the tangles:
\s

\begin{figure}[ht]
\begin{center}
\epsfxsize=4.6in
{\ \epsfbox{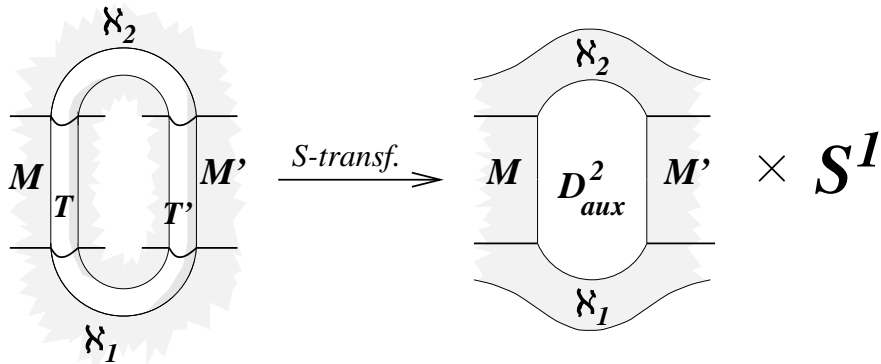}}
\end{center}

\caption{Glue-$\otimes$ Surgery}
\label{fig-repl}
\end{figure}

A  pair $(T,T')$ of the $K$ cylindrical pieces in the boundaries of $M$
and $M'$ are now joined by the ribbons in $\aleph_i$'s. Thus
$M\tilde\otimes M'\,$ is a morphism in $\ct(N+M-2K)\,$, with $K$ solid
tori removed from the inside. A vicinity of such a torus is depicted on the
left of Figure~\ref{fig-repl}. 
\s

Using an `$S$-transformation', we find it to
be homeomorphic to the region on the right side of Figure~\ref{fig-repl}
times a circle. If we fill in a $\,D^2_{aux}\times S^1\,$, we obtain the glued
tensor $M\otimes_{glue}M'\,$ with an ``irrelevant'' 
$\,D^2_{aux}\times S^1\,\cong\,[0,1]\times
T\,$  inserted between $T$ and $T'\,$.
\s

Thus if we simply {\em reinterpret} the $K$ closed puncture ribbons as
surgery
ribbons we obtain a presentation of $M\otimes_{glue}M'\,$. Using similar 
moves as for the case of closed surfaces this can be brought again into a
standard form.
\s

The fact that the family of $\aleph\,$'s for all pairs of standard 
surfaces is by no means canonical, which is the source of some 
problems:
\s

For example one has to verify that the $\aleph\,$'s are chosen, such
that composites for different orders of sewing the same surface
together yield equivalent cobordisms between the different pieces and
the glued standard surface. 
\s

More importantly, in the construction of extended TQFT's for surfaces
with boundaries, we can at most expect to be able to construct 
{\em pseudo functors} from the double category of relative cobordisms 
into a corresponding algebraic double category. This functor will
respect the $\,\otimes_{glue}\,$-product only up to an equivalence
depending on the choice of $\aleph\,$'s. 

Specifically, in [KL] we consider the algebraic  double category,
in which the 0-Objects are $n$-fold tensor products, 
$\,{\cal C}\odot\ldots\odot{\cal C}\,$,  of an abelian category
$\,{\cal C}\,$, the horizontal and vertical 1-arrows are functors,
and the 2-arrows are natural transformations between functors. 
We can, however, obtain an honest functor between the two double
categories,
by admitting not one but a finite, combinatorial set of surfaces
as 1-arrows for each homeomorphism class (characterized by number of
components, genera, and holes).

\nopagebreak

%% file: BLbibl.tex
\section*{References}
\ben 

\item[[At]] Atiyah, M.: On Framings of 3-Manifolds. {\em Topology}, Vol. {\bf
29} (1990) 1-7.

\item[[Bi]] Birman, J.: Braids, Links and Mapping Class Groups. {\em Princeton Univ.
Press} (1974).

\item[[Ce]] Cerf, J.: La Stratification Naturelle des Espace de Fonctions
Diff\'erentiables r\'eeles et la Th\'eor\`eme de la Pseudoisotopie.
{\em Publ. Math I.H.\'E.S.} {\bf 39} (1970).

\item[[FR]] Fenn, R., Rourke, C.: On Kirby's Calculus of Links. {\em Topology}, Vol.
{\bf 18} (1979) 1-15.      
 
\item[[HW]] Hatcher, A., Wagoner, J.: Pseudo-Isotopies of Compact Manifolds.
{\em Ast\'erisque} {\bf 6}. (1973).

\item[[HKRL]] Hennings,M.A.: Invariants of Links and 3-Manifolds Obtained
from Hopf Algebras, Cambridge Preprint, 1990. 

Kauffman, L., Radford, D.:
Invariants of 3-Manifolds derived from Finite Dimensional Hopf Algerbras,
Chicago Preprint, 1994. 

Lyubashenko, V.: Invariants of 3-Manifolds and 
Projective Representations of Mapping Class Groups via Quantum Groups at
Roots of Unity,  {\em Commun. Math. Phys} {\bf 172} (1995) 467-516.

\item[[Ig]] Igusa, K.: The Space of Framed Functions. {\em Trans. \hspace*{-.3cm}
Amer. \hspace*{-.3cm} Math. \hspace*{-.3cm} Soc.} Vol. {\bf 301-2}  (1987) 431-477.

\item[[J]] J\"anich, K.: Charakterisierung der Signatur von Mannigfaltigkeiten
durch eine Additivit\"atseigenschaft. {\em Invent. Math.} {\bf 6}(1968), 35-40.

\item[[Ke1]]Kerler,T.: Mapping Class Group Actions on Quantum Doubles. 
 {\em Commun. Math. Phys.} {\bf 168} (1995) 353-388.

\item[[Ke2]]Kerler, T.: A Topological Hopf Algebra.(Talk presented at UNC,
April 1994, paper in preparation)

\item[[Ke3]]Kerler, T.: Equivalence of a Bridged Link
Calculus and Kirby's Calculus of Links on Non-Simply Connected 3-Manifolds.
{\em Topology Appl.} to appear.\newline
Available at {\tt http://www.math.ohio-state.edu/$\sim$kerler/papers/KC/}

\item[[Ki]] Kirby, R.: A Calculus for framed Links in $S^3\,$, {\em Invent. Math.} {\bf
45} (1978), 35-56.

\item[[KL]] Kerler, T., Lyubashenko,
V.: Non-semisimple Topological Quantum Field Theories for Connected
Surfaces. Preprint.

\item[[Li]] Lickorish, W.B.R.: A Representation of Orientable 3-Manifolds.
 {\em Ann. Math.} {\bf 76} (1962) 531-540.

\item[[Mi]] Milnor, J.: Morse Theory. Annals of Mathematical Studies.
{\em Princeton University Press} {\bf 51} (1969).

\item[[MP]] Matveev, S.,  Polyak, M.: A Geometrical Presentation of the Surface 
Mapping Class Group and Surgery. {\em Commun. Math. Phys.} {\bf 160} (1994) 537.
 
\item[[RM]] Roberts, J.D., Masbaum, G.: On Central Extensions of Mapping
Class Groups. Mathematica Gottingensis, {\bf 42} (1993). 

\item[[Ro]] Rohlin, V.A.: A Three Dimensional Manifold is the Boundary of a
Four Dimensional One, {\em Dokl.Akad. Nauk.} SSSR (N.S.) {\bf 114} (1951)
355-357.

\item[[Rf]] Rolfsen, D.: Knots and Links. Mathematics Lecture Series {\bf
7}. {\em Publish or Perish} (1976).

\item[[RT0]] Reshetikhin,N., Turaev,V.: Ribbon Graphs and Their Invariants
Derived from Quantum Groups. {\em Commun. Math. Phys.} {\bf 127}(1990), 1-26. 

\item[[RT]] Reshetikhin,N., Turaev,V.: Invariants of 3-Manifolds via Link
Polynomials and Quantum Groups. {\em Invent. Math.} {\bf 103} (1991) 547-598.

\item [[Sw]] Sawin, S.: Extending the Kirby Calculus to Manifolds with
Boundary, MIT-preprint, 1994. 

\item[[Tu]] Turaev, V.: Quantum Invariants of 3-Manifolds. {\em Walter de Gryter},
 Berlin (1994) (to appear).

\item[[Wj]] Wajnryb,B.:  A Simple Presentation for the Mapping Class Group of an
Orientable Surface. {\em Israel J. Math.} {\bf 45} (1982) 157-174. 

\item[[Wa]] Wall, C.T.C.: Non-Additivity of the Signature, {\em Invent. Math.} {\bf 7}
(1969), 269-274.

\item[[Wa2]] Wall, C.T.C.: Determination of the Cobordism Ring, {\em Ann. Math.}
{\bf 72} (1960), 292-311.

\item[[Wc]] Wallace, A.H.: Modifications and Cobounding Manifolds, {\em Cam. J. Math.}
{\bf 12 } (1960), 503-528.

\item[[Wk]] Walker, K.: On Witten's 3-Manifold Invariants. UCSD-Preprint (1991).

\een

\bigskip

{\sc Harvard University, Cambridge, MA, USA}

present address:

{\sc The Ohio State University, Columbus, OH, USA}

{\tt kerler@math.ohio-state.edu}